\documentclass[12pt]{amsart}
\usepackage{amssymb}
\usepackage{graphicx}

\newtheorem{thm}{Theorem}[section]
\newtheorem{lem}{Lemma}[section]
\newtheorem{cor}{Corollary}[section]

\newtheorem{prop}{Proposition}[section]

\theoremstyle{remark}
\newtheorem{rem}{Remark}[section]
\newtheorem{claim}{Claim}[section]

\theoremstyle{definition}
\newtheorem{defn}{Definition}[section]

\numberwithin{equation}{section}

\allowdisplaybreaks[4]

\begin{document}

\newcommand{\claimref}[1]{Claim \ref{#1}}
\newcommand{\thmref}[1]{Theorem \ref{#1}}
\newcommand{\propref}[1]{Proposition \ref{#1}}
\newcommand{\lemref}[1]{Lemma \ref{#1}}
\newcommand{\coref}[1]{Corollary \ref{#1}}
\newcommand{\remref}[1]{Remark \ref{#1}}
\newcommand{\conjref}[1]{Conjecture \ref{#1}}
\newcommand{\secref}[1]{Sec. \ref{#1}}
\newcommand{\ssecref}[1]{\ref{#1}}
\newcommand{\sssecref}[1]{\ref{#1}}
\newcommand{\figref}[1]{Figure \ref{#1}}

\def \d#1{\displaystyle{#1}}
\def \mult{\mathop{\mathrm{mult}}\nolimits}
\def \rank{\mathop{\mathrm{rank}}\nolimits}
\def \codim{\mathop{\mathrm{codim}}\nolimits}
\def \Ord{\mathop{\mathrm{Ord}}\nolimits}
\def \Var{\mathop{\mathrm{Var}}\nolimits}
\def \Ext{\mathop{\mathrm{Ext}}\nolimits}
\def \ext{\mathop{{{\mathcal E}xt}}\nolimits} 
\def \Pic{\mathop{\mathrm{Pic}}\nolimits}
\def \Spec{\mathop{\mathrm{Spec}}\nolimits}
\def \dsum{\mathop{\oplus}}
\def \mapright#1{\smash{\mathop{\longrightarrow}\limits^{#1}}}
\def \mapleft#1{\smash{\mathop{\longleftarrow}\limits^{#1}}}
\def \mapdown#1{\Big\downarrow\rlap{$\vcenter{\hbox{$\scriptstyle#1$}}$}}
\def \smapdown#1{\downarrow\rlap{$\vcenter{\hbox{$\scriptstyle#1$}}$}}
\def \A{{\mathbb A}}
\def \I{{\mathcal I}}
\def \J{{\mathcal J}}
\def \CO{{\mathcal O}}
\def \C{{\mathcal C}}
\def \W{{\mathcal W}}
\def \BC{{\mathbb C}}
\def \m{{\mathcal M}}
\def \n{{\mathcal N}}
\def \H{{\mathcal H}}
\def \S{{\mathcal S}}
\def \Z{{\mathcal Z}}
\def \BZ{{\mathbb Z}}
\def \X{{\mathcal X}}
\def \Y{{\mathcal Y}}
\def \T{{\mathcal T}}
\def \P{{\mathbb P}}
\def \PGL{{\mathbb P}{\mathrm{GL}}}
\def \G{{\mathbb G}}
\def \F{{\mathbb F}}
\def \BR{{\mathbb R}}
\def \CR{{\mathcal R}}
\def \Hom{{\mathrm{Hom}}}
\def \hom{\mathop{{{\mathcal H}om}}\nolimits}
\def \nilrad{{\mathrm{nilrad}}}
\def \Supp{{\mathrm{Supp}}}
\def \Ann{{\mathrm{Ann}}}
\def \closure#1{\overline{#1}}
\def \EQ{\Leftrightarrow}
\def \imply{\Rightarrow}
\def \isom{\cong}
\def \embed{\hookrightarrow}
\def \tensor{\otimes}
\def \wt#1{{\widetilde{#1}}}
\def \Jac{{\mathrm{Jac}}}
\def \red{{\mathrm{red}}}
\def \tor{{\mathrm{tor}}}
\def \lb{{\mathrm{lb}}}
\def \Res{{\mathrm{Res}}}
\def \sing{{\mathrm{sing}}}
\def \intsc{{\mathrm{Intsc}}}
\def \Aut{{\mathrm{Aut}}}

\begin{abstract}
We proved that every rational curve in the primitive class of a general K3
surface is nodal.
\end{abstract}

\title{Singularities of Rational Curves on K3 surfaces}
\author{Xi Chen}
\address{UCLA Department of Mathematics\\
6363 Math Sciences\\
Box 951555\\
Los Angeles, CA 90095-1555}
\email{xchen@math.ucla.edu}
\date{\today}
\maketitle

\section{Introduction}

The purpose of this paper is to prove the following theorem.

\begin{thm}\label{T:NOD}
For $g\ge 3$,
all rational curves in the linear system $|\CO_S(1)|$ on a general 
primitive K3 surface $S$ in ${\mathbb P}^g$ are nodal.
\end{thm}

The motivation to study this problem has been explained in
\cite{C}. Basically, we want to justify the beautiful formula of Yau and
Zaslow \cite{Y-Z}, which counts the number of rational curves
in $|\CO_S(1)|$ on a primitive K3 surface $S\subset \P^n$. The primary
consequence
of \thmref{T:NOD} is that the formula of Yau and Zaslow actually gives the
number of rational curves in $|\CO_S(1)|$ on a general K3 surface $S\subset
\P^n$.

It has been proved in \cite{C} that \thmref{T:NOD} is true for $g \le 9$
and $g = 11$ by degenerating a general K3 surface to a trigonal K3
surface. However, for $g$ large, we have to further degenerate a trigonal
K3 surface. The complexities involved in this process prevent us carrying
out the proof for any $g$. Although here we still use a degeneration
argument, our approach is entirely different. Here is a rough
sketch of the proof.

We start with the degeneration of a K3 surface to the union of two rational
surfaces. Let $X\to \Delta$ be a family of K3 surfaces of genus $g$
over the disk $\Delta$ whose central fiber $X_0 = R = R_1\cup R_2$ is
the union of two rational surfaces $R_1$ and $R_2$ which meet transversely
along an elliptic curve $E = R_1\cap R_2$. We may choose $R_i$ to be
$\P^1\times\P^1$ if $g$ is odd and choose $R_i$ to be $\F_1$ 
if $g$ is even for $i=1,2$. Let us consider the case that $g = 2k+1$ is
odd. We may construct $X$ in such a way that the limit of primitive line
bundles $\CO_{X_t}(1)$ on $X_t$ is the line bundle $\CO_R(1)$ on $X_0 = R$,
whose restriction to each $R_i\isom \P^1\times\P^1$ is the line bundle of type
$(1, k)$. So if we have a family $\Upsilon\subset X$ of rational curves over
$\Delta$, whose general fiber $\Upsilon_t$ is a rational curve
in the linear series $|\CO_{X_t}(1)|$ for each $t$, the central fiber
$\Upsilon_0$ will be a curve in the linear series $|\CO_R(1)|$ and hence
$\Upsilon_0 = \Sigma_1\cup \Sigma_2$ where $\Sigma_i$ is a curve of type
$(1, k)$ on $R_i$. Our trivial observation is that $\Upsilon_t$ is nodal if
$\Upsilon_0$ is nodal. However, $\Upsilon_0$ could fail to be nodal where

\begin{enumerate}
\item it has a reduced (i.e. isolated) singularity other than a node; or
\item it is nonreduced.
\end{enumerate}

The first case turns out much easier to handle than the second. This is
basically due to the fact that each $\Sigma_i$ is a curve of type $(1, k)$
on $\P^1\times\P^1$. So all the isolated singularities of $\Sigma_i$ are
nodes. If $\Upsilon_0 = \Sigma_1\cup \Sigma_2$ has an isolated singularity
other than a node, it must be one of the intersections between $\Sigma_1$ and
$\Sigma_2$ on $E$. The deformation of such singularities has been studied
in \cite{C}. With a bit more care, we are able to show that these
singularities deform to nodes on the general fiber $\Upsilon_t$. However,
if $\Upsilon_t$ is a rational curve in a multiple of the primitive class,
$\Sigma_i$ might have isolated singularities other than nodes which have to
be taken care of. This is one of the major obstacles to generalize
\thmref{T:NOD} to all rational curves on K3 surfaces.

To handle the second case, i.e., to handle the nonreduced components of
$\Upsilon_0$, we first divide them into three types, which we will call
Type I, II or III chain (see \secref{sec1}), respectively. The deformation
of $\Upsilon_0$ along a Type I chain is studied in \secref{sec2}. The basic
technique used there is to normalize the total family along the Type
I chain after a suitable base change. The deformation of $\Upsilon_0$ along
a Type II chain is studied in \secref{sec3}, where we build our argument
upon a lower bound estimation on the $\delta$-invariant of $\Upsilon_t$ in the
neighborhood of a Type II chain. The deformation of $\Upsilon_0$ along a
Type III chain is studied in \secref{sec4}. This turns out to be the
hardest case among the three. A two-stage degeneration is used, First, we
degenerate a general K3 surface to an elliptic K3 (see \secref{sec4}); and
then we degenerate an elliptic K3 to the union of two rational surfaces
described above.
The degeneration of a K3 surface to an elliptic K3 is also an
important step in Bryan and Leung's work \cite{B-L}, although the elliptic
K3 surfaces they used are different from the ones we use.

As a side note, there have been several progresses made on the enumeration
problems on K3 surfaces following Yau and Zaslow's
work. A. Beauville pointed out that the numbers Yau and Zaslow obtained are
the numbers of rational curves in $|\CO_S(1)|$ with each curve
counted with certain multiplicity \cite{B}, the multiplicity of a rational
curve only depends on its singularities and is 1 if the curve is nodal. He gave
an algebraic definition of the multiplicy. Later B. Fantechi, L. G\"ottsche
and D. Straten proved that the multiplicy assigned by Beauville to a
rational curve is positive. Recently, J. Bryan and N.C. Leung obtained
Yau-Zaslow's formula via a completely different approach \cite{B-L}.

\medskip\noindent{\bf Conventions.}

\begin{enumerate}
\item Throughout the paper, we will work exclusively over 
${\mathbb C}$.
\item Here a general K3 surface $S$ in $\P^g$ refers to 
a general primitive K3 surface in $\P^g$, where the number $g$ is called
the genus of $S$ by convention.
\item Since we are working over ${\mathbb C}$, we will use analytic geometry
whenever possible. Hence we will use analytic neighborhoods of points
instead of Zariski open neighborhoods in most cases,
while you may always replace them by formal or etale neighborhoods.
\end{enumerate}

\medskip\noindent{\bf Acknowledgments.}

I would like to thank Joe Harris for suggesting the problem and helping me
through the initial stage of the proof. I also
benefited greatly from the discussions with Mark Green, James McKernan
and Ziv Ran.

\section{Degeneration of K3 Surfaces and Limiting Rational Curves}
\label{sec1}

\subsection{Degeneration of K3 surfaces}

V. Kulikov classified all the possible degenerations of K3 surfaces
in \cite{K}. Here we only need one of the simplest cases: 
the degeneration of K3 surfaces into a union of two rational scrolls as
used in \cite{CLM}. This is also the degeneration used
in the proof of the existence of rational curves on K3 surfaces in \cite{C}. 

Following the notations in \cite{CLM}, let $R = R_1 \cup R_2$ be the 
union of two rational surfaces $R_1$ and $R_2$, which meet transversely
along a smooth elliptic curve $E = R_1\cap R_2$. For our purpose, we only
need the cases that either $R_1, R_2\isom \P^1\times \P^1$ or $R_1,
R_2\isom \F_1$, where $\F_1$ is the rational ruled surface
$\P(\CO \oplus \CO(-1))$ over $\P^1$.

We may represent such $R$ by the
tuple $(E, i_1, i_2)$ where $i_1: E\to R_1$ and $i_2: E\to R_2$ are
the embeddings of $E$ to $R_1$ and $R_2$, respectively. Two unions $R$ and
$R'$ represented by $(E, i_1, i_2)$ and $(E', i_1', i_2')$ are isomorphic
if and only if $E\isom E'$ and there exist isomorphisms $\varphi : E'\to E$,
$\phi_1: R_1\to R_1$ and $\phi_2: R_2\to R_2$ such that $i_1' = \phi_1\circ
i_1 \circ \varphi$ and $i_2' = \phi_2\circ i_2 \circ \varphi$. Then it is not
hard to see that such $R$'s form an irreducible moduli space of dimension $4$.

Let $C_i$ and $F_i$ be two generators
of $\Pic(R_i)$ with $C_i\cdot F_i = 0$, 
$F_i^2 = 0$ and $C_i^2 = 0$ if $R_i\isom
\P^1\times \P^1$; $C_i^2 = -1$ if $R_i\isom \F_1$, for $i=1,2$.
We use the notation $\CO_R(a C + b F)$ to denote the
line bundle on $R$ whose restrictions to $R_i$ are $\CO_{R_i}(a C_i + b F_i)$,
if such line bundle exists.

It is not hard to see that
the dualizing sheaf $\omega_R$ of $R$ is trivial and 
$H^1(R, \BZ) = 0$. So
it is expected that a general deformation of $R$, say $X\to \Delta$ with $X_0 =
R$, is a family of K3 surfaces. But since a general $R$ is not projective,
the general fiber of $X$ is not algebraic.
It is easy to see that $R$ is
projective if and only if there exists two ample line bundles $L_1\in \Pic
R_1$ and $L_2 \in \Pic R_2$ such that $L_1 |_E = L_2 |_E$. 
One obvious choice of $L_i$ is $L_i = \CO_{R_i}(C_i + kF_i)$. 
So we are considering the unions $R$ of the following type.

Let $R = R_1\cup R_2$ be a union of rational surfaces described as above,
which further satisfies $\CO_{E}(C_1+kF_1) = \CO_E(C_2+kF_2)$
for some $k\ge 1$ if
$R_i\isom \P^1\times\P^1$ and $k\ge 2$ if $R_i\isom \F_1$. 
We will call such $R$ {\it a union of scrolls of genus $g$\/}, 
where $g = 2k+1$ if $R_i\isom \P^1\times\P^1$ and $g = 2k$ if
$R_i\isom \F_1$.

Notice that the relation $\CO_E(C_1 + kF_1) = \CO_E(C_2 + kF_2)$
imposes one extra condition on the tuple $(E, i_1, i_2)$
which represents $R$. So the moduli space of
unions of scrolls of fixed genus $g$ has dimension $3$. With a little extra
effort, one can see that the moduli space is irreducible.

We are also interested in a special union of scrolls $R = R_1\cup R_2$
which satisfies $\CO_{E}(C_1) =
\CO_{E}(C_2)$ and $\CO_{E}(F_1) = \CO_{E}(F_2)$. 
We will call such $R$ {\it a degenerated (or special) union of scrolls\/}.
It follows from the similar argument as before that
degenerated unions of scrolls form an irreducible moduli
space of dimension $2$.

It was proved in \cite{CLM} that a general union of scrolls $R$ of genus
$g$ lies on
the boundary of a complete family of K3 surfaces of genus $g$. The
construction is straightforward to carry out by embedding $R$ to
$\P^g$ by the complete linear series $|\CO_R(C + k F)|$.
Then $R$ lies on the
component of the Hilbert scheme whose general point is a primitive
K3 surfaces in $\P^g$.

However, we need a little bit more. We want to find a complete family of
K3 surfaces of genus $g$ whose boundary also contains degenerated unions of
scrolls. The
previous construction fails since $\CO_R(C+kF)$ is ample but not very
ample on a degenerated union of scrolls $R$:
the morphism $R\to \P^g$ given by $|\CO_R(C+kF)|$
maps $R$ to a double scroll. The remedy to this situation is
trivial. Instead of using $|\CO_R(C+kF)|$, we embed $R$ to a projective
space by the complete linear series $|\CO_R(l(C+kF))|$ for some large
$l$. Actually, it is enough to take $l = 2$. Namely, we embed $R$
to $\P^{4g-3}$ by $|\CO_R(2C + 2kF)|$. 
And we can show that $R$ lies on the component of the
Hilbert scheme whose general point is a K3 surface of genus $g$ in
$\P^{4g-3}$. The argument for this statement
is identical to that in \cite[Sec. 2.2]{CLM} and
we will only formulate it in the following proposition without
the proof.

\begin{prop}\label{prop1a}
Let $R$ be a union of scrolls of genus $g$ which is
embedded into $\P^{4g - 3}$ by $|\CO_R(2C+2kF)|$, where $k = \lfloor g/2
\rfloor$. 
Let $N_R$ be the normal bundle of $R\subset \P^{4g-3}$ and $T_R^1 =
\ext^1(\Omega_R, \CO_R)$. Then
\begin{enumerate}
\item $H^1(N_R) = 0$;
\item $\dim H^0(N_R) = \dim \PGL (4g-2) + 19$;
\item the natural map from $H^0 (N_R)$ to $H^0(T_R^1)$ is surjective;
\item $R$ is represented by a smooth point in the component $\H_g$ 
of the Hilbert scheme whose general point is a K3 surface of genus $g$ in
$\P^{4g-3}$.
\end{enumerate}
\end{prop}

Let $\CR_g\subset \H_g$ be the locus in $\H_g$ consisting
of points representing unions of scrolls and 
let $\CR_g^0\subset\CR_g$ be the locus consisting of points representing
degenerated unions of scrolls. Every automorphism of a union of scrolls
$R = R_1 \cup R_2$ induces an automorphism of $R_1$ which maps the double
curve $E = R_1\cap R_2$ to itself. Obviously, there are only finitely many
automorphisms of $R_1$ with this property. So the automorphism group of $R$
is finite. Therefore, $\dim \CR_g= \dim \CR_g^0 +1 = \dim \PGL (4g-2) + 3$.

Let $\wt{\H_g}$ be the blowup of $\H_g$ along the closure of $\CR_g$ 
and let $\wt{\CR_g}$ and $\wt{\CR_g^0}$
be the total transforms of $\CR_g$ and $\CR_g^0$ under the map
$\wt{\H_g} \to \H_g$, respectively. 

For any $[R]\in  \CR_g$, it is not hard to see that
the tangent space $T_{\CR_{g}, [R]}$ of $\CR_{g}$ at
$[R]$ lies inside the kernel of the surjection
$T_{\H_g, [R]}= H^0(N_R)\to H^0(T_R^1)$. On the other hand, $\dim \CR_g = 
\dim H^0(N_R) - \dim H^0(T_R^1)$. Therefore, $\CR_g$ is smooth
everywhere.
Hence $\wt{\H_g}$ is smooth in the neighborhood of $\wt{\CR_g}$.

Let $\S_g = \{ ([X], p): p\in X\} \subset\H_g\times\P^{4g-3}$
be the universal family over $\H_g$ and $\wt{\S_g} =
\S_g\times_{\H_g} \wt{\H_g}$.

It is not hard to see that every
point of $\wt{\CR_g}$ can be uniquely
represented by a pair $([R], s)$ with $[R]\in
\CR_g$ and $s\in \P H^0(T_R^1)$, which we will write as $[R^s]$.
Notice that $T_R^1$ is a sheaf supported
along $E$ whose restriction to $E$ is the line bundle $N_{E/R_1}\tensor
N_{E/R_2}$, where $N_{E/R_1}$ and $N_{E/R_2}$ are the normal bundle of $E$
in $R_1$ and $R_2$, respectively. A geometric interpretation of this blowup
process can be put as follows. Let $\pi: \Delta\to \wt{\H_g}$ be a morphism
from the disk $\Delta$ to $\wt{\H_g}$,
where $\pi(0) = [R^s] \in \wt{\CR_g}$ and $\pi(\Delta)\not\subset
\wt{\CR_g}$. Then the corresponding family 
$X =\wt{\S_g} \times_\wt{\H_g} \Delta \to \Delta$ has sixteen
rational double points which are the vanishing loci of $s$.

\begin{rem}\label{rem1}
Here we use the term ``rational double points'' in a broad sense. Let $X\to
\Delta$ be a one-parameter family of surfaces over the disk $\Delta$.
A point $p\in X$ is called a rational double point of the family $X$ over
$\Delta$ if $X$
is locally isomorphic to
$\Spec \BC[[x, y, z, t]]/(xy - t^\alpha z)$ at $p$, 
where $t$ parameterizes the disk $\Delta$ and
$\alpha$ is a positive integer. If $\alpha = 1$, $p$ is a rational
double point of the threefold $X$ in the usual sense, which can be resolved
by blowing up $X$ at $p$ and then blowing down along one of
the rulings of the
exceptional divisor $\P^1\times\P^1$. If $\alpha > 1$,
we may resolve $p$ in the same way as in \cite[Appendix C, p. 39]{G-H}.
But we will take a more direct approach here by choosing an neighborhood
$U$ of $p$ and simply letting $\wt{U}\subset U\times \P^1$ be defined by
\begin{equation}\label{e1z}
\frac{x}{t^\alpha} = \frac{z}{y} = \frac{Z_0}{Z_1}
\end{equation}
where $(Z_0, Z_1)$ are the homogeneous coordinates of $\P^1$. It is trivial
to glue $\wt{U}$ with $X\backslash \{p\}$
and arrive at a ``resolution'' $\wt{X}$ of $X$ (rigorously, $\wt{X}$ is not
a desingularization of $X$ since $\wt{X}$ is still singular in dimension
one; but now $\wt{X}$ is isomorphic to $\Spec \BC[[x,
y, z, t]]/(xy - t^\alpha)$ everywhere along its singular locus, which is
all we need).
Of course, we may switch $x$ and $y$ in \eqref{e1z} and arrive at another
resolution of $X$. Obviously, this corresponds to the ``flip'' phenomenon of
the resolutions of an ordinary three-fold double point.
\end{rem}

Let $\W_g$ be the incidence correspondence
\[
\begin{split}
\W_g = & \big\{([S], [C]): [S]\in \wt{\H_g} \text{ general},\\
& \quad\quad\text{$C$ is a rational curve in the primitive
class of $S$}\big\}\\
& \subset \wt{\H_g} \times |\CO_{\P^{4g-3}}(1)|,
\end{split}
\]
where we construct $\W_g$ as a subvariety of $\wt{\H_g} \times
|\CO_{\P^{4g-3}}(1)|$ by identifying $[C]$ with
$2C\in |\CO_S(2) | = |\CO_{\P^{4g-3}}(1)|$ (let $\CO_S(1)$ be the primitive
line bundle on $S$).
Let $\wt{\W_g}$ be the closure of $\W_g$ in $\wt{\H_g} \times
|\CO_{\P^{4g-3}}(1)|$.

Obviously, \thmref{T:NOD} is equivalent to the following statement: 
for every point $([S], [C])$ on the fiber of $\wt{\W_g}\to
\wt{\H_g}$ over a general point $[S]\in \wt{\H_g}$,
$C$ is a nodal curve.
Obviously, it suffices  to verify this statement for every
irreducible component of $\wt{\W_g}$ which dominates $\wt{\H_g}$. So, from
now on, we will pretend that $\wt{\W_g}$ is irreducible and dominates
$\wt{\H_g}$.

It is obvious that the map $\wt{\W_g}\to \wt{\H_g}$ is generically
finite.
One important step to prove \thmref{T:NOD} is to figure out what kind of
curves $[C]$ lie on the fiber $(\wt{\W_g})_{[R^s]}$
of $\wt{\W_g}\to \wt{\H_g}$ over a point $[R^s]\in \wt{\CR_g}$.
Of course, $C\in |\CO_R(C+kF)|$. Besides, $C$ must be the limit of
rational curves on K3 surfaces. More precisely, $([R^s], [C])\in
(\wt{\W_g})_{[R^s]}$ if and only if there exist a family $X$ of K3 surfaces
over $\Delta$ and a family $\Upsilon\subset X$ of rational curves where
$X_0 = R$, $X$ has sixteen rational double
points which are the vanishing loci of
$s$, $\Upsilon_0 = C$ and $\Upsilon_t$ is a rational curve in the
primitive class of $X_t$. It turns out that
there are only finitely many curves $C$
in $|\CO_R(C+kF)|$ with this property.

\begin{defn}
Let $R$ be a union of scrolls and $s\in \P H^0(T_R^1)$. A curve $C\in
|\CO_R(C+kF)|$ with $([R^s], [C])\in \wt{\W_g}$ is called a limiting
rational curve on $R^s$. When there is no possibility for confusion, we
will abbreviate $R^s$ to $R$ and simply call $C$ a limiting rational
curve on $R$.
\end{defn}

\subsection{Classifications of limiting rational curves}

Let $\pi: \Delta \to \wt{\H_g}$ be a morphism from the disk $\Delta$ to
$\wt{\H_g}$ where $\pi(0) = [R^s]\in \wt{\CR_g}$ and $\pi(t)$
represents a general K3 surface of genus $g$ for $t\in \Delta$ general.

Let $X =\wt{\S_g} \times_\wt{\H_g} \Delta \to \Delta$ be the 
one parameter family of K3 surfaces
corresponding to $\pi$. 
Obviously, $X$ has exactly $16$ rational double points
$p_1, p_2, ..., p_{16}$ lying on the double curve $E = R_1\cap R_2$,
where $\{p_1, p_2, ...,
p_{16}\}$ is the vanishing locus of $s\in \P H^0(T_R^1)$.
For a general choice of $s$, $p_1, p_2, ...,
p_{16}$ satisfy the only relation
\begin{equation}\label{e2}
\CO_E(p_1+p_2+...+p_{16}) = N_{E/R_1}\tensor N_{E/R_2}
\end{equation}
where
$N_{E/R_i} = \CO_E(2C_i + 2F_i)$ if $R_i\isom \P^1\times\P^1$ and 
$N_{E/R_i} = \CO_E(2C_i + 3F_i)$ if $R_i\isom \F^1$.

Let $\Upsilon\subset X$ be a family of curves over $\Delta$
whose general fiber $\Upsilon_t$ is a rational curve in the primitive class
of $X_t$ for each $t$.
Let $\wt{\Upsilon}$ be the nodal reduction of $\Upsilon$.
Namely, $\wt{\Upsilon}$ is a family of stable maps to $X$
such that $\wt{\Upsilon}\to X$ factors through $\Upsilon$ and the general
fiber $\wt{\Upsilon}_t$ of $\wt{\Upsilon}$ is the normalization of the
general fiber $\Upsilon_t$ of $\Upsilon$. Notice that such
$\Upsilon$ and $\wt{\Upsilon}$ exist after a base change.

One of the key lemmas
we use to classify limiting rational curves on $X_0$ is Lemma
2.2 in \cite{C}. We need a slightly stronger version, which is formulated
and proved as follows.

\begin{lem}\label{lem1}
Let $X\subset \Delta_{x_1x_2...x_n}^n\times \Delta_t$ ($n\ge 3$) be 
the hypersurface cut out by $x_1x_2 =
t^\alpha$ for some $\alpha > 0$, where $\Delta_{x_1x_2...x_n}^n$ is the
$n$-dimensional polydisk parameterized by
$(x_1, x_2, ..., x_n)$ and $\Delta_t$ is the disk parameterized by $t$. Let
$X_0$ be the central fiber of $X$ over $\Delta_t$, let
$X_0 = R_1 \cup R_2$ where $R_1 = \{x_1 = t = 0\}$
and $R_2 = \{x_2 = t = 0\}$ and let $E = R_1 \cap R_2$.
Let $S$ be a flat family of analytic curves over
$\Delta_t$ and $\pi: S\to X$ be a proper morphism preserving
the base $\Delta_t$. Suppose
that the image $\pi(S_0)$ of the central fiber $S_0$ of $S$ meets 
$E$ properly along a 0-dimensional scheme supported at the origin.
Let $S_0 = \Gamma_1 \cup
\Gamma_2$ with $\pi(\Gamma_1)\subset R_1$ and $\pi(\Gamma_2)\subset R_2$.
Then $\pi(\Gamma_1) \cdot R_2 = \pi(\Gamma_2) \cdot R_1$
holds on $\{t = 0\}\isom \Delta_{x_1x_2...x_n}^n$.
\end{lem}

\begin{proof}
Without the loss of generality,
we may assume that $S$ is irreducible and smooth; otherwise, we apply the
following argument to each irreducible component of its desingularization.
Let $\eta: X\to X' = \{ x_1x_2 = t\} \subset \Delta_{x_1x_2...x_n}^n\times
\Delta_t$ be the morphism sending $(x_1, x_2, ..., x_n, t)$
to $(x_1, x_2, ..., x_n, t^\alpha)$ and let $p = \eta\circ
\pi$. Obviously, it suffices to show that $p_*(\Gamma_1) \cdot R_2 = 
p_*(\Gamma_2) \cdot R_1$ on $X'$.

Since $\Gamma_1\cdot p^*(R_1 + R_2) = 0$ and $(\Gamma_1 + \Gamma_2) \cdot
p^*(R_1) =0$ on $S$,
\[
\Gamma_1 \cdot p^*(R_2) = \Gamma_2\cdot
p^*(R_1).
\]
Therefore, $p_*(\Gamma_1) \cdot R_2 = 
p_*(\Gamma_2) \cdot R_1$ by the projection formula.
\end{proof}

The following proposition deals with the case that $X$ has a rational
double point at the origin, which can be viewed as a corollary of
\lemref{lem1}.

\begin{cor}\label{cor3}
Let $X, R_1, R_2, E, \pi, S$ and 
$S_0$ be defined as in \lemref{lem1} except that
$X$ is cut out by $x_1x_2 = t^\alpha x_3$ for some $\alpha > 0$.
Suppose that $\pi(S_0)$ contains an irreducible component $\Gamma_1\subset
R_1$ such that the tangent cones of $\Gamma_1$ and $E$ at the origin $p$
do not meet properly in the tangent space of $R_1$ at $p$.
Then $\pi(S_0)$ also contains a component $\Gamma_2\subset R_2$.
\end{cor}

\begin{proof}
We may resolve the double point of $X$ as in \remref{rem1}. Let
$\wt{X}\subset X\times \P^1$ be given by
\[
\frac{x_1}{t^\alpha} = \frac{x_3}{x_2} = \frac{Z_0}{Z_1}
\]
where $(Z_0, Z_1)$ are the homogeneous coordinates of $\P^1$. Let $\wt{S} =
S\times_X \wt{X}$. Obviously, $\wt{X}_0 = \wt{R_1} \cup R_2$, where
$\wt{R_1}$ is the blowup of $R_1$ along the subscheme $\{x_2 = x_3 = 0\}$.
Obviously, $\wt{R_1}$ and $R_2$ meet along $\wt{E} = \{x_1 = x_2 = z = t =
0\}$ (let $z = Z_1/Z_0$) and $\wt{X}$ is locally defined by
\begin{equation}\label{cor3-e1}
x_1 z = t^\alpha \text{ and } x_2 = z x_3
\end{equation}
in the neighborhood of $\wt{E}$. Let $\wt{p}\in \wt{E}$ be the point $\{x_1
= x_2 = ... = x_n = z = t = 0\}$ and let $\wt{\Gamma_1}\subset \wt{S}_0$ 
be the proper transform of $\Gamma_1$ under the map $\wt{R_1}\to R_1$.
The assumption that the tangent cones of $\Gamma_1$ and $E$ at $p$ do not
meet properly implies that their proper transforms
$\wt{\Gamma_1}$ and $\wt{E}$ under the blowup $\wt{R_1}\to R_1$ still meet
at a point on the exception divisor, i.e., they meet at $\wt{p}$.
Now by \eqref{cor3-e1}, $\wt{X}$ is locally defined by $x_1 z = t^\alpha$ in
$\Delta_{x_1 z x_3 ... x_n}^n\times \Delta_t$ at $\wt{p}$. And
$\wt{\pi}(\wt{S}_0)$ has a component $\wt{\Gamma_1}$ passing through
$\wt{p}$ and lying on $\wt{R_1}$, where $\wt{\pi}: \wt{S}\to \wt{X}$ is the
map induced by $\pi$. So by \lemref{lem1}, $\wt{\pi}(\wt{S}_0)$ contains a
component $\Gamma_2\subset R_2$, i.e., $\pi(S_0)$ contains the component
$\Gamma_2\subset R_2$.
\end{proof}

Let $\Upsilon_0 = \Sigma_1 \cup \Sigma_2$ with $\Sigma_i \subset R_i$ for
$i=1,2$. And let $\Gamma_i\subset \Sigma_i$ be the irreducible component of
$\Sigma_i$ in $|C_i + k_i F_i|$ for some $k_i \le k$ ($i = 1,2$). 

Let $r_{ij}$ be all the points on $E$ satisfying $\CO_E(2 r_{ij}) =
\CO_E(F_i)$ for $i=1,2$ and $j = 1,2,3,4$. Notice that if $R$ is a
degenerated union of scroll, we have $\{r_{1j}\} = \{r_{2j}\}$. So we will
simply let $r_{1j} = r_{2j}$ for $j = 1,2,3,4$ if $[R]\in \CR_g^0$.

For two points $p$ and $q$ on $E$ satisfying $\CO_E(p+q)=\CO_E(F_1)$ or
$\CO_E(F_2)$, we use the notation $\overline{pq}$ to denote the curve in
$|F_1|$ or $|F_2|$ passing through $p$ and $q$. If $R$ is general in
$\CR_g$, there is no ambiguity; otherwise, we use $\overline{pq}^{(1)}$ and
$\overline{pq}^{(2)}$ to distinguish on which of $R_1$ and $R_2$ this curve
lies.
And we use $\overline{r_{ij}}$ to denote
the curve in $|F_i|$ passing through $r_{ij}$.

Since there
are exactly two components on $\wt{\Upsilon}_0$ dominates $\Gamma_1$
and $\Gamma_2$, respectively,
we will continue to use $\Gamma_1$ and $\Gamma_2$ to denote
these two components on $\wt{\Upsilon}_0$.

\begin{defn}
Let $C$ be a stable curve. The ``dual graph'' $G$ of $C$ is a graph
constructed by
representing each component $A$ of $C$ by the vertex $[A]$ and drawing an edge
between two vertices $[A]$ and $[B]$ if the corresponding curves $A$ and
$B$ meet at a point. We allow multiple edges between $[A]$ and $[B]$ if
they meet at more than one point; and if $A$ has a node, we will draw a
loop around $[A]$. Let $\deg([A])$ denote the degree of the vertex $[A]$
in $G$.

A sequence of components $C_1, C_2, ..., C_n$ of $C$ form a chain in $G$
if there is an edge between $[C_i]$ and $[C_{i+1}]$ for $i =
1,2, ..., n-1$.     
\end{defn}

First, we will show that there are only finitely many limiting rational curves
on $X_0$. Namely, there are only finitely many possible configurations for
$\Upsilon_0$.

\begin{prop}\label{prop1b}
Let $[R^s]\in \wt{\CR_g^0}$ with
$s\in \P H^0(T_R^1)$ general. Then
there are only finitely many limiting rational curves on $R^s$. Namely,
the fiberation of $\wt{\W_g}\to \wt{\H_g}$ is finite over $[R^s]$.
\end{prop}

\begin{proof}
Since $[R^s]\in \wt{\CR_g^0}$, $\CO_E(C_1) = \CO_E(C_2) = \CO_E(C)$ and
$\CO_E(F_1) = \CO_E(F_2) = \CO_E(F)$. It is not hard to see that if
$\Sigma_1$ contains a component $\overline{pq}^{(1)}$, 
$\Sigma_2$ must contain $\overline{pq}^{(2)}$ with the same multiplicity
and vice versa. Therefore, we necessarily have $\Gamma_1 \cap E = 
\Gamma_2 \cap E$.

Let $q_i$ be the point on $E$ such that $\CO_E(p_i + q_i) = \CO_E(F)$
for $i=1,2,...,16$. Let $r_j = r_{1j} = r_{2j}$ for $j = 1,2,3,4$.

By \eqref{e2}, the only relation among $p_1, p_2, ..., p_{16}$ is
\(\CO_E(p_1+p_2+...+p_{16}) = \CO_E(4C + 4F)\) or \(\CO_E(4C + 6F)\). 
Therefore, the subgroup of $\Pic(E)$ generated by $p_l$, $q_l$ and $r_j$ is
\[
\dsum_{l=1}^{15} \BZ p_l \oplus \BZ r_1 \oplus \BZ (4C) \oplus (\BZ_2
\oplus \BZ_2)
\]
where $\BZ_2\oplus \BZ_2$ is the subgroup consisting of $r_j - r_1$ for $j
= 1,2,3,4$. This group clearly
does not contain the divisor $C$ when $p_1, p_2, ..., p_{15}$ and $E$
are general. Therefore, $\Gamma_i$ must meet $E$ at (at least)
one point other than $p_l$, $q_l$ and $r_j$ for $l = 1,2,..., 16$ and
$j=1,2,3,4$.

Suppose that $p\in \Gamma_i\cap E$ and $p\not\in \{p_1, p_2, ..., p_{16}, 
q_1, q_2, ..., q_{16}\}$. Let $q$ be
the point on $E$ such that $\CO_E(p+q) = \CO_E(F)$. We claim that 
on $\wt{\Upsilon}_0$, the curves $\Gamma_1$ and $\Gamma_2$ are joined by a
chain of curves which either contract to one of the two points $p$ and $q$ 
or dominate one of the two curves $\overline{pq}^{(1)}$ and
$\overline{pq}^{(2)}$.

Let $D_1, D_2, ..., D_\gamma$ be the components of $\wt{\Upsilon}_0$ which 
either contract to one of $p$ and $q$ or dominate one of $\overline{pq}^{(i)}$.
Let $G$ be the dual graph of $\Gamma_1\cup \Gamma_2\cup D_1\cup D_2\cup ...
\cup D_\gamma$ with each curve $A$ being represented by the vertex $[A]$.

If $p\not\in \{r_1, r_2, r_3, r_4\}$,
then $\deg([\Gamma_i])\ge 1$ and $\deg([D_j])\ge 2$ for
$i=1,2$ and $j=1,2,...,\gamma$ by \lemref{lem1}. So $G$ either contains a
cycle, which is impossible, or $G$ is connected, which implies our claim.

If $p = q = r_j$ for some $j$, then $\Gamma_i$ must meet $E$ transversely
at $r_j$; otherwise, $\Gamma_i$ will meet $\overline{r_{ij}}$ at $r_j$ with
multiplicity at least $2$, which is impossible. 
Locally at $r_j$, $\Upsilon_0$ consists of $\Gamma_1\cup \Gamma_2$ and a
multiple of $\overline{r_{1j}}\cup \overline{r_{2j}}$. Notice that
the local intersection number between $\Gamma_i$ and $E$ at $r_j$ is 1
(odd) and the local intersection number between
$\overline{r_{ij}}$ and $E$ at $r_j$ is 2 (even). Therefore,
$[\Gamma_1]$ and $[\Gamma_2]$ must lie on the
same connected component of $G$ by \lemref{lem1}, which also implies our
claim.

In conclusion, $\Gamma_i$ meet $E$ at (at least) one point
other than $p_l$, $q_l$ and $r_j$. If $\Gamma_i$ meet $E$ at $p$ and $p\not\in
\{p_1, p_2, ..., p_{16}, q_1, q_2, ..., q_{16}\}$,
then $\Gamma_1$ and $\Gamma_2$ are joined by a chain of
curves on $\wt{\Upsilon}_0$ whose image lie in
$\overline{pq}^{(1)}\cup \overline{pq}^{(2)}$, where $q\in E$
is the point satisfying that $\CO_E(p+q) = \CO_E(F)$.
So it is not hard to see that $\Gamma_i$ cannot meet $E$ at more than one
point other than $p_l$ and $q_l$; otherwise, $\Gamma_1$ and $\Gamma_2$ will
be joined by two different chains of curves. In conclusion, $\Gamma_i$
meets $E$ at exactly one point $p$ other than $p_l$ and $q_l$ for
$l=1,2,...,16$ and $p\ne r_j$ for $j = 1,2,3,4$. Obviously, there are at
most finitely many curves in $|C_i + k_i F_i|$ ($k_i \le k$) with this
property.

It remains to show that there are only finitely many possible
configurations for a component $\overline{pq}^{(i)}\subset
\Sigma_i$. Actually, we claim that if $\overline{pq}^{(i)}\subset
\Sigma_i$, then

\begin{enumerate}
\item $p, q\in \{p_1, p_2, ..., p_{16}, q_1, q_2, ..., q_{16}\}$; OR
\item $p = q\in \{r_1, r_2, r_3, r_4\}$; OR
\item $p$ or $q$ lies on $\Gamma_i$.
\end{enumerate}

Suppose that $\overline{pq}^{(i)}\subset \Sigma_i$,
$p,q\not\in \{p_1, p_2, ..., p_{16}, q_1, q_2, ..., q_{16},
r_1, ..., r_4\}$ and
$p, q\not\in \Gamma_i$. Let $D_1, D_2, ..., D_\gamma$ be the
components of $\wt{\Upsilon}_0$ which either contract to one of $p$ and $q$
or dominate one of $\overline{pq}^{(i)}$. Let $G$ be the dual graph of
$D_1\cup D_2\cup ...\cup D_\gamma$ with each curve $A$ represented by the
vertex $[A]$.
Since $p\ne q$, $p, q\not\in \{p_1, p_2, ...,p_{16}\}$ and $\Gamma_i$ does
not pass through $p$ and $q$, $\deg([D_j])\ge 2$ for $j=1,2,...,\gamma$ by
\lemref{lem1}. So $G$ must contain a cycle, which is a contradiction.

In summary, there are at most finitely many possible configurations for
$\Upsilon_0 = \Sigma_1 \cup \Sigma_2$. So $\wt{\W_g}$ is finite over
$[R^s]\in \wt{\CR_g^0}$.
\end{proof}

\begin{rem}
Since $\wt{\CR_g^0}\subset \wt{\CR_g}$, \propref{prop1b} implies
that the fiberation of $\wt{\W_g}\to\wt{\H_g}$ is finite over 
a general point $[R^s] \in\wt{\CR_g}$.
\end{rem}

In the rest of this section, we will assume that the choice of $R^s$
in the construction of $X$ is general in $\wt{\CR_g}$.

\begin{defn}
For each point $p\in \Sigma_i\cap E$ ($i = 1,2$), we use
$m_i(p)$ to denote the local
intersection multiplicity between $\Gamma_i$ and $E$ at
$p$ (let $m_i(p) = 0$ if $\Gamma_i$ does not pass through $p$).

For each $\overline{pq}\subset \Upsilon_0$,
let $\mu(\overline{pq})$ denote the multiplicity of
$\overline{pq}$ in $\Upsilon_0$.
\end{defn}

\begin{defn}
An $F$-chain on $R$ is a union of $m$ distinct curves 
$C = \overline{q_0q_1}\cup \overline{q_1 q_2} ...\cup \overline{q_{m-1} q_m}$,
where, as the notation suggests, $\overline{q_l q_{l+1}}$ is a curve in
either $|F_1|$ or $|F_2|$. We call $m$ the length of $C$. If $m = 0$,
we let $C$ be the point $q_0$.

A maximal $F$-chain $C$ in $\Upsilon_0$ is an $F$-chain
$C \subset \Upsilon_0$ and it is maximal in the sense that there does not
exist an $F$-chain $C'\subset \Upsilon_0$ containing $C$ as a proper subset.
\end{defn}

Since we will always deal with maximal
$F$-chains in $\Upsilon_0$, we will simply
call them $F$-chains in $\Upsilon_0$.

We divide the $F$-chains in $\Upsilon_0$ into three types.

\begin{enumerate}
\item A chain $C$ is called a {\it Type I\/} chain if
$p_l\not\in C$ for $l=1,2,...,16$ 
and $r_{ij}\not\in C$ for $i=1,2$ and
$j=1,2,3,4$. 
\item A chain $C$ is called a {\it Type II\/} chain if
$p_l\in C$ for some $1\le l\le 16$.
\item A chain $C$ is called a {\it Type III\/} chain if
$r_{ij}\in C$ for some $1\le i\le 2$ and $1\le j\le 4$.
\end{enumerate}

Notice that a Type I or II chain could consist of a single point.

\begin{prop}\label{prop1}
There is exactly one Type I chain in $\Upsilon_0$, whose length is even.
Let $\overline{q_0q_1}\cup \overline{q_1 q_2}\cup ...
\cup \overline{q_{2l-1} q_{2l}}\subset \Upsilon_0$
be the Type I chain and assume that
$\overline{q_0q_1}\subset R_1$ without the loss of generality.
Then
\begin{gather*}
\mu(\overline{q_0q_1}) = 
\mu(\overline{q_1q_2}) = ... = \mu(\overline{q_{2l-1} q_{2l}}) =\beta =
m_2(q_0) = m_1(q_{2l}) = \beta\text{ and }\\
\begin{split}
m_1(q_0) &= m_1(q_1) = ... = m_1(q_{2l-1})\\
&= m_2(q_1) = m_2(q_2) = ... = m_2(q_{2l}) = 0
\end{split}
\end{gather*}
for some $\beta > 0$. There are exactly $2l$ components $D_1, D_2, ...,
D_{2l}$ on $\wt{\Upsilon}_0$ where each $D_i$ dominates $\overline{q_{i-1}
q_{i}}$ with a degree $\beta$ map totally ramified at $q_{i-1}$ and $q_{i}$. 
The components $\Gamma_2, D_1, D_2, ..., D_{2l}, \Gamma_1$ form a chain in
the dual graph of $\wt{\Upsilon}_0$.
\end{prop}

\begin{proof}
Obviously, every point in $C\cap E$ for a Type II or Type III chain
$C$,
as an element in the
Picard group $\Pic(E)$, lies in the subgroup generated by $p_l$
($l=1,2,...,16$) and $r_{ij}$ ($i=1,2$ and $j=1,2,3,4$), which, 
by \eqref{e2}, is
\[
\left(\dsum_{l=1}^{15} \BZ p_l\right)
\oplus \left(\dsum_{i=1}^2 \BZ r_{i1}\right) 
\oplus \BZ (2C_1 + 2C_2)
\oplus (\BZ_2\oplus \BZ_2)^{\oplus 2}
\]
where two copies of $\BZ_2\oplus \BZ_2$ are the subgroups consisting of 
$r_{ij} - r_{i1}$ ($j = 1,2,3,4$) for $i = 1,2$, respectively.
Obviously, $C_1 + kF_1 = C_2 +k F_2$ does not lie in this group, i.e., it
is not generated by $p_l$ and $r_{ij}$, for $R$ and $s$ general.  
Therefore there is at least one Type I chain.

Let $\overline{q_0q_1}\cup \overline{q_1q_2}\cup ...
\cup \overline{q_{m-1} q_m}$ form a Type I chain in $\Upsilon_0$.
It is not hard to see that $q_0, q_m\in \Gamma_1\cup\Gamma_2$.

We let
$D_1, D_2, ..., D_\gamma$ be the components on $\wt{\Upsilon}_0$
dominating the curves
$\overline{q_0q_1}, \overline{q_1q_2}, ..., \overline{q_{m-1}q_m}$ and
$H_1, H_2, ..., H_\mu$ be the components contracting to the points $q_0,
q_1, q_2, ... q_m$.

Let $G$ be the dual graph of $\Gamma_1, D_1, D_2, ...,
D_\gamma, \Gamma_2, H_1, H_2, ..., H_\mu$ with each curve $A$ being
represented by the vertex $[A]$.
Obviously, $G$ contains no circuit and is consequently a ``forest'' (a
disjoint union of trees).
By \lemref{lem1}, $\deg([D_i]) \ge 2$ for $i=1,2,...,\gamma$ and
$\deg([\Gamma_1]) + \deg([\Gamma_2]) \ge 2$.
And since $\deg([H_j])\ge 3$, $G$ has at
least $\gamma + 1 + 3\mu/2$ edges with $\gamma + 2+\mu$ vertices.
On the other hand, 
$G$ has at most $\gamma + 1 + \mu$ edges since it is a forest. So we must
have
\begin{enumerate}
\item $\mu = 0$, i.e., there are no components of $\wt{\Upsilon}_0$
contracting to the points $q_0, q_1, q_2, ..., q_m$;
\item $G$ is a tree;
\item $\deg([D_i]) = 2$ for $i=1,2,...,\gamma$;
\item $\deg([\Gamma_1]) = \deg([\Gamma_2]) = 1$.
\end{enumerate}
And since we assume that $\overline{q_0q_1}\subset R_1$, we must have
$q_0\in \Gamma_2$ and $q_m\in\Gamma_1$. Hence $m$ must be even, say $m=2l$.

The picture of $G$ is very clear now. The vertices of $G$ form a chain
after some ordering.
Without the loss of generality, we may assume that
$[\Gamma_2], [D_1], [D_2], ..., [D_\gamma], [\Gamma_1]$
form a chain in that order. 
Obviously, $D_1$ must dominate $\overline{q_0q_1}$ and there
is no other $D_i$ ($i\ne 1$) dominating $\overline{q_0q_1}$; otherwise,
$\deg([\Gamma_2]) \ge 2$ by \lemref{lem1}. Similarly, $D_2$ must dominate
$\overline{q_1q_2}$ and there is no other $D_i$ dominating
$\overline{q_1q_2}$; otherwise,
$\deg([D_1]) \ge 3$ by \lemref{lem1}. This line of argument goes on and
finally shows that $\gamma = 2l$ and 
each $D_i$ dominates $\overline{q_{i-1}q_{i}}$
for $i=1,2,...,2l$. Also, the map from $D_i$ to $\overline{q_{i-1}q_{i}}$
must be totally ramified at $q_{i-1}$ and $q_{i}$; 
otherwise, $\deg([D_i]) \ge 3$ by \lemref{lem1}. 
Since $D_i$ is rational, the map
from $D_i$ to $\overline{q_{i-1}q_{i}}$ is only ramified at $q_{i-1}$ and
$q_{i}$.

Let $\beta = m_2(q_0)$. Obviously, $\Gamma_2$ does not pass through $q_1, q_2,
..., q_{2l}$ and $\Gamma_1$ does not pass through $q_0, q_1, ...,
q_{2l-1}$; otherwise, either $\deg([\Gamma_2])\ge 2$ or
$\deg([\Gamma_1])\ge 2$. So
$m_2(q_1) = m_2(q_2) = ... = m_2(q_{2l}) = m_1(q_0) = m_1(q_1) = ... =
m_1(q_{2l-1}) = 0$.
Hence each $\overline{q_{i-1}q_{i}}$ have multiplicity exactly
$\beta$ in $\Upsilon_0$ for $i=1,2,..., 2l$. 
Therefore, the map from $D_i$ to 
$\overline{q_{i-1}q_{i}}$ has degree $\beta$.
This also implies that $m_1(q_{2l}) = \beta$.

Hence $\Gamma_1$ and $\Gamma_2$ are joined by a chain of curves $D_1\cup
D_2 \cup ...\cup D_{2l}$ on
$\wt{\Upsilon}_0$ whose images are contained in $\overline{q_0q_1}\cup
\overline{q_1 q_2}\cup ... \cup \overline{q_{2l-1} q_{2l}}$.
Therefore, there is only
one Type I chain in $\Upsilon_0$; otherwise, $\Gamma_1$ and
$\Gamma_2$ will be joined by two different chains of curves
on $\wt{\Upsilon}_0$.
\end{proof}

To study the behavior of a family $\Upsilon$ of curves near a subscheme
$S\subset \Upsilon_0$, it is usually very convenient to take an analytic
neighborhood $U$ of $S$ in $\Upsilon$ and study $U$ instead of $\Upsilon$
(or alternatively, study the formal completion of $\Upsilon$ along
$S$). Notice that even if $\Upsilon$ is irreducible, $U$ could be reducible
after a base change. We call a component of $U$ which is irreducible
under any base changes ``a locally irreducible component of $\Upsilon$
around $S$''. And if $U$ is irreducible under any base changes, we call
$\Upsilon$ is locally irreducible around $S$.


Here is a trivial remark. Let $S\subset \Upsilon_0$ be a closed subscheme
of $\Upsilon_0$ and $\wt{S}$ be the total transform of $S$ under the map
$\wt{\Upsilon} \to \Upsilon$. Then $\Upsilon$ is locally irreducible around
$S$ if and only if $\wt{S}$ is connected.

Let us write a Type II chain in the form $\cup_{i=-l}^{m-1} 
\overline{q_i q_{i+1}}$, where $q_0 \in \{p_1, p_2, ..., p_{16}\}$.

\begin{prop}\label{prop2a}
Let $\cup_{i=-l}^{m-1} \overline{q_i q_{i+1}}$ be a Type II chain in
$\Upsilon_0$ with $q_0 \in  \{p_1, p_2, ..., p_{16}\}$. 
Without the loss of generality, assume that $q_0=p_1$ and
$\overline{q_0q_1}\subset R_1$.

Let $Y$ be a locally irreducible component of $\Upsilon$ around
$\cup_{i=-l}^{m-1} \overline{q_i q_{i+1}}$.

Then

\begin{enumerate}
\item $\mu(\overline{q_{-l} q_{-l+1}})\le \mu(\overline{q_{-l+1} q_{-l+2}})
\le ...\le \mu(\overline{q_{-1} q_0})$ and
$\mu(\overline{q_{0} q_{1}})\ge \mu(\overline{q_{1} q_{2}})
\ge ...\ge \mu(\overline{q_{m-1} q_{m}})$;

\item $m_1(q_0), m_2(q_0)\le 1$;

\item for each $q_j$ ($j\ne 0$), either $m_1(q_j) = 0$ or $m_2(q_j) = 0$;

\item $Y_0$ is one of the following

\begin{enumerate}
\item $Y_0 = \Delta(w)$
where $w$ is the intersection between
$\Gamma_1\cup \Gamma_2$ and $\overline{q_i q_{i+1}}$ other than the points
$q_i$ and $q_{i+1}$ for some $i$ and $\Delta(w)\subset \Gamma_1\cup \Gamma_2$
is a disk centered at $w$ regarded as an analytic neighborhood of $w$ on
$\Gamma_1$ or $\Gamma_2$;

\item $Y_0 = \Delta(q_0)$ where $\Delta(q_0)\subset \Gamma_1\cup \Gamma_2$
is a disk centered at $q_0$ regarded as an analytic neighborhood of $q_0$ on
$\Gamma_1$ or $\Gamma_2$;

\item $(Y_0)_\red = \cup_{i=0}^{n-1} \overline{q_i q_{i+1}} 
\cup \Delta(q_n)$ for some $n > 0$,
where $\Delta(q_n)$ is a disk centered at $q_n$ regarded as an analytic
neighborhood of $q_n$ on $\Gamma_1$ or $\Gamma_2$,
$\Delta(q_n)$ and $\overline{q_{n-1} q_{n}}$ lie on the
different $R_j$'s and $\overline{q_i q_{i+1}}$ has multiplicity $E\cdot
\Delta(q_n)$ in $Y_0$;

\item $(Y_0)_\red = \cup_{i=-n}^{-1}
\overline{q_i q_{i+1}}\cup \Delta(q_{-n})$ for some $n > 0$,
where $\Delta(q_{-n})$ is a disk centered at $q_{-n}$ regarded as an analytic
neighborhood of $q_{-n}$ on $\Gamma_1$ or $\Gamma_2$,
$\Delta(q_{-n})$ and $\overline{q_{-n} q_{-n+1}}$ lie on the different
$R_j$'s and $\overline{q_i q_{i+1}}$ has multiplicity $E\cdot \Delta(q_{-n})$
in $Y_0$.
\end{enumerate}
\end{enumerate}
\end{prop}

\begin{rem}
The statements in \propref{prop2a} may need some further explanation. Let
$U$ be an analytic neighborhood of $\cup_{i=-l}^{m-1} \overline{q_i
q_{i+1}}$ in $\Upsilon$. Then the central fiber $U_0$ of $U$ consists of
$\cup_{i=-l}^{m-1} \overline{q_i q_{i+1}}$ plus a few disks which are
``pieces'' of $\Gamma_1$ and $\Gamma_2$. If $\Gamma_i$ meets $\overline{q_j
q_{j+1}}$ at $w \ne q_j, q_{j+1}$, there is a disk $\Delta(w)\subset
\Gamma_i$ centered at $w$ on $U_0$; if $\Gamma_i$ passes through $q_j$,
there is a disk $\Delta(q_j)\subset \Gamma_i$ centered at $q_j$ on $U_0$.
Let $\wt{U} = U\times_{\Upsilon} \wt{\Upsilon}$ be the nodal reduction of
$U$. Each locally irreducible component $Y$ of $U$ corresponds to a connected
component of the dual graph of $\wt{U}_0$ and vice versa.

First, locally at $w = \Gamma_i\cap \overline{q_j q_{j+1}}$ with $w\ne q_j,
q_{j+1}$, our statements about $Y_0$ show that $\Gamma_i$ ($\Gamma_i\supset
\Delta(w)$) is not joined to any component dominating $\overline{q_j
q_{j+1}}$ by a chain of curves contracting to $w$ on $\wt{\Upsilon}_0$. In
a more intuitive language, $\Gamma_i$ is ``separated'' from $\overline{q_j
q_{j+1}}$ at $w$ after the nodal reduction. 


Second, locally at $q_0$, our statements about $Y_0$ show that all branches
of $\Upsilon_0$ at $q_0$
are separated from each other after the nodal reduction. By that we mean
among the components dominating $\Gamma_1$, $\Gamma_2$, $\overline{q_{-1}
q_0}$ or $\overline{q_0 q_1}$ on $\wt{\Upsilon}_0$, no two are joined by a
chain of curves contracting to $q_0$. 


Third, locally at $q_n$ ($n\ne 0$), suppose that $\Gamma_i$ passes through
$q_n$ and $Y$ is the locally irreducible component of $U$ containing the
disk $\Delta(q_n)$. Then our statements about $Y_0$ show that
$\overline{q_n q_{n+1}}$ (if $n < 0$) or
$\overline{q_{n-1} q_n}$ (if $n > 0$) must lie on $R_{3-i}$. This
also implies that $\Gamma_1$ and $\Gamma_2$ cannot both pass through $q_n$
for any $n\ne 0$. 
Let $\wt{Y} = Y\times_{\Upsilon} \wt{\Upsilon}$. Then $\wt{Y}_0$ is the
union of the components of $\wt{U}_0$ which form a connected component of
the dual graph of $\wt{U}_0$ corresponding to $Y$. By \propref{prop2a},
$\wt{Y}_0$ consists of $\Delta(q_n)$ and curves over $\cup_{i=0}^{n+1}
\overline{q_{i-1} q_{i}}$ (if $n < 0$) or $\cup_{i=1}^{n} \overline{q_{i-1}
q_{i}}$ (if $n > 0$) and each $q_{i-1} q_i$ is dominated by the components
of $\wt{Y}_0$ through maps whose total degree is $E\cdot
\Delta(q_n)$. Later, we will prove in \thmref{thm2} that $E\cdot \Delta(q_n)
= 1$, i.e., if $\Gamma_i$ passes through $q_n$, $\Gamma_i$ and $E$ must
meet transversely at $q_n$.
\end{rem}

\begin{proof}[Proof of \propref{prop2a}]
Let us first prove the statements concerning $Y_0$. The rest will follow
more or less immediately.

Obviously, there is at most one disk $\Delta\subset Y_0$ with
$\Delta\subset\Gamma_1\cup\Gamma_2$ by \propref{prop1}; otherwise,
either $\Gamma_i$ is joined to itself by a chain of curves over 
$\cup_{j=-l}^{m-1} \overline{q_j q_{j+1}}$ for some $i$ or $\Gamma_1$ and
$\Gamma_2$ are joined by a chain of curves over $\cup_{j=-l}^{m-1}
\overline{q_j q_{j+1}}$, but by \propref{prop1}, $\Gamma_1$ and $\Gamma_2$
are already joined by the curves over a Type I chain.

If $\overline{q_i q_{i+1}}\subset Y_0$ for some $i \ge 0$, by
\lemref{lem1}, either $\overline{q_{i+1} q_{i+2}}\subset Y_0$ or there is a
disk $\Delta(q_{i+1})\subset \Gamma_1\cup \Gamma_2$ such that
$\Delta(q_{i+1})\subset Y_0$ and
$\Delta(q_{i+1})$ and $\overline{q_i q_{i+1}}$ lie on the different
$R_j$'s. If $\overline{q_{i+1} q_{i+2}}\not\subset Y_0$, we are done;
otherwise, we apply the same argument to $\overline{q_{i+1} q_{i+2}}$
again. And eventually, this sequence of curves will end up at some disk
$\Delta(q_n)\subset Y_0$ for some $n > i$.
Obviously, $\overline{q_n q_{n+1}}\not\subset
Y_0$; otherwise, we may continue to apply the above argument to show that
there exists another disk $\Delta(q_{n'})\subset Y_0$ for some $n' > n$.
On the other hand, since $\overline{q_i q_{i+1}}\subset Y_0$, we have
either $\overline{q_{i-1} q_i}\subset Y_0$ or $\Delta(q_{i-1})\subset
Y_0$ by \lemref{lem1}. Since $\Delta(q_n)\subset Y_0$, we necessarily have
$\overline{q_{i-1} q_i}\subset Y_0$. Apply the same argument to
$\overline{q_{i-1} q_i}$ and we obtain that $\overline{q_{i-2}
q_{i-1}}\subset Y_0$. So eventually, we have $\overline{q_0 q_1},
\overline{q_1 q_2}, ..., \overline{q_{n-1} q_n}, \Delta(q_n)\subset
Y_0$. It is impossible that $\overline{q_{-1} q_0}\subset Y_0$; otherwise,
we may apply the same line of argument to show subsequently that
$\overline{q_{-1} q_{0}}, \overline{q_{-2} q_{-1}}, ..., \overline{q_{-l}
q_{-l + 1}} \subset Y_0$ and eventually $\Delta(q_{-l})\subset Y_0$. So
$\overline{q_{-1} q_0} \not\subset Y_0$, $\overline{q_n q_{n+1}}\not\subset
Y_0$ and hence $(Y_0)_\red = \cup_{j=0}^{n-1} \overline{q_j q_{j+1}} \cup
\Delta(q_n)$. By \lemref{lem1}, all $\overline{q_j q_{j+1}}$ in this
sequence has the same multiplicity $\Delta(q_n) \cdot E$ in $Y_0$. So if
$\overline{q_i q_{i+1}}\subset Y_0$ for some $i \ge 0$, we will
necessarily end up in case (c).

The same argument shows if $\overline{q_{i} q_{i+1}}\subset Y_0$ for
some $i < 0$, we will end up in case (d).

It is obvious that if there is no $\overline{q_{i} q_{i+1}} \subset Y_0$,
we will necessarily end up in case (a) or (b).

So if $\overline{q_{i} q_{i+1}}\subset Y_0$ for $i \ge 0$, we necessarily
have $\overline{q_j q_{j+1}}\subset Y_0$ for any $0\le j \le i$ and
$\overline{q_j q_{j+1}}$ has the same multiplicity in $Y_0$ as
$\overline{q_{i} q_{i+1}}$. Therefore, $\mu(\overline{q_{0} q_{1}})\ge
\mu(\overline{q_{1} q_{2}}) \ge ...\ge \mu(\overline{q_{m-1}
q_{m}})$. Similarly, $\mu(\overline{q_{-l} q_{-l+1}})\le
\mu(\overline{q_{-l+1} q_{-l+2}}) \le ...\le \mu(\overline{q_{-1} q_0})$.

It follows from \coref{cor3} that $m_1(q_0) \le 1$ and $m_2(q_0) \le 1$.
Otherwise, suppose that $m_1(q_0) > 1$. Let $Y$ be the locally irreducible
component of $\Upsilon$ around $\cup_{i=-l}^{m-1} \overline{q_i
q_{i+1}}$ such that $\Delta(q_0)\subset \Gamma_1$ and 
$\Delta(q_0)\subset Y_0$. By \coref{cor3} and the fact that $\Delta(q_0)$
is the only ``piece'' of $\Gamma_1\cup \Gamma_2$ in $Y_0$, we must have
$\overline{q_{-1} q_0}\subset Y_0$, which contradicts our results on
possible $Y_0$'s.

Finally, it is impossible that $m_1(q_n) > 0$ and $m_2(q_n) > 0$ for
$n \ne 0$. Otherwise, suppose that both $\Gamma_1$ and $\Gamma_2$ pass
through $q_n$ for some $n > 0$. Since $\Gamma_1$ passes through $q_n$,
there exists a locally irreducible component $Y$ of $\Upsilon$ around
$\cup_{i=-l}^{m-1} \overline{q_i q_{i+1}}$ such that
$(Y_0)_\red = \cup_{i=0}^{n-1} \overline{q_i q_{i+1}} \cup \Delta(q_n)$
where $\Delta(q_n)\subset \Gamma_1$ and $\Delta(q_n)$ and
$\overline{q_{n-1} q_n}$ lie on the different $R_j$'s. So
$\overline{q_{n-1} q_n}$ lies on $R_2$. The same argument shows that
$\overline{q_{n-1} q_n}$ lies on $R_1$ since $\Gamma_2$ passes through
$q_n$. Contradiction.
\end{proof}

As \propref{prop1} and \ref{prop2a} for Type I and Type II chains, we
have a similar statement for Type III chains. However, we do not really
need it in our proof. So we will state the proposition without a proof.
Interested readers could follow the same line of argument as in the
\propref{prop1} and \ref{prop2a} and give a proof themselves.

\begin{prop}\label{prop2}
A type III chain in $\Upsilon_0$ containing the point $r_{ij}$ must also
contain the curve $\overline{r_{ij}}$.

Let $\overline{r_{ij}}\cup \overline{r_{ij} q_1}\cup \overline{q_1q_2}\cup
... \cup \overline{q_{m-1} q_m}$ be a Type III chain in
$\Upsilon_0$. Assume that $r_{ij} = r_{11}$ without the loss of
generality.

Let $Y$ be a locally irreducible component of $\Upsilon$ around
$\overline{r_{11}}\cup \overline{r_{11} q_1}\cup \overline{q_1q_2}\cup
... \cup \overline{q_{m-1} q_m}$.  Let $q_0 = r_{11}$.

Then

\begin{enumerate}
\item $2\mu(\overline{r_{11}})\ge \mu(\overline{r_{11} q_1}) \ge
\mu(\overline{q_1q_2})\ge ...\ge \mu(\overline{q_{m-1} q_m})$;
\item $m_i(q_j)$ are even for all $i$ and $j$;
\item $m_i(q_j) = 0$ if $i+j$ is odd;
\item $Y_0$ is one of the following

\begin{enumerate}
\item $Y_0 = \Delta(w)$ where $w$ is one of the intersections between
$\Gamma_1\cup \Gamma_2$ and $\overline{r_{ij}}\cup \overline{r_{ij}
q_1}\cup \overline{q_1q_2}\cup ... \cup \overline{q_{m-1} q_m}$
other than the points $r_{11}, q_1, q_2, ..., q_m$ and
$\Delta(w)$ is an analytic neighborhood of $w$ on $\Gamma_1$ or $\Gamma_2$;

\item $(Y_0)_\red = \overline{r_{11}}\cup\overline{r_{11} q_1}\cup
\overline{q_1 q_2} \cup ...\cup \overline{q_{n-1} q_n} \cup \Delta(q_n)$
for some $n \ge 0$, where $\Delta(q_n)$ is an analytic neighborhood of
$q_n$ on $\Gamma_1$ or $\Gamma_2$,
$\Delta(q_n)$ and $\overline{q_{n-1} q_n}$ lie on the
different $R_j$'s, $\overline{q_i q_{i+1}}$ has multiplicity
$\Delta(q_n)\cdot E$ in $Y_0$ and $\overline{r_{11}}$ has multiplicity
$(\Delta(q_n)\cdot E)/2$ in $Y_0$.
\end{enumerate}
\end{enumerate}
\end{prop}

\section{Deformation around a Type I Chain}
\label{sec2}

In this section,
we will study the behavior of $\Upsilon_t$ at the neighborhood of a Type I
chain in order to show that $\Upsilon_t$ has only
nodes as singularities in the neighborhood of a Type I chain.

For a one-parameter family of curves $S$ over $\Delta$ and
a reduced subscheme $B\subset S_0$, we use the notation
$\delta(S_t, B)$ to denote the total $\delta$-invariant 
of the general fiber $S_t$ in a neighborhood of $B$.

The main theorem of this section is

\begin{thm}\label{thm1}
Let $\cup_{i=1}^{2l} \overline{q_{i-1} q_i}$ be the Type I chain in
$\Upsilon_0$. Assume that $\overline{q_0 q_1}\subset R_1$. Then

\begin{enumerate}
\item $\delta(\Upsilon_t, \cup_{i=1}^{2l} \overline{q_{i-1} q_i}) =
(2l+1)\beta - 1$, where $\beta = m_2(q_0) = m_1(q_{2l})$;
\item $\Upsilon_t$ has exactly $\beta$ nodes in the neighborhood of 
each point $w_j$ for $j = 1, 2,...,2l$, where $w_j$ is the intersection
between $\Gamma_1\cup \Gamma_2$ and $\overline{q_{j-1}
q_j}$ other than the points $q_{j-1}$ and $q_j$; 
\item $\Upsilon_t$ has exactly $\beta - 1$ nodes in the neighborhood of
$q_l$.
\end{enumerate}
\end{thm}

\begin{rem}
Notice that
\[
\begin{split}
\delta(\Upsilon_t, \cup_{i=1}^{2l} \overline{q_{i-1} q_i})
&= (2l+1)\beta - 1\\
&=  \sum_{i =0}^{2l} \intsc(\Sigma_1, E; q_i) -1\\
&= \sum_{i =0}^{2l} \intsc(\Sigma_2, E; q_i) - 1,
\end{split}
\]
where $\intsc(A, B; C)$ denotes the local intersection multiplicities
between $A$ and $B$ along $C$.
\end{rem}

The main obstacle here is that $\Upsilon_0$ is nonreduced along
$\cup_{i=1}^{2l} \overline{q_{i-1} q_i}$ if $\beta > 1$.
To study the deformation of a
nonreduced curve, we introduce a method called ``patching
technique''. Actually, this is a very commonly used method in deformation
theory. To study the
deformation of a projective variety or compact complex manifold $M$, 
we cover $M$ by affine or analytic open sets, study
the deformation of each piece separately and then ``patch'' these
deformations together. We will show how this can be done for a nonreduced
curve. 

\subsection{Deformation of a Nonreduced Planary Curve}
 
\begin{defn}
Let $C$ be a nonreduced scheme and $C_\red$ be its reduced subscheme. Let
$\I$ be the ideal sheaf of $C_\red$ in $C$. Then $\n_{C_\red} = 
\hom(\I/\I^2, \CO_{C_\red})$ is a coherent
sheaf over $C_\red$, which we will call the 
``intrinsic normal sheaf'' of $C$.
\end{defn}

\begin{defn}
We call a nonreduced curve $C$ ``planary'' if it can be locally 
embedded to the 2-dimensional polydisk $\Delta_{xy}^2$ everywhere. We say
$C$ has irreducible support if $C_\red$ is irreducible.

If $C$ is a planary nonreduced curve with irreducible support, we can
cover $C$ with analytic open sets $U_\alpha$ such that each 
\[
U_\alpha\isom \Spec \BC[[x, y]]/(f^m(x, y)),
\]
where $f(x, y) = 0$ defines
$C_\red$ in $U_\alpha$ and $m$ is called the multiplicity of $C$.
And it is obvious in this case that the intrinsic normal sheaf
$\n_{C_\red}$ of $C$ is a line bundle over $C_\red$.

Obviously, $C$ is planary if it lies on a smooth surface.
If $C$ lies on a smooth surface $S$ and has irreducible support,
its intrinsic normal sheaf
is simply the normal sheaf $\n_{C_\red/S}$ of $C_\red\subset S$.
\end{defn}

Now let $C$ be a nonreduced planary curve with irreducible support and 
assume that $C_\red$ is smooth. Then $C$ is locally isomorphic to
$\Spec\BC[[x, y]]/(y^m)$, where $m$ is the multiplicity of $C$.
Let $\C$ be a one-parameter family of
curves over disk $\Delta$ whose central fiber is $C$. Cover $C$ with open
sets $U_\alpha$ of $\C$ such that $U_\alpha\isom \Spec \BC[[x, 
y, t]]/
(F_\alpha(x, y, t))$ where $F_\alpha(x, y, t)$ can be put into the form
\[
F_\alpha(x, y, t) = y^m + \sum_{\genfrac{}{}{0pt}{}{i>0}{0\le j\le m-2}} 
f_{\alpha ij}(x) t^i y^j.
\]
For each $\alpha$, we let
\[
\gamma(\alpha) = \min\left\{\frac{i}{m - j}: f_{\alpha ij}(x) 
\not\equiv 0\right\}.
\]
Alternatively, we may define $\gamma(\alpha)$ as the largest number such
that
\[
t^{m \gamma(\alpha)} \big | F_\alpha(x, t^{\gamma(\alpha)} y, t).
\]

After a base change, we may assume that $\gamma(\alpha)\in \BZ$ for each
$\alpha$. Then we may put $F_\alpha(x, y, t)$ into the form
\[
F_\alpha(x, y, t) = y^m + \sum_{j = 0}^{m-2} t^{(m - j) \gamma(\alpha)}
\phi_{\alpha j} (x, t) y^j
\]
where $\phi_0(x, 0), \phi_1(x, 0), ..., \phi_{m-2}(x, 0)$ are not all zeros
by the choice of $\gamma(\alpha)$.

Next, we will ``patch'' $U_\alpha$ together. We need the following
``patching'' lemma.

\begin{lem}[Patching Lemma]\label{lem3}
Let $S_1, S_2\subset \Delta_{xy}^2\times \Delta_t$ be two families of
curves over disk $\Delta_t$. Suppose that $S_i$ ($i = 1,2$)
is cut out by
\[
y^m + \sum_{j=0}^{m-2} 
t^{(m - j) \gamma_i} \phi_{ij}(x, t) y^j
= 0
\]
in $\Delta_{xy}^2\times \Delta_t$, 
where $\gamma_i\in \BZ$, $\gamma_i > 0$,
and $\phi_{i, 0}(x, 0), \phi_{i, 1}(x, 0), ..., 
\phi_{i, m-2}(x, 0)$ are not all zero.
If there is an isomorphism between $S_1$ and $S_2$ which preserves the base
and induces the identity map on the central fibers,
then $\gamma_1 = \gamma_2$ and 
$\phi_{1j}(x, 0) = \phi_{2j}(x, 0)$ for $j=0,1,...,m-2$.
\end{lem}

\lemref{lem3} is more or less obvious and we will leave its proof to the
readers. It follows from \lemref{lem3} that all $\gamma(\alpha)$ are equal
and for each $j$, $\{ \phi_{\alpha j}(x, 0)\}$ defines a global section of
the line bundle $\n_{C_\red}^{\tensor (m - j)}$, i.e.,
there exists $s\in H^0(\n_{C_\red}^{\tensor (m - j)})$ such that 
$s |_{U_\alpha} = \phi_{\alpha j}(x, 0)$, 
where $\n_{C_\red}$ is the intrinsic normal sheaf of $C$.

Let $\gamma = \gamma(\alpha)$. Then $\{y = t^\gamma = 0\}$ is a well-defined
closed
subscheme of $\C$ supported at the central fiber. Let $\wt{\C}$ be
the blowup of $\C$ along the subscheme $\{ y = t^\gamma = 0\}$. It is not
hard to see that $\wt{\C}_0$ is a curve in the linear series $|\CO_\P(m)|$ on
$\P = \P\left(\CO_{C_\red} \oplus \n_{C_\red}\right)$. Now $\wt{\C}_0$ is
``less'' nonreduced than $C$, i.e., each component of $\wt{\C}_0$ has
multiplicity strictly less than $m$ in $\wt{\C}_0$ due to our choice of
$\gamma$. This makes $\wt{\C}$ easier to investigate than $\C$.

\subsection{Outline of the Proof of \thmref{thm1}}

Let $I_j$ be the nonreduced component of $\Upsilon_0$ supported on
$\overline{q_{j-1} q_j}$ with multiplicity $\beta$.
Obviously, $I_j$ is a nonreduced planary curve with trivial intrinsic
normal sheaf.

The proof of \thmref{thm1} is carried out by repeatedly
blowing up $\Upsilon$ along $\cup I_j$. To be precise, we
will construct a sequence of families
\begin{equation}\label{e3b}
\Upsilon^{(l)} \to 
\Upsilon^{(l-1)} \to ... \to \Upsilon^{(0)} = \Upsilon
\end{equation}
where

\begin{enumerate}
\item the morphisms $\Upsilon^{(j)}\to \Upsilon^{(j-1)}$
are isomorphisms on the general fibers for $j=1,2,...,l$;
\item on the central fiber, $\Upsilon_0^{(j)}$ contains
$\wt{I}_{1}\cup \wt{I}_2 \cup...\cup \wt{I}_{j-1}\cup \wt{I}_j \cup I_{j+1}
\cup ... \cup I_{2l - j} \cup \wt{I}_{{2l-j+1}} \cup \wt{I}_{{2l-j+2}}
\cup ... \cup \wt{I}_{{2l-1}}\cup \wt{I}_{2l}$ for $j=0,1,2,...,l$,
where $\wt{I}_{1}, \wt{I}_{2}, ...,\wt{I}_{2l}$ are the curves dominating 
$(I_1)_\red = \overline{q_0 q_1}, (I_2)_\red = \overline{q_1q_2}, ..., 
(I_{2l})_\red = \overline{q_{2l-1}q_{2l}}$,
respectively; recall from \propref{prop1} that $\wt{I}_j$ dominates
$\overline{q_{j-1}q_{j}}$ with a degree $\beta$ map totally ramified at
$q_{j-1}$ and $q_j$;
\item let $w_{j1}, w_{j2}, ..., w_{j\beta}$ be the points over $w_i$ in the
map $\wt{I}_j \to (I_j)_\red$ for $j=1,2,...,2l$; there is a contractible 
curve $H_j$ meeting $\wt{I}_j$ transversely at
$w_{j1}, w_{j2}, ..., w_{j\beta}$ and meeting $\Gamma_1$ or $\Gamma_2$
transversely at another point;
\item
each pair of curves $\wt{I}_j$ and $\wt{I}_{j+1}$ meet
transversely at a point over $q_j$, which we still denote by $q_j$, except
that $\wt{I}_{l}$ and $\wt{I}_{{l+1}}$ meet at the point $q_l$ with
multiplicity $\beta$ on $\Upsilon^{(l)}$; more
precisely, $\Upsilon^{(l)}\isom \Spec \BC[[y, z, t]]/(y(y + z^\beta +
O(t)) - t^\alpha)$ in an analytic neighborhood of $q_l$, where
we use the notation $O(f_1, f_2, ..., f_n)$ to denote an element generated
by $f_1, f_2, ..., f_n$ in some ring (it should be clear from the context
which ring we are talking about).
\end{enumerate}

It is obvious from the above description that $\Upsilon$ has exactly
$\beta$ nodes in the neighborhoods of $w_j$. The statement
that $\Upsilon_t$ has exactly $\beta - 1$ nodes in the neighborhood of
$q_l$ follows from a theorem of L. Caporaso and J. Harris on the
deformation of tacnodes \cite[Lemma 4.1]{CH} and the fact that
the curves $\wt{I}_l$ and $\wt{I}_{l+1}$ meet transversely on the
nodal reduction $\wt{\Upsilon}$. 

So the proof of \thmref{thm1} boils down to the construction of the blowup
sequence \eqref{e3b}.

\subsection{Construction of the Blowup Sequence}

We will do the construction inductively. For that purpose, we will work
on an analytic neighborhood of $\Upsilon$ around $(\cup I_j)_\red$, which
we will still refer to by $\Upsilon$ and can be described as follows.

First, $\Upsilon$ is a one-parameter family of curves over disk $\Delta_t$ with
irreducible general fibers $\Upsilon_t$. Second, $\Upsilon_0 = I_1 \cup I_2
\cup ...\cup I_{2l} \cup \Delta(q_0) \cup \Delta(q_{2l}) \cup \Delta(w_1) \cup
\Delta(w_2) \cup... \cup \Delta(w_{2l})$, where $I_j\isom \P^1\times \Spec
\BC[z]/(z^\beta)$ and $\Delta(q_0)$, $\Delta(q_{2l})$ and $\Delta(w_i)$ are
disks centered at $q_0$, $q_{2l}$ and $w_i$ regarded as the analytic
neighborhoods of $q_0$, $q_{2l}$ and $w_i$ on $\Gamma_1$ or $\Gamma_2$.
Finally, these curves are ``patched'' up in the following way:

\begin{enumerate}
\item $q_{j} = (I_j)_\red\cap (I_{j+1})_\red$ 
and
\begin{equation}\label{e4}
\Upsilon\isom \Spec \BC[[x, y, z, t]]/(xy - t^\alpha, z^\beta + O(t)) 
\end{equation}
in the neighborhoods of $q_{j}$ for $j = 1, ..., 2l-1$, 
\item $q_0 = \Delta(q_0)\cap (I_1)_\red$, 
$q_{2l} = \Delta(q_{2l}) \cap (I_{2l})_\red$ and
\begin{equation}\label{e5}
\Upsilon\isom \Spec \BC[[x, y, z, t]]/(xy - t^\alpha, y - z^\beta + O(t))
\end{equation}
in the neighborhoods of $q_0$ and $q_{2l}$;
\item $w_j = \Delta(w_j)\cap (I_j)_\red$ and
\begin{equation}\label{e6}
\Upsilon\isom \Spec \BC[[x, z, t]]/(x z^\beta + O(t))
\end{equation}
in the neighborhoods of $w_j$ for $j=1,2,...,2l$.
\end{enumerate}


We also know that the nodal reduction $\wt{\Upsilon}$ of $\Upsilon$ has
the properties described in \propref{prop1} with $\Delta(q_0)$,
$\Delta(q_{2l})$ and $\Delta(w_i)$ regarded as ``pieces'' of $\Gamma_1$ or
$\Gamma_2$.


By \eqref{e6}, $\Upsilon$ can be locally embedded into $\Delta_{xz}^2\times
\Delta_t$ at $w_j$. After applying some automorphism of
$\Delta_{xz}^2\times \Delta_t$ which preserves the base and the central fiber,
we may put its local defining equation at $w_j$ into the form 
\begin{equation}\label{e7}
x z^\beta + t^{m_{1j}} z^{\beta-1} + \sum_{i=2}^\beta 
t^{m_{ij}} f_{ij}(x, t) z^{\beta -i} = 0
\end{equation}
where $f_{ij}(x, t)\in \BC[[x, t]]$, $f_{ij}(x, 0) \ne 0$ and we put
$m_{ij} = \infty$ if the corresponding term $z^{\beta-1}$ or $f_{ij}(x, t)
z^{\beta -i}$ does not appear in the defining equation. Let
\[
\gamma_j = \min\left\{\frac{m_{ij}}{i}: 1\le i\le \beta\right\} \text{ and
}
\gamma = \min\left\{\frac{\alpha}{\beta}, \gamma_1, \gamma_2, ...,
\gamma_{2l}
\right\}.
\]
Obviously, we may assume that $\gamma, \gamma_j\in\BZ$ for $j =1,2,..., 2l$
after a base change.

We can write \eqref{e7} as
\begin{equation}\label{e8}
x z^\beta + \sum_{i=1}^\beta t^{i\gamma} F_{ij}(x, t) z^{\beta -i} = 0
\end{equation}
where $F_{ij}(x, t)\in \BC[[x, t]]$ and $F_{1j}(x, t) = F_j(t)\in \BC[[t]]$. 
Notice that $F_{ij}(x, 0) = 0$ if $\gamma_j>\gamma$.

\begin{claim}\label{claim1}
We claim that
\begin{enumerate}
\item $F_{ij}(0, 0) = 0$ for $i = 1, 2, ..., \beta$ and $j = 1,2,...,2l$
and especially, $F_{1j}(x, 0) = F_j(0) = 0$;
\item $x^{-1} F_{ij}(x, 0)$ extends to a meromorphic function $G_{ij}$ on
$(I_j)_\red$;
\item each $G_{ij}$ is holomorphic everywhere on $(I_j)_\red$ except that
$G_{\beta, 1}$ and $G_{\beta, 2l}$ have a simple pole at $q_0$
and $q_{2l}$, respectively, if $\alpha = \beta\gamma$;
\item at each $q_j$ for $j= 1,2,...,2l-1$, $G_{i,j}(q_j) = G_{i,j+1}(q_j)$.
\end{enumerate}
\end{claim}

The first statement follows from the assumption that $\Delta(w_j)$
and $\wt{I}_j$ are disjoint on $\wt{\Upsilon}_0$.

Let $U$ be an analytic neighborhood of $w_j$ on $\Upsilon$. Since $\Delta(w_j)$
and $\wt{I}_j$ are disjoint on $\wt{\Upsilon}_0$, $U$ is reducible and $U =
U^{(1)}\cup U^{(2)}$ after a base change
where the central fibers of $U^{(1)}$ and $U^{(2)}$
are given by $x = t= 0$ and $z^\beta = t = 0$,
respectively. Correspondingly, we may factor the LHS of \eqref{e8} in the
following way
\begin{equation}\label{e9}
x z^\beta + \sum_{i=1}^\beta t^{i\gamma} F_{ij}(x, t) z^{\beta -i}
= (x + O(t)) (z^\beta + O(t)).
\end{equation}
It follows immediately from \eqref{e9} that $F_{ij}(0, 0) = 0$ which also
implies $F_{1j}(x, 0) = F_{j}(0) = 0$.

It follows from \lemref{lem3} that
$x^{-1} F_{ij}(x, 0)$ can be analytically extended to a section in
\[
\Gamma\left((I_j)_\red \backslash \{q_{j - 1}, q_j\}, \n_j^{\tensor
i}\right)
\]
where $\n_j$ is the intrinsic normal sheaf of $I_j$. Of course, $\n_j$ is
trivial. So each $G_{ij} = x^{-1} F_{ij}(x, 0)$ is a meromorphic function
on $(I_j)_\red\isom \P^1$, which is holomorphic everywhere except at
the points $q_0, q_1, q_2,..., q_{2l}$.

Actually, $G_{ij}$ can be extended over $q_1, q_2, ..., q_{2l-1}$.
By \eqref{e4}, $\Upsilon$ can be locally embedded into $\Spec \BC[[x,y,z,
t]]/(xy - t^\alpha)$ at $q_1, q_2,...,q_{2l-1}$ and its local defining
equation at $q_j$ can be put into the form
\begin{equation}\label{e9a}
z^\beta + \sum_{i=2}^\beta t^{i \delta_j} \phi_{ij}(x,y, t) z^{\beta - i}
= 0
\end{equation}
where $\delta_j > 0$, $\phi_{ij}(x, y, t)\in\BC[[x, y, t]]$,
for $i=2,3,...,\beta$ and $j = 1,2,...,2l-1$ and $\phi_{2j}(x, y, 0),
\phi_{3j}(x, y, 0), ..., \phi_{\beta j}(x, y, 0)$ do not all lie in the
ideal $(xy)\subset \BC[[x, y]]$.

By comparing \eqref{e9a} with \eqref{e8} and applying \lemref{lem3}, we
have $\delta_j = \min(\gamma_j, \gamma_{j+1})$. So

\begin{enumerate}
\item $\delta_j\ge \gamma$;
\item $\gamma_j > \gamma$ and $\gamma_{j+1} > \gamma$ if $\delta_j > \gamma$.
\end{enumerate}

Hence if $\delta_j > \gamma$, we have
$G_{i, j} = G_{i, j+1} = 0$; otherwise, if $\delta_j = \gamma$, it is not
hard to see that $G_{i, j}(q_j) = G_{i, j+1} (q_j) = \phi_{ij}(0,
0, 0)$. Hence $G_{ij}$ are holomorphic at $q_j$ and $G_{i, j} (q_j) =
G_{i, j+1}(q_j)$ for $j=1,2,...,2l-1$.

So $G_{ij}$'s are holomorphic everywhere except at $q_0$ and $q_{2l}$. Next,
we will try to find out what kind of singularities $G_{ij}$'s could
have at $q_0$ and $q_{2l}$.
 
By \eqref{e5}, $\Upsilon$ can be locally
embedded into $\Spec \BC[[x, z, t]]$ at $q_0$ and $q_{2l}$. We can put
its local defining equations at $q_j$ into the form
\begin{equation}\label{e10}
x\left(z^\beta + \sum_{i=2}^\beta t^{i\delta_j} \phi_{ij}(x, t)
z^{\beta-i} \right) = t^\alpha
\end{equation}
where $\delta_j > 0$, $\phi_{ij}(x, t)\in \BC[[x, t]]$
for $i = 2,3,...,\beta$ and $j =0, 2l$ and
$\phi_{2j}(x, 0), \phi_{3j}(x, 0), ..., 
\phi_{\beta j}(x, 0)$ are not all zero.

By comparing \eqref{e10} with \eqref{e8} and applying \lemref{lem3}, we
have

\begin{enumerate}
\item $\delta_j \ge \gamma$, $\alpha/\beta \ge\gamma_1$ and $\alpha/\beta
\ge \gamma_{2l}$;
\item $G_{i, 1}$ and $G_{i, 2l}$ are holomorphic at $q_0$ and $q_{2l}$,
respectively, for $i=2, ..., \beta - 1$;
\item $G_{\beta, 1}$ and $G_{\beta, 2l}$ 
have a simple pole at $q_0$ and $q_{2l}$, respectively, if
$\alpha = \beta \gamma$ and is holomorphic otherwise.
\end{enumerate}

So we have justified every statement in \claimref{claim1} and we are ready
to construct a blowup map $\Upsilon' \to \Upsilon$
where $\Upsilon_0'\supset \wt{I}_1\cup I_2 \cup ...\cup I_{2l-1} \cup
\wt{I}_{2l}$.

Notice that $z = t^\gamma = 0$ defines subschemes of $\Upsilon$ locally at
$w_j$ by \eqref{e8}.
We can extend these subschemes and patch them together to obtain a closed
subscheme of $\Upsilon$ due to our choice of the number $\gamma$.
This subscheme is
obviously supported on $(I_1)_\red\cup (I_2)_\red\cup ...\cup (I_{2l})_\red$.
The family $\Upsilon'$ is obtained by blowing up $\Upsilon$
along this closed subscheme, which is locally cut out by $z = t^\gamma = 0$.

First we claim that $\alpha = \beta \gamma$. If not, we necessarily have
$\alpha > \beta\gamma$. Hence each $G_{ij}$ is a holomorphic function over
$I_j$ and consequently a constant. By the equality $G_{i,j-1}(q_j) =
G_{i, j}(q_j)$, $G_{i1} = G_{i2} = ... = G_{i, 2l}$ and we let $G_i =
G_{i1}$.

Obviously, $\gamma_1 = \gamma_2 = ... =
\gamma_{2l} = \gamma$. Otherwise, if $\gamma_k > \gamma$ for some $k$, then
$G_{ik} = 0$, which implies $G_i = G_{ij} = 0$ for each $i$ and $j$.
And this is a contradiction to the choice of $\gamma$.

Let us examine the behavior of $\Upsilon'$ over the point $q_0$. Since
$\gamma_1 = \gamma < \alpha/\beta$, by \eqref{e10},
$\Upsilon_0'$ consists of the curve
\begin{equation}\label{e11}
x \left(z_1^\beta + \sum_{i=2}^\beta G_i z_1^{\beta-i}\right) = 0,
\end{equation}
where $z_1 = z/t^\gamma$. Obviously, the curve $x = 0$
dominates $\Delta(q_0)$ and the curve
\begin{equation}\label{e12}
z_1^\beta + \sum_{i=2}^\beta G_i z_1^{\beta-i} = 0
\end{equation}
maps to $I_1$. Since $G_2, G_3,..., G_\beta$ are not all zero, the LHS of
\eqref{e12} has at least two distinct roots. Hence there are at least two
different components over $I_1$ at $q_0$. This contradicts the fact that on
$\wt{\Upsilon}_0$ there is a single component $\wt{I}_1$ dominating $I_1$
with a map totally ramified at $q_0$. 

Therefore, we have $\alpha = \beta\gamma$. Then $G_{ij}$ are still
constants except that $G_{\beta, 1}$ and $G_{\beta, 2l}$ have a simple pole
at $q_0$ and $q_{2l}$, respectively. 

Let us examine the behavior of $\Upsilon'$ over the point $q_0$ again. Since
$\alpha = \beta\gamma$, by our previous analysis, $\gamma = \gamma_1 =
\alpha/\beta$. By \eqref{e10}, $\Upsilon_0'$ consists of the curve
\begin{equation}\label{e13}
x \left( 1+ \sum_{i=2}^{\beta-1} G_{i1} z_2^i\right) +
 (xG_{\beta,1})(q_0) z_2^\beta = z z_2 = 0
\end{equation}
where $z_2 = t^\gamma / z$ and $(xG_{\beta,1})(q_0)$ reads as the value of the
function $xG_{\beta,1}$ at point $q_0$.
Obviously, the curve $x = z_2 = 0$ dominates
$\Delta(q_0)$ and the irreducible curve
\begin{equation*}
x \left( 1+ \sum_{i=2}^{\beta-1} G_{i1} z_2^i\right) +
 (xG_{\beta,1})(q_0) z_2^\beta  = z = 0
\end{equation*}
dominates $(I_1)_\red$ with a degree $\beta$ map. So $\Upsilon_0'$ contains
an irreducible curve $\wt{I}_1$ dominating $(I_1)_\red$ with a degree
$\beta$ map. And the map $\wt{I}_1\to (I_1)_\red$ is totally ramified at
$q_0$ and $q_1$.

Let us examine the behavior of $\Upsilon'$ over the point $q_1$, where
$\Upsilon$ is defined by \eqref{e9a} in $\Spec \BC[[x, y, z, t]]/(xy -
t^\alpha)$. Without the loss of generality, we may assume that $I_1\subset
\{ y = 0 \}$ locally at $q_1$.

Since $G_{\beta, 1}$ is a meromorphic function with a simple pole at $q_1$,
$G_{\beta, 1}\not\equiv 0$ in the neighborhood of $q_1$. Hence we must have
$\phi_{\beta, 1}(x, 0, 0) = G_{\beta, 1}$ and $\delta_1 = \gamma$ in
\eqref{e9a}.
Hence the curve $\Upsilon_0'$ is given by
\begin{equation}\label{e14}
z_1^\beta + \sum_{i=2}^\beta G_{i1} z_1^{\beta - i} = xy = 0
\end{equation}
in the neighborhood of $q_1$, 
where $z_1 = z / t^\gamma$. By \eqref{e14}, we necessarily have
$G_{i1}(q_1) = 0$ for $i = 2, 3, ..., \beta$; otherwise,
the map $\wt{I}_1\to I_1$ will not be totally ramified at $q_1$.
Therefore, $G_{i2} = G_{i2}(q_1) = G_{i1}(q_1) = 0$ and $\gamma_2 >
\gamma$, which implies the curve in $\Upsilon_0'$ dominating $I_2$ is still
a nonreduced curve isomorphic to $\P^1\times \BC[[z]]/(z^\beta)$. So we
will use the same notation $I_2$ to denote the curve in $\Upsilon_0'$ over
$I_2\subset \Upsilon_0$.

The same argument can be carried out by studying the behavior of
$\Upsilon'$ over $q_2, ..., q_{2l-1}$. Finally we obtained that
$G_{ij} = 0$, $\gamma_j > \gamma$ and $\Upsilon_0'$ contains $I_j\isom
\P^1\times \BC[[z]]/(z^\beta)$ for $i=2,3,...,\beta$ and $j = 2, 3, ...,
2l - 1$. And by symmetry, $\Upsilon_0'$ contains the irreducible curve
$\wt{I}_{2l}$ dominating $I_{2l}$ with a degree $\beta$ map.
So we may take $\Upsilon^{(1)} = \Upsilon'$ in \eqref{e3b}.

Next, we may take a neighborhood of $\Upsilon'$ around $I_2\cup I_3 \cup
... \cup I_{2l-1}$ and go through this procedure again.
Of course, we have to check
that $\wt{I}_1, I_2, ..., I_{2l-1}, \wt{I}_{2l}$ are ``patched'' up at
$q_1, q_2, ..., q_{2l-1}$ as required at the beginning of the
construction. 

Since $\wt{I}_1$ and $\wt{I}_{2l}$ are smooth everywhere, $\Upsilon'$ is
locally given by \eqref{e5} at $q_1$ and $q_{2l-1}$. It is easy to check that
$\Upsilon'$ is locally given by \eqref{e4} at $q_2, q_3, ..., q_{2l-2}$ and
by \eqref{e6} at $w_2, w_3, ..., w_{2l-1}$.
So we may go through the same procedure for $\Upsilon^{(1)} = \Upsilon'$ by
blowing up $\Upsilon'$ along a subscheme of $\Upsilon_0'$ supported on
$(I_2)_\red\cup (I_3)_\red\cup ...\cup (I_{2l-1})_\red$.
The resulting family will be
$\Upsilon^{(2)}$ and the blowup sequence \eqref{e3b} is constructed
inductively in this way.

This finishes the construction of the blowup sequence \eqref{e3b} with only
one thing left to check. We need to check that $\wt{I}_1$ and
$\wt{I}_2$ meet at $q_1$ with multiplicity $\beta$ if $l = 1$, i.e., we
want to show that $\wt{I}_{l}$ and $\wt{I}_{{l+1}}$ meet at the point $q_l$
with multiplicity $\beta$ on $\Upsilon^{(l)}$.

Suppose that $l = 1$. By our previous argument, we still have $G_{i1} =
G_{i2} = 0$ for $i=2,3,..., \beta -1$ and $G_{\beta, 1}(q_1) = G_{\beta,
2}(q_1) = 0$. Since $G_{\beta, 1}$ and $G_{\beta, 2}$ are
meromorphic functions over $\P^1$ with exactly one simple pole at $q_0$ and $q_2$, respectively, each of them has a simple zero at $q_1$. 
Notice that $\phi_{\beta, 1}(x, 0, 0)$ and $\phi_{\beta, 1}(0, y, 0)$ are the
localizations of $G_{\beta, 1}$ and $G_{\beta, 2}$ at $q_1$ in \eqref{e9a}
since $\delta_1 = \gamma$.
Hence $\phi_{\beta, 1}(x, y, 0) = a x + b y + O(x^2, xy, y^2)$ 
for some constants $a, b \ne 0$.
Therefore, $\Upsilon'$ is cut out by
\begin{equation*}
z_1^\beta + ax + by + O(x^2, xy, y^2, t) = 0
\end{equation*}
in $\Spec \BC[[x, y, z_1, t]]/(xy - t^\alpha)$ locally at $q_1$.

\section{Deformation Around a Type II Chain}
\label{sec3}

Our main theorem of this section is the following.

\begin{thm}\label{thm2}
Let $\cup_{i=-l}^{m-1} \overline{q_i q_{i+1}}$ be a Type II chain in
$\Upsilon_0$ and let $q_0 \in \{p_1, p_2, ..., p_{16}\}$. Then
\[
\delta\left(\Upsilon_t, \cup_{i=-l}^{m-1} \overline{q_i q_{i+1}}\right)
\ge \sum_{i=-l}^{m} \intsc(\Sigma_1, E; q_i).
\]
If the equality holds, then

\begin{enumerate}
\item $|\mu(\overline{q_{i-1}q_i}) - 
\mu(\overline{q_i q_{i+1}})| \le 1$ for any $i$ 
(let $\mu(q_iq_{i+1}) = 0$ if $i < -l$ or $i \ge m$);
\item all singularities of $\Upsilon_t$ are nodes in the neighborhood
of $\cup_{i=-l}^{m-1} \overline{q_i q_{i+1}}$.
\end{enumerate}
\end{thm}

\subsection{Some Basic Results on Curve Singularities}\label{s1}

Most of the following results on curve singularities are well known. But we
will prove them here for the lack of a definite reference.

\begin{prop}\label{prop7}
Let $C = \cup_{i=1}^n C_i$ be a reduced curve in $\Delta^2$, 
where $C_i$ are distinct curves in $\Delta^2$.
Then
\begin{equation}\label{prop7-e0}
\delta(C) = \sum_{i=1}^n \delta(C_i) + \sum_{1\le r < s\le n} C_r\cdot C_s,
\end{equation}
where $\delta(C)$ and $\delta(C_i)$ are the $\delta$-invariants of $C$ and
$C_i$ at the origin.
\end{prop}

\begin{proof}
Let $\wt{\Delta^2}$ be the blowup of $\Delta^2$ at the origin and let
$\wt{C}$ and $\wt{C_i}$
be the proper transforms of $C$ and $C_i$ for $i=1,2,..., n$.
Let $\wt{C}$ meet the exception divisor at
points $p_1, p_2, ..., p_l$ and let $\Delta^2(p_j)$ be the neighborhood of
$p_j$ in $\wt{\Delta^2}$.

Let $\wt{C_i} = \cup_{j=1}^l C_i^j$ where $C_i^j\subset \Delta^2(p_j)$ (let
$C_i^j = \emptyset$ if $\wt{C_i}$ does not pass through $p_j$). Let 
$m_i$ be the multiplicity of $C_i$ at the origin and let $m = \sum_{i=1}^n
m_i$ be the multiplicity of $C$ at the origin. We argue by induction on
$\delta(C)$. It is obvious when $\delta(C) = 0$. Suppose that $\delta(C) >
0$. Then $m > 1$.

First, we have
\begin{equation}\label{prop7-e1}
\delta(C_i) = \frac{m_i(m_i -1)}{2} + \delta(\wt{C_i})
= \frac{m_i(m_i -1)}{2} + \sum_{j=1}^l \delta(C_i^j).
\end{equation}
Second, by induction hypothesis, we have
\begin{equation}\label{prop7-e2}
\begin{split}
\delta(C) &= \frac{m(m-1)}{2} +
\sum_{j=1}^l \delta\left(\cup_{i=1}^n C_i^j\right)\\
&= \frac{m(m-1)}{2} + \sum_{j=1}^l \left( \sum_{i=1}^n \delta(C_i^j) + 
\sum_{1\le r < s\le n} C_r^j \cdot C_s^j\right).
\end{split}
\end{equation}
Finally, we have
\begin{equation}\label{prop7-e3}
C_r \cdot C_s = m_r m_s + \sum_{j=1}^l C_r^j \cdot C_s^j.
\end{equation}
Combining \eqref{prop7-e1}, \eqref{prop7-e2} and \eqref{prop7-e3}, we
obtain \eqref{prop7-e0}.
\end{proof}

The next is a parameterized version of \propref{prop7}.

\begin{prop}\label{prop8}
Let $\Upsilon = \cup_{i=1}^n \Upsilon^{(i)}
\subset \Delta^2\times \Delta_t$ be a reduced flat family of curves over
$\Delta_t$ where $\Upsilon^{(i)}$ are distinct flat families of curves
over $\Delta_t$. Then
\begin{equation}\label{prop8-e0}
\delta(\Upsilon_t) = 
\sum_{i=1}^n \delta(\Upsilon_t^{(i)}) + \sum_{1\le r < s\le n} 
\Upsilon_t^{(r)}\cdot \Upsilon_t^{(s)},
\end{equation}
where $\delta(\Upsilon_t)$ and $\delta(\Upsilon_t^{(i)})$ are the total
$\delta$-invariants of the general fibers of $\Upsilon$ and
$\Upsilon^{(i)}$ and the intersection between $\Upsilon_t^{(r)}$ and
$\Upsilon_t^{(s)}$ is taken on the general fiber of $\Delta^2\times
\Delta_t\to \Delta_t$.
\end{prop}

\begin{proof}
After a base change, we may assume that each singular point of $\Upsilon_t$
is given by a section $p: \Delta_t\to \Delta^2\times \Delta_t$. At each
point $p = p(t)$, we have
\begin{equation}\label{prop8-e1}
\delta(\Upsilon_t, p) = 
\sum_{i=1}^n \delta(\Upsilon_t^{(i)}, p) + \sum_{1\le r < s\le n} 
\left(\Upsilon_t^{(r)}\cdot \Upsilon_t^{(s)}\right)_p
\end{equation}
by \propref{prop7}, where $\delta(\Upsilon_t, p)$ and
$\delta(\Upsilon_t^{(i)}, p)$ are the $\delta$-invariants of $\Upsilon_t$
and $\Upsilon_t^{(i)}$ at $p$ and $\left(\Upsilon_t^{(r)}\cdot
\Upsilon_t^{(s)}\right)_p$ is the local intersection number between
$\Upsilon_t^{(r)}$ and $\Upsilon_t^{(s)}$ at $p$ (take
$\left(\Upsilon_t^{(r)}\cdot \Upsilon_t^{(s)}\right)_p = 0$ if
$\Upsilon_t^{(r)}$ and $\Upsilon_t^{(s)}$ do not meet at $p$). 
Obviously, each intersection between $\Upsilon_t^{(r)}$ and
$\Upsilon_t^{(s)}$ is necessarily a singularity of $\Upsilon_t$. So
summing \eqref{prop8-e1} over all sections $p$ of singularities 
yields \eqref{prop8-e0}.
\end{proof}

\begin{cor}\label{cor4}
Let $\Upsilon = \cup_{i=1}^n \Upsilon^{(i)} \subset \Delta^2\times
\Delta_t$ be a reduced flat family of curves over $\Delta_t$ where
$\Upsilon^{(i)}$ are distinct flat families of curves over $\Delta_t$.
If $\Upsilon_0^{(r)}$ and $\Upsilon_0^{(s)}$ meet properly on the central
fiber of $\Delta^2\times \Delta_t\to \Delta_t$ for any $1\le r < s\le n$, then
\begin{equation}\label{cor4-e0}
\delta(\Upsilon_t) \ge 
\sum_{1\le r < s\le n} \Upsilon_0^{(r)}\cdot \Upsilon_0^{(s)},
\end{equation}
where the intersection between $\Upsilon_0^{(r)}$ and
$\Upsilon_0^{(s)}$ is taken on the central fiber of $\Delta^2\times
\Delta_t\to \Delta_t$.
\end{cor}

\begin{proof}
Since $\Upsilon_0^{(r)}$ and $\Upsilon_0^{(s)}$ meet properly,
$\Upsilon_0^{(r)}\cdot \Upsilon_0^{(s)} = \Upsilon_t^{(r)} \cdot
\Upsilon_t^{(s)}$. Then \eqref{cor4-e0} follows from \eqref{prop8-e0}.
\end{proof}

The following is a special case of \coref{cor4}, which is directly
applicable to our situation.

\begin{cor}\label{cor5}
Let $\Upsilon\subset \Delta^2\times \Delta_t$ be a 
reduced flat family of curves over $\Delta_t$ whose central fiber 
$\Upsilon_0$ consists of $n$ irreducible components $\Gamma_1, \Gamma_2, ...,
\Gamma_n$ with multiplicities $\mu_1, \mu_2, ..., \mu_n$, respectively.

Let $\pi: \wt{\Upsilon}\to
\Upsilon$ be the nodal reduction of $\Upsilon$ and 
let $\wt{\Upsilon} = \cup_{i=1}^\alpha \wt{\Upsilon}^{(i)}$ where
$\wt{\Upsilon}^{(i)}$ are the connected components of $\wt{\Upsilon}$.
Suppose that each $\pi(\wt{\Upsilon}_0^{(i)})$ is
supported on $\Gamma_j$ for some $1\le j\le n$. 
Then
\begin{equation}\label{cor5-e0}
\delta(\Upsilon_t) \ge \sum_{1\le r < s \le n} \mu_r\mu_s \left(\Gamma_r\cdot 
\Gamma_s\right).
\end{equation}
\end{cor}

\begin{proof}
By our assumption on $\wt{\Upsilon}$, we can write $\Upsilon$ as $\Upsilon
= \cup_{j=1}^n \Upsilon^{(j)}$ where $\Upsilon_0^{(j)}$ consists of the
component $\Gamma_j$ with multiplicity $\mu_j$, i.e., we let
\[
\Upsilon^{(j)} = \pi\left(\bigcup_{\pi(\wt{\Upsilon}_0^{(i)})_\red \subset \Gamma_j}
\wt{\Upsilon}^{(i)}\right).
\]
Since $\Upsilon_0^{(r)}$ and $\Upsilon_0^{(s)}$ meet properly for any $r\ne s$,
\eqref{cor5-e0} follows from \coref{cor4}.
\end{proof}

\begin{prop}\label{prop9}
Let $\Upsilon, \Upsilon_0, 
\Gamma_j, \mu_j, \pi, \wt{\Upsilon}$ and $\wt{\Upsilon}^{(i)}$ be
defined as in \coref{cor5} except that we further assume
each $\pi(\wt{\Upsilon}_0^{(i)})$ to be reduced, i.e.,
$\pi(\wt{\Upsilon}_0^{(i)}) = \Gamma_j$ for some $j$.
If
\begin{equation}\label{prop9-e0}
\delta(\Upsilon_t) = \sum_{1\le r < s \le n} \mu_r\mu_s,
\end{equation}
then $\Upsilon_t$ only has nodes as singularities.
\end{prop}

\begin{proof}
By our assumptions on $\wt{\Upsilon}$, we have
\[
\Upsilon = \bigcup_{i=1}^n \bigcup_{j=1}^{\mu_i} \Upsilon^{(i, j)}
\]
where $\Upsilon_0^{(i, j)} = \Gamma_i$.
By \propref{prop8},
\[
\delta(\Upsilon_t) = \sum_{i, j} \delta(\Upsilon_t^{(i, j)}) + 
\sum_{(r, p) < (s, q)} \Upsilon_t^{(r, p)}\cdot \Upsilon_t^{(s, q)}
\]
where we define $(r, p) < (s, q)$ if either $r < s$ or $r = s$ and $p < q$.
Notice that
\[
\begin{split}
\sum_{(r, p) < (s, q)} \Upsilon_t^{(r, p)}\cdot
\Upsilon_t^{(s, q)} & \ge \sum_{p, q}\sum_{r < s} \Upsilon_t^{(r, p)}\cdot
\Upsilon_t^{(s, q)}\\
&= \sum_{p, q} \sum_{r < s} \Upsilon_0^{(r, p)}\cdot
\Upsilon_0^{(s, q)} = \sum_{r < s} \mu_r\mu_s
\left(\Gamma_r \cdot \Gamma_s\right)\\
&\ge \sum_{r < s} \mu_r\mu_s.
\end{split}
\]
If the equality holds,
we will necessarily have that $\delta(\Upsilon_t^{(i, j)}) = 0$,
$\Upsilon_t^{(r, p)}\cdot \Upsilon_t^{(s, q)} = 0$ if $r = s$ and $1$ if 
$r\ne s$. This implies that $\Upsilon_t^{(r, p)}$ and $\Upsilon_t^{(s, q)}$
meet transversely for any $r \ne s$ and 
$\Upsilon_t^{(r, p)}\cap \Upsilon_t^{(s, q)} \ne 
\Upsilon_t^{(r', p')}\cap \Upsilon_t^{(s', q')}$ for any $r < s$, $r' < s'$
and $(p, q, r, s)\ne (p', q', r', s')$. Therefore, $\Upsilon_t$ 
has exactly $\sum_{r<s} \mu_r \mu_s$ nodes as singularities, which are the
intersections $\Upsilon_t^{(r, p)}\cap \Upsilon_t^{(s, q)}$ for $r < s$.
\end{proof}

\begin{prop}\label{prop10}
Let $X\subset \Delta_{xyz}^3 \times \Delta_t$ be defined by $xy = t^\alpha
z$ for some $\alpha > 0$ and let $X_0 = R_1\cup R_2$ and $E = R_1\cap R_2$,
where $R_1 = \{x = t = 0\}$ and $R_2 = \{y = t = 0\}$.

Let $\Upsilon\subset X$ be a reduced flat family of curves whose central
fiber $\Upsilon_0$ consists of $m+n$ 
irreducible components $\Gamma_1^{(1)},
\Gamma_2^{(1)}, ..., \Gamma_m^{(1)}, \Gamma_1^{(2)}, \Gamma_2^{(2)}, ...
\Gamma_n^{(2)}$ where $\Gamma_j^{(i)}\subset R_i$ and $\Gamma_j^{(i)} \ne
E$. Let $\mu_{ij}$ be the multiplicity of $\Gamma_j^{(i)}$ in
$\Upsilon_0$.

Let $\pi: \wt{\Upsilon}\to
\Upsilon$ be the nodal reduction of $\Upsilon$ and
let $\wt{\Upsilon} = \cup \wt{\Upsilon}^{(k)}$ where $\wt{\Upsilon}^{(k)}$
are the connected components of $\wt{\Upsilon}$. Suppose that each
$\pi(\wt{\Upsilon}_0^{(k)})$ is supported on $\Gamma_j^{(i)}$ for some $i$
and $j$. Then
\begin{equation}\label{prop10-e0}
\delta(\Upsilon_t) \ge \left(\sum_j \mu_{ij}
\left(\Gamma_j^{(i)}\cdot E\right)\right)
\left(\sum_j \mu_{3-i, j}\right),
\end{equation}
for $i=1,2$.
\end{prop}

\begin{proof}
We resolve the double point $p$ in the same way as in the proof of
\coref{cor3}. Let $\wt{X}\subset X\times \P^1$ be defined by
\[
\frac{x}{t^\alpha} = \frac{z}{y} = \frac{W_1}{W_0}
\]
where $(W_0, W_1)$ are the homogeneous coordinates of $\P^1$.
Let $\wt{X}_0 = \wt{R_1}\cup R_2$ and $\wt{E} = \wt{R_1}\cap R_2$, where
$\wt{R_1}$ is the blowup of $R_1$ at the origin $p$ and 
$\wt{E} = \{ x = y = t = W_0/W_1 = 0\}$. Let $\wt{p} = \{ x = y = z = t =
W_0/W_1 = 0\}$ and let $P = \{ x = y = z = t = 0\}$ be the exceptional
curve of $\wt{X}\to X$. Let $\wt{\Gamma_j^{(1)}}$ be the proper transform of
$\Gamma_j^{(1)}$ and $p_j = \wt{\Gamma_j^{(1)}}\cap P$ for $j =
1,2,...,m$. Our assumptions on $\wt{\Upsilon}$ guarantee 
$p_j\ne \wt{p}$; otherwise, $\wt{\Gamma_j^{(1)}}$ will pass through
$\wt{p}$ and by \lemref{lem1}, each component of $\wt{\Upsilon}_0$ that
dominates $\wt{\Gamma_j^{(1)}}$ will be joined by a chain of curves to a
component dominating $\Gamma_k^{(2)}$ for some $k$.

Let $\Y = \Upsilon\times_X \wt{X}$, $\wt{\Y} = \wt{\Upsilon}\times_X
\wt{X}$ and $\wt{\pi}: \wt{\Y} \to \Y$ be the map induced by
$\pi$. By \lemref{lem1}, $\Y_0$ contains $P$ with multiplicity
\[
\sum_{j=1}^n \mu_{2j}\left(\Gamma_j^{(2)}\cdot E\right).
\]
Let $\wt{\Y} = \cup \wt{\Y}^{(k)}$ where $\wt{\Y}^{(k)} =
\wt{\Upsilon}^{(k)}\times_X \wt{X}$.

It is not hard to see that $\wt{\pi}(\wt{\Y}_0^{(k)})_\red
\subset \wt{\Gamma_j^{(1)}}\cup P$ if $\pi(\wt{\Upsilon}_0^{(k)})_\red\subset
\Gamma_j^{(1)}$ for some $j$. By \lemref{lem1}, 
$P\not\subset \wt{\pi}(\wt{\Y}_0^{(k)})$ if 
$\pi(\wt{\Upsilon}_0^{(k)})_\red\subset\Gamma_j^{(1)}$ for some $j$.
Therefore, we may apply \coref{cor5} to each points $p_j$ for $j=1,2,
..., m$ and obtain
\begin{align*}
\delta(\Upsilon_t) &= \delta(\Y_t) \ge \sum_{j=1}^m \delta(\Y_t, p_j)\\
&\ge \sum_{j=1}^m \left(\mu_{1j} \sum_{l=1}^n
\mu_{2l}\left(\Gamma_l^{(2)}\cdot E\right)\right)\\
&= \left(\sum_{l=1}^n \mu_{2l}
\left(\Gamma_l^{(2)}\cdot E\right)\right)
\left(\sum_{j=1}^m \mu_{1j}\right).
\end{align*}
Our argument needs some trivial change if $p_j$'s fail to be distinct,
which we will leave it to the readers.
\end{proof}

\begin{prop}\label{prop11}
Let $X, \Upsilon, \Gamma_j^{(i)}, \mu_{ij}, \pi, \wt{\Upsilon}$ and
$\wt{\Upsilon}^{(k)}$ be defined as in \propref{prop10} except that we
further assume each $\pi(\wt{\Upsilon}_0^{(k)})$ to be reduced, i.e.,
$\pi(\wt{\Upsilon}_0^{(k)}) = \Gamma_j^{(i)}$ for some $i$ and $j$.
If
\begin{equation}\label{prop11-e0}
\delta(\Upsilon_t) = \left(\sum_{j=1}^m \mu_{1j}\right)
\left(\sum_{j=1}^n \mu_{2j}\right),
\end{equation}
then $\Upsilon_t$ only has nodes as singularities.
\end{prop}

\begin{proof}
The proof is the same as the proof of \propref{prop10} except that we need
to apply \propref{prop9} to each point $p_j$ at the last step.
\end{proof}

\subsection{Proof of \thmref{thm2}}

Without the loss of generality, let us assume that
$\overline{q_0q_1}\subset R_1$. Let $\alpha_1 = m_1(q_0)$, $\alpha_2 =
m_2(q_0)$ and $\mu_i = \mu(\overline{q_i q_{i+1}})$ for each $i$.

Let $w_i$ be the intersection between $\overline{q_i q_{i+1}}$ and
$\Gamma_1\cup \Gamma_2$ other than the points $q_i$ and $q_{i+1}$,
if such intersection exists.

Basically, $\Upsilon_t$ has singularities in the neighborhoods of $q_0$ and
$w_i$. \propref{prop2a} tells us the configurations of $\wt{\Upsilon}$ over
these neighborhoods, while the series of results we have obtained in \ref{s1}
can be used to estimate the $\delta$-invariants of $\Upsilon_t$ in these
neighborhoods.

By \propref{prop2a} and \ref{prop10},
\begin{equation}\label{thm2-e1}
\delta(\Upsilon_t, q_0) \ge (\alpha_1 + \mu_0)^2 = (\alpha_2 + \mu_{-1})^2.
\end{equation}

If $\alpha_1 = 0$, the point $w_0$ exists.
By \propref{prop2a} and \coref{cor5},
\begin{equation}\label{thm2-e2}
\delta(\Upsilon_t, w_0)\ge \mu_0
\end{equation}
if $\alpha_1 = 0$.
Since $\alpha_1 = 0$ or $1$, we may write \eqref{thm2-e2} in the form
\begin{equation}\label{thm2-e3}
\delta(\Upsilon_t, w_0) \ge (1-\alpha_1)\mu_0.
\end{equation}
Similarly, we have
\begin{equation}\label{thm2-e4}
\delta(\Upsilon_t, w_{-1}) \ge (1-\alpha_2)\mu_{-1}.
\end{equation}

Obviously, for $i > 0$, the point $w_i$ exists if $\mu_{i-1} = \mu_{i}$.
By \propref{prop2a} and \coref{cor5}, for $i > 0$,
\begin{equation}\label{thm2-e5}
\delta(\Upsilon_t, w_i) \ge \mu_i
\end{equation}
if $\mu_{i-1} = \mu_i$. Similarly, for $i > 1$,
\begin{equation}\label{thm2-e6}
\delta(\Upsilon_t, w_{-i}) \ge \mu_{-i}
\end{equation}
if $\mu_{-i} = \mu_{-i + 1}$.

Let $0\le a_0 < a_1 < a_2 < ... < a_n < ...$ 
be the sequence of integers such that
\[
\begin{split}
\mu_0 = ... = \mu_{a_0} & > \mu_{a_0+1} = \mu_{a_0+2} = ... = \mu_{a_1} > 
\mu_{a_1+1} = \mu_{a_1+2}
= ... = \mu_{a_2}\\
& > ... > \mu_{a_{n-1} + 1} = \mu_{a_{n-1} + 2} = ... = \mu_{a_n} > ... .
\end{split}
\]

Then by \eqref{thm2-e5},
\begin{equation}\label{thm2-e7}
\sum_{i>0} \delta(\Upsilon_t, w_i) \ge a_0 \mu_0 + 
\sum_{i > 0} (a_i - a_{i-1} - 1)\mu_{a_i}.
\end{equation}
And since
\[
\sum_{i > 0} \intsc(\Sigma_1, E; q_i) = (a_0 + 1)\mu_0 +\sum_{i>0} (a_i -
a_{i-1}) \mu_{a_i},
\]
\begin{equation}\label{thm2-e8}
\begin{split}
\sum_{i>0} \delta(\Upsilon_t, w_i) &- \sum_{i > 0} \intsc(\Sigma_1, E; q_i)
\ge -\left(\mu_0 + \sum_{i>0} \mu_{a_i}\right)\\
&\ge -\left(\mu_0 + (\mu_0 - 1) + (\mu_0 - 2) + ... + 1\right) =
-\frac{\mu_0(\mu_0+1)}{2}.
\end{split}
\end{equation}
Similarly,
\begin{equation}\label{thm2-e9}
\sum_{i>1} \delta(\Upsilon_t, w_{-i}) - \sum_{i \ge 1} \intsc(\Sigma_1, E; 
q_{-i})\ge -\frac{\mu_{-1}(\mu_{-1}+1)}{2}.
\end{equation}

Combining \eqref{thm2-e1}, \eqref{thm2-e3}, \eqref{thm2-e4},
\eqref{thm2-e8} and \eqref{thm2-e9}, we get
\begin{align}\label{thm2-e10}
& \delta(\Upsilon_t,  \cup_{i=-l}^{m-1} \overline{q_i q_{i+1}})
- \sum_{i=-l}^m \intsc(\Sigma_1, E; q_i)\\
& \ge \delta(\Upsilon_t, q_0)
+ \left(\sum_{i = 0, -1} \delta(\Upsilon_t, w_i)  -  
\intsc(\Sigma_1, E; q_0)\right)\notag\\
& \quad + \left(\sum_{i > 0} \delta(\Upsilon_t, w_i) - \sum_{i>0}
\intsc(\Sigma_1, E; q_i)\right)\notag\\
& \quad + \left(\sum_{i > 1} \delta(\Upsilon_t, w_{-i}) - \sum_{i\ge 1}
\intsc(\Sigma_1, E; q_{-i})\right)\notag\\
& \ge \left(\frac{1}{2}(\alpha_1 + \mu_0)^2 + 
\frac{1}{2}(\alpha_2 + \mu_{-1})^2\right)\notag\\
&\quad + \left( (1-\alpha_1) \mu_0 + (1-\alpha_2)\mu_{-1} 
- \frac{1}{2}(\alpha_1 +
\mu_0) - \frac{1}{2} (\alpha_2 + \mu_{-1})\right)\notag\\
&\quad - \frac{\mu_0(\mu_0 + 1)}{2} -
\frac{\mu_{-1}(\mu_{-1}+1)}{2}\notag\\
&= \frac{1}{2}(\alpha_1^2 - \alpha_1) + \frac{1}{2}(\alpha_2^2 - \alpha_2)
= 0.\notag
\end{align}

Suppose that the equality in \eqref{thm2-e10} holds. Then the equalities in
\eqref{thm2-e8} and \eqref{thm2-e9} have to hold, which implies $|\mu_{i-1} -
\mu_i| = 0$ or $1$ for any $i$. Under this condition,
it is not hard to see by \propref{prop2a} that each component of
$\wt{\Upsilon}_0$ dominating a component
$\overline{q_i q_{i+1}} \subset \cup_{i=-l}^{m-1} \overline{q_i q_{i+1}}$ must
dominate $\overline{q_i q_{i+1}}$ with a degree one map. So we may apply
\propref{prop9} to $\Upsilon$ at $w_i$ and apply \propref{prop11} to
$\Upsilon\subset X$ at $q_0$ to conclude that $\Upsilon_t$ only has nodes
as singularities in the neighborhoods of $w_i$ and $q_0$.
On the other hand, since the equality in
\eqref{thm2-e10} holds, $\Upsilon_t$ does not have singularities anywhere
else in the neighborhood of $\cup_{i=-l}^{m-1} \overline{q_i q_{i+1}}$.
Therefore, $\Upsilon_t$ only has nodes as singularities in the
neighborhood of $\cup_{i=-l}^{m-1} \overline{q_i q_{i+1}}$ if the equality
in \eqref{thm2-e10} holds.

\section{Deformation Around a Type III Chain}
\label{sec4}

The deformation of $\Upsilon$ around a Type III chain is much more
complicated than the other two cases. Instead of attempting a direct
analysis as we did in the previous two sections, 
we will employ a completely different approach. Our strategy can be briefly
described as follows.

We are trying to move $[R^s]$ around in $\wt{\H_g}$,
or equivalently, apply a monodromy action to $[R^s]$ such that a limiting
rational curve $\Upsilon_0$ on $R^s$ will become another limiting
rational curve $\Upsilon_0'$ and meanwhile a Type III chain in $\Upsilon_0$
will be ``transformed'' to a Type II chain in $\Upsilon_0'$. 

Our main theorem of this section is the following.

\begin{thm}\label{thm3}
Let $\overline{r_{ij}}\cup \overline{r_{ij} q_1}\cup ...\cup
\overline{q_{m-1} q_m}$ be a Type III chain in $\Upsilon_0$. Then
\[
\delta(\Upsilon_t, \overline{r_{ij}}\cup \overline{r_{ij} q_1}\cup
...\cup \overline{q_{m-1} q_m})\ge \sum_{l=0}^m \intsc(\Sigma_1, E; q_l)
\]
(let $q_0 = r_{ij}$) and if the equality holds, all singularities of
$\Upsilon_t$ are nodes
in the neighborhood of $\overline{r_{ij}}\cup
\overline{r_{ij} q_1}\cup ...\cup \overline{q_{m-1} q_m}$.
\end{thm}

\subsection{Notations and Definitions}

In order to make our ideas precise, we need to introduce a few new objects.

Since we assume that $\wt{\W_g}$ is irreducible, 
$C$ has the same number of singularities for all general $([S], [C])\in
\wt{\W_g}$. Let $C$ have exactly $\beta$ singular points for a general
$([S], [C])\in \wt{\W_g}$.

Let $\Z_g$ be the incidence correspondence
\[
\begin{split}
\Z_g = & \big\{ ([S], [C], w_1, w_2, ..., w_\beta, \delta_1, \delta_2,
..., \delta_\beta, \nu_1, \nu_2, ..., \nu_\beta):\\ 
& \quad\quad ([S], [C]) \in \W_g \text{ general},\\
& \quad\quad C_\sing = \{ w_1, w_2, ..., w_\beta \},\\
& \quad\quad \text{$\delta_i$ is $\delta$-invariant of $C$ at $w_i$},\\
& \quad\quad \text{$\nu_i = 1$ if $w_i$ is a node of $C$; $\nu_i = 0$ if it is not}
\big\}\\
& \subset \wt{\W_g}\times(\P^{4g-3})^\beta\times \BZ^{2\beta}.
\end{split}
\]
Again, we first define $\Z_g$ for general K3 surfaces $S$ and
then take its closure $\wt{\Z_g}$
in $\wt{\W_g}\times(\P^{4g-3})^\beta\times \BZ^{2\beta}$.

Notice that the last $2\beta$ coordinates
$(\delta_1, \delta_2, ..., \delta_\beta, \nu_1, \nu_2,
..., \nu_\beta)$ are the same throughout an irreducible component of
$\wt{\Z_g}$. \thmref{T:NOD} is equivalent to the statement that $\nu_1 =
\nu_2 = ... = \nu_\beta = 1$ on each irreducible component of
$\wt{\Z_g}$ which dominates $\wt{\H_g}$.
Without the loss of generality, let us assume that $\wt{\Z_g}$ is
irreducible and dominates $\wt{\H_g}$.

\begin{defn}
Let $([S], [C], w_1, w_2, ..., w_\beta, ...)\in \wt{\Z_g}$. We will call
$w_i$'s ``limiting singularities'' of $C$.
\end{defn}

We observe that

\begin{prop}\label{prop3a}
Let $[R^s]\in \wt{\CR_g}$ general and let $\Upsilon_0$ be a limiting
rational curve on $R^s$.
Then each limiting singularity of $\Upsilon_0$
lies on an $F$-chain in $\Upsilon_0$.
\end{prop}

\begin{proof}
Notice that $\Upsilon_0$ is smooth outside of its $F$-chains.
\end{proof}

The same holds for a general $[R^s]\in \wt{\CR_g^0}$.

\begin{prop}\label{prop3b}
Let $[R^s]\in \wt{\CR_g^0}$ general and let
$\Upsilon_0$ be a limiting rational curve on $R^s$.
Then each limiting singularity of $\Upsilon_0$ is
either a point on a component $\overline{pq}\subset
\Upsilon_0, \overline{pq}\in |\CO_{R}(F)|$ or an
point among $\Gamma_1\cap E = \Gamma_2\cap E$.
\end{prop}

\begin{proof}
Again, the reason for this fact is trivially that $\Upsilon_0$ is smooth
outside of the locus described.
\end{proof}

Let us consider the K3 surfaces $S$ with Picard lattice
\begin{equation}\label{thm3-e1}
\begin{pmatrix}
0 & 2\\
2 & 0
\end{pmatrix}
\end{equation}
or
\begin{equation}\label{thm3-e2}
\begin{pmatrix}
-2 & 2\\
2 & 0
\end{pmatrix}
\end{equation}
or equivalent, 
K3 surfaces whose Picard groups are generated by two divisors $C$
and $F$ with $C\cdot F = 2$, $F^2 = 0$ and $C^2 = 0$ or $C^2 = -2$. 
For lacking a good name for such surfaces, we will simply call them
{\it elliptic K3\/}'s. Although elliptic K3 surfaces usually refer to all K3
surfaces that admit an elliptic fiberation, we will use the term in our
context to refer to K3 surfaces with Picard lattice \eqref{thm3-e1} or
\eqref{thm3-e2}.
It is not hard to see that such K3 surfaces $S$ can be realized as
double covers
of $\P^1\times\P^1$ or $\F_1$ ramified along a curve in $|4C + 4F|$ or $|4C
+ 6F|$. 

Notice that the divisor $C+kF$ is ample but not very ample on an elliptic
K3 surface $S$.
To embed $S$ into a projective space, we need to use the divisor
$2C+2kF$.
Embed $S$ into $\P^{4g-3}$ ($g\ge 3$) by $|2C + 2kF|$ and we observe that
$[S]\in \H_g$, where $g = 2k+1$ if $S$ has Picard lattice \eqref{thm3-e1}
and $g = 2k$ if $S$ has Picard lattice \eqref{thm3-e2}.

Let $\Y_g\subset \H_g$ be the locus in
$\H_g$ consisting of elliptic K3's with Picard lattice \eqref{thm3-e1} or
\eqref{thm3-e2} embedded into $\P^{4g-3}$ by $|2C + 2kF|$, where $k =
\lfloor g/2 \rfloor$.
And let $\wt{\Y_g}$ be the closure of
$\Y_g$ in $\wt{\H_g}$, where we may regard $\Y_g$ as a subscheme of
$\wt{\H_g}$ since $\Y_g$ is disjoint from $\CR_g$.
Obviously, $\Y_g$ and $\wt{\Y_g}$ are irreducible
and have codimension 1 in $\H_g$ and $\wt{\H_g}$, respectively. And

\begin{prop}\label{prop3c}
We have $\wt{\Y_g}\cap \wt{\CR_g} = \wt{\CR_g^0}$.
\end{prop}

\begin{proof}
We realize $S\in\Y_g$ as double covers of $\P = \P^1\times \P^1$ or $\F_1$
totally ramified along a curve $C\in |-2K_\P|$, where $K_\P$ is the canonical
divisor of $\P$. We want to construct a double cover $Y$ of 
$X = |-2K_\P|\times \P$ whose restriction to a point $C\in |-2K_\P|$ is the
double cover of $\P$ ramified along $C$.

Let $B\subset X = |-2K_\P|\times \P$ be the universal family of curves in
$|-2K_\P|$ and let $\pi_1$ and $\pi_2$ be the projections of $X$ to
$|-2K_\P|$ and $\P$, respectively.
Since
\[
\CO_X(B) = \pi_1^* \CO_{|-2K_\P|}(1) \times \pi_2^* \CO_\P(-2K_\P) \ne
L^2
\]
for any line bundle $L$, we cannot construct $Y$ directly as a double cover
of $X$ along $B$. The remedy for this situation is
trivial. We take a general $N$-dimensional linear subsystem of 
$|\CO_{|-2K_\P|}(2)|$, where $N= \dim |-2K_\P|$, and use it to map
$|-2K_\P|$ to itself. For example, after fixing homogeneous coordinates
$(Z_0, Z_1, ..., Z_N)$ of $|-2K_\P|$, we may simply take the map to be
sending $(Z_0, Z_1, ..., Z_N)$ to $(Z_0^2, Z_1^2, ..., Z_N^2)$. This will
induce a map $f: X\to X$. Obviously,
\[
\begin{split}
f^* \CO_X(B) &= \pi_1^* \CO_{|-2K_\P|}(2) \times \pi_2^* \CO_\P(-2K_\P)
\\
&= \left(\pi_1^* \CO_{|-2K_\P|}(1) \times \pi_2^* \CO_\P(-K_\P)\right)^2.
\end{split}
\]
So there exists a double cover $Y$ of $X$ ramified along $f^{-1} (B)$. It
is trivial to check that the fiber of $Y\to |-2K_\P|$ over an irreducible 
$C\in |-2K_\P|$ is the double cover of $\P$ ramified along $C$, while the
fiber of $Y\to |-2K_\P|$ over a double curve $C = 2D$ with $D\in |-K_\P|$ is
a surface $R\in \CR_g^0$. Therefore, $\CR_g^0$ lies on the closure of
$\Y_g$ in $\H_g$. 

Notice that both $\wt{\Y_g}$ and $\wt{\CR_g}$ has codimension 1 in
$\wt{\H_g}$. Obviously, $\wt{\Y_g}\cap \wt{\CR_g}\subset \wt{\CR_g^0}$. The
previous argument shows that $\wt{\Y_g}\cap \wt{\CR_g}\ne \emptyset$. And
since $\wt{\CR_g^0}$ is irreducible with codimension 2 in $\wt{\H_g}$,
$\wt{\Y_g}\cap \wt{\CR_g} = \wt{\CR_g^0}$.
\end{proof}

\subsection{Rational Curves on Elliptic K3's}

As mentioned at the beginning of this section, our strategy is to move
$[R^s]$ around in $\wt{\H_g}$. The way we move $[R^s]$ can now be described
as follows:
\begin{equation}\label{thm3-e3}
[R^s]\in \wt{\CR_g}\imply [R^s]\in \wt{\CR_g^0} \imply [S]\in \Y_g\imply
[R^s] \in \wt{\CR_g^0}\imply [R^s]\in \wt{\CR_g}.
\end{equation}
Namely, we first degenerate a $R^s$ in $\wt{\CR_g}$ to a $R^s$ in
$\wt{\CR_g^0}$, then move a $R^s$ in $\wt{\CR_g^0}$ to an $S$ in $\Y_g$ and
so on.

It is clear what happens to the limiting rational curves on a $R^s$ in
$\wt{\CR_g}$ when it degenerates to a $R^s$ in $\wt{\CR_g^0}$. However, we
do not have much idea about what happens to the rational curves on
an $S$ in $\Y_g$ when it degenerates to a $R^s$ in $\wt{\CR_g^0}$. So we
will prove a collection of results regarding rational curves on an elliptic
K3 and their degeneration as an elliptic K3 degenerates to a special union
of scrolls.

On each $S$ in $\Y_g$, there is a pencil of elliptic curves in $|F|$. There
are exactly $24$ rational curves in $|F|$ for $S$ general. We are
interested in their monodromy group as $S$ varies in $\Y_g$.

\begin{prop}\label{prop3}
The monodromy group of the $24$ rational curves in $|F|$ as $S$ varies in
$\Y_g$ is the full symmetric group.
\end{prop}

\begin{proof}
We may realize $S$ as the double cover of $\P = \P^1\times \P^1$ or
$\F_1$ ramified along a curve $C\in |-2 K_{\P}|$, where $K_\P$ is the
canonical divisor of $\P$. Let $\pi: S\to \P$ be the covering map. We
use the same notation $F$ to denote the divisor on $\P$ such that $\pi^*(F)
= F$ on $S$. Obviously, the $24$ rational curves in $|\CO_S(F)|$ correspond
to the $24$ curves in $|\CO_\P(F)|$ tangent to $C$. Then the statement of the
proposition is equivalent to saying that the monodromy group $G$ of the
$24$ curves in $|\CO_\P(F)|$ tangent to $C$ is the full symmetric group as
$C$ varies in $|-2 K_{\P}|$.

Following the same line of argument as in \cite{H}, we show that $G$ is the
symmetric group by arguing that

\begin{enumerate}
\item $G$ is twice transitive;
\item $G$ contains a simple transposition.
\end{enumerate}

To see that $G$ is twice transitive, let $W$ be the incidence correspondence
\[
\begin{split}
W = & \big\{(C, I_1, I_2): C\in |-2 K_{\P}|,\\
& \quad \text{$I_1$ and $I_2$ are two different curves in 
$|\CO_\P(F)|$ tangent to $C$}\big\}\\
& \subset |-2 K_{\P}|\times |\CO_\P(F)|\times |\CO_\P(F)|.
\end{split}
\]
Saying that $G$ is twice transitive is equivalent to saying that $W$ is
irreducible. To see that $W$ is irreducible, it suffices to project $W$ to
\(|\CO_\P(F)|\times |\CO_\P(F)|\). 
It is obvious that $W$ dominates \(|\CO_\P(F)|\times |\CO_\P(F)|\) and 
$p(W) = \{(I_1, I_2): I_1 \ne I_2\in |\CO_\P(F)|\}$, 
where \(p: W\to |\CO_\P(F)|\times |\CO_\P(F)|\) is the projection. 
And the fibers of \(p: W\to |\CO_\P(F)|\times |\CO_\P(F)|\) 
are irreducible and have the same dimension
everywhere. Therefore, $W$ is irreducible and $G$ is twice transitive.

To see that $G$ contains a simple transposition, let $C_0$ be a curve in
$|-2K_\P|$ having exactly one node $p$ and smooth everywhere else. 
If there is a family of curves
$C_t \in |-2K_\P|$ whose central fiber is $C_0$, two
out of the $24$ curves in $|\CO_\P(F)|$ tangent to $C_t$ will degenerate to
the curve in $|\CO_\P(F)|$ passing through the node $p$ as $t\to 0$.
It is easy to see that a loop around $C_0$ will be lifted to
a simple transposition which transposes these two curves.
\end{proof}

Let $\pi: \Delta\to \wt{\Y_g}\subset \wt{\H_g}$ be a morphism
from the disk $\Delta$ to $\wt{\Y_g}$, where $\pi(0) = [R^s]\in
\wt{\CR_g^0}$ and $\pi(t)$ is a general point in $\Y_g$ for $t\in \Delta$
general. Let $\{p_1, p_2, ..., p_{16}\}$ be the vanishing locus of $s$.

Let $X = \wt{\S_g}\times_\wt{\H_g} \Delta$ be the family of K3 surfaces
corresponding to $\pi$. As before, $X$ has sixteen rational double points
$p_1, p_2, ..., p_{16}$ lying on the double curve $E = R_1\cap R_2$, where
$X_0 = R = R_1\cup R_2$.

Let $r_{ij}$ be the points on $E$ defined as before, i.e., $\CO_E(2r_{ij})
= \CO_E(F_i)$. Since $R\in \CR_g^0$, $r_{1j} = r_{2j}$. Let $q_i$ be the
point on $E$ such that $\CO_E(p_i + q_i) = \CO_E(F)$ for $i=1,2,...,16$.

\begin{prop}\label{prop4}
As $X_t$ degenerates to $X_0 = R$, the $24$ rational curves in $|\CO_{X_t}(F)|$
on $X_t$ behave in the following way

\begin{enumerate}
\item two of them degenerate to
$\overline{r_{1j}}\cup \overline{r_{2j}}$ for each $1\le j \le 4$;
\item one of them degenerates to $\overline{p_i q_i}^{(1)}\cup \overline{p_i
q_i}^{(2)}$ for each $1\le i \le 16$.
\end{enumerate}
\end{prop}

\begin{proof}
The first statement follows directly from the 
Caporaso-Harris theorem on tacnodes \cite[Lemma 4.1]{CH}.

Let us resolve the sixteen double points $p_1, p_2, ..., p_{16}$ of $X$
by blowing up $R_1$ at these points. Let $\wt{X}$ be the resulting family
and $\wt{X}_0 = \wt{R_1}\cup R_2$, where $\wt{R_1}$ is the blowup of $R_1$
at the sixteen points. The total transform of $\overline{p_i q_i}^{(1)}\cup
\overline{p_i q_i}^{(2)}$ consists of the proper transform $I_i$ of
$\overline{p_i q_i}^{(1)}$, the exceptional divisor $P_i$ and $\overline{p_i
q_i}^{(2)}$. Obviously, $I_i$ and $P_i$ meet transversely at a point
$w_i$. By the same line of argument as in \cite{C}, we can deform $I_i\cup
P_i\cup \overline{p_i q_i}^{(2)}$ to a curve on the general fiber while
preserving a node in the neighborhood of $w_i$. This concludes the second
statement. 
\end{proof} 

Let $\Upsilon\subset X$ be a one-parameter family of curves whose general
fiber $\Upsilon_t\in |\CO_{X_t}(C + kF)|$ is irreducible and rational. Let
$\Sigma_1, \Sigma_2, \Gamma_1, \Gamma_2$ be the components of $\Upsilon_0$
defined as before.

Let $\overline{pq} = \overline{pq}^{(1)}\cup \overline{pq}^{(2)}\in
|\CO_R(F)|$ for $p, q\in E$ and $\CO_E(p+q) = \CO_E(F)$. 
If $\overline{pq}\subset \Upsilon_0$,
from the proof of \propref{prop1b}, we see that $\overline{pq}$
can only be one of the following

\begin{enumerate}
\item $p = p_i$ and $q = q_i$ for some $i$;
\item $p = q = r_{1j} = r_{2j}$ for some $j$; and if
this is the case, $r_{1j}$ does not lie on $\Gamma_1\cup
\Gamma_2$;
\item $p, q\not\in \{r_{ij}, p_l, q_l\}$ and either $p$ or $q$ lies on
$\Gamma_1\cup \Gamma_2$.
\end{enumerate}

The following two propositions say that we can be more specific in this case
because $\Upsilon_0$ is a limit of rational curves on elliptic K3's, 
while all the analysis we did in the proof of \propref{prop1b} is solely
based on the fact that $\Upsilon_0$ is a limit of rational curves on smooth
K3 surfaces.

\begin{prop}\label{prop5}
Let $\overline{pq} = \overline{pq}^{(1)}\cup 
\overline{pq}^{(2)} \in|\CO_R(F)|$ and
$\overline{pq} \subset \Upsilon_0$.
If $\overline{pq} = \overline{p_i q_i}$ for some $i$, then

\begin{enumerate}
\item either $p = p_i\in \Gamma_1\cup \Gamma_2$ or
$q = q_i\in \Gamma_1\cup \Gamma_2$; 
\item $\delta(\Upsilon_t, \overline{pq}) = \intsc(\Sigma_1, E; p) +
\intsc(\Sigma_1, E; q)$.
\end{enumerate}
\end{prop}

\begin{prop}\label{prop5z}
We have $\overline{r_{1j}}\cup \overline{r_{2j}}\not\subset \Upsilon_0$ for
any $1\le j\le 4$.
\end{prop}

We need a few algebraic lemmas.

Let $A$ be an integral domain. We use the notation $\overline{A}$ to denote
the integral closure of $A$ in the algebraic closure of its quotient field.
It is well known that $\overline{\BC[[t]]} = \BC[[t, \sqrt{t},
\sqrt[3]{t},..., \sqrt[n]{t}, ...]]$.

\begin{lem}\label{lema1}
Let $A = \BC[x, [y]]\tensor_\BC \overline{\BC[[t]]}$ and $W\subset A$ be
\[
W = \left\{ \sum_{i=0}^\infty \left( a_i(t) y^i + b_i(t) xy^i\right):
a_i(t), b_i(t) \in \overline{\BC[[t]]}\right\}.
\]
Let $f(x, y, t)\in W$ satisfying $f(x, y, 0)\ne 0$.
If $f(x, y, t)$ is reducible in $A$, then there exists $\alpha(t)\in
\overline{\BC[[t]]}$ such that $\alpha(0) = 0$ and
$f(x, \alpha(t), t) = 0$. 
\end{lem}

\begin{proof}
Let $f(x, y, t) = g(x, y, t) h(x, y, t)$ where neither of $g(x, y, t)$ 
and $h(x, y, t)$ is a unit in $A$.
Let $\deg_x g(x, y, t)$ and $\deg_x h(x, y, t)$ be the degrees of $g(x, y,
t)$ and $h(x, y, t)$ as polynomials in $x$. Obviously, since $f(x, y, t)\in
W$, at least one of $\deg_x g(x, y, t)$ and $\deg_x h(x, y, t)$ is zero.
Assume that $\deg_x g(x, y, t) = 0$, i.e., $g(x, y, t) = r(y, t)$ for some
$r(y, t) \in \BC[[y]]\tensor_\BC \overline{\BC[[t]]}$. Since $g(x, y, t)$
is not a unit in $A$, $r(0, 0) = 0$. And $f(x,y,0)\ne 0\imply r(y, 0)\ne
0$. So using Weierstrass Preparation Theorem, we can show that
there exists $\alpha(t)\in \overline{\BC[[t]]}$ such that $\alpha(0) = 0$
and $r(\alpha(t), t) = 0$, which is what we want.
\end{proof}

A corollary of \lemref{lema1} is

\begin{cor}\label{cor2}
Let $A = \BC[x, y]\tensor_\BC \overline{\BC[[t]]}$ and $W\subset A$ be
\[
W = \left\{ \sum \left( a_i(t) y^i + b_i(t) xy^i\right):
a_i(t), b_i(t) \in \overline{\BC[[t]]}\right\}.
\]
Let $f(x, y, t)\in W$ satisfying $f(x, y, 0)\ne 0$.
If $f(x, y, t)$ is irreducible in $A$, then
$f(x, y, t)$ is also irreducible in
$\BC[x, [y]]\tensor_\BC \overline{\BC[[t]]}$.
\end{cor}

\begin{proof}
By \lemref{lema1}, if $f(x, y, t)$ is reducible in $\BC[x, [y]]\tensor_\BC
\overline{\BC[[t]]}$, $f(x, \alpha(t), t) = 0$ for some
$\alpha(t)\in \BC[[t]]$ and $\alpha(0) = 0$. On the other hand, $f(x, y, t)
= g(x, y, t)(y - \alpha(t)) + a(t) x + b(t)$ for some $g(x, y, t)\in A$
and $a(t), b(t)\in \overline{\BC[[t]]}$. So we necessarily have $a(t) =
b(t) = 0$ and $f(x, y, t)$ must be reducible in $A$. Contradiction.
\end{proof}

\begin{lem}\label{lema2}
Let $C\subset \Delta_{xyz}^3 / (xy = 0) = R_1\cup R_2$ be a reduced curve
cut out by $f(x,y,z) = 0$ for some $f(x, y, z) \in \BC[[x, y, z]]$, where
$R_1 = \Delta_{xyz}^3/(x = 0)$ and $R_2 = \Delta_{xyz}^3/(y = 0)$. Suppose
that $C$ does not contain the double curve $E =R_1\cap R_2$. Then
\begin{equation}\label{lema2-e0}
\delta(C) = \delta(C_1) + \delta(C_2) + C_1\cdot E = \delta(C_1) +
\delta(C_2) + C_2\cdot E
\end{equation}
where $C = C_1\cup C_2$, $C_i \subset R_i$ and the intersection $C_i\cdot
E$ is taken on $R_i$ for $i=1,2$.
\end{lem}

\begin{proof}
Let $\wt{\CO_C}$ and $\wt{\CO_{C_i}}$ be the normalizations of the coordinate
rings of $C$ and $C_i$, respectively.
By the definition of $\delta$-invariants,
$\delta(C) = l(\wt{\CO_C}/\CO_C)$, where $l(M)$ is the length of an $\CO_C$
module $M$.

Obviously, $\wt{\CO_C} = \wt{\CO_{C_1}}\times \wt{\CO_{C_2}}$. Therefore,
\begin{align*}
\delta(C) &= l(\wt{\CO_C}/\CO_C) = 
l\left((\wt{\CO_{C_1}}\times \wt{\CO_{C_2}})/\CO_C\right)\\
&= l\left((\wt{\CO_{C_1}}\times \wt{\CO_{C_2}})/(\CO_{C_1}\times
\CO_{C_2})\right) + l((\CO_{C_1}\times \CO_{C_2})/\CO_C)\\
&= l(\wt{\CO_{C_1}}/\CO_{C_1}) + l(\wt{\CO_{C_2}}/\CO_{C_2}) + C_i\cdot E,
\end{align*}
for $i=1,2$.
\end{proof}

\begin{lem}\label{lema3}
Let \(p, q, f\in \BC[z]\) be the polynomials in \(z\) such that
\(\deg f \le 1\), \(\gcd (p, q) = 1\) and \( (p^2 - f q^2)^2 - q^4 \) is a
perfect square in \(\BC[z]\), where \(\gcd (p, q)\) is the greatest common
divisor of \(p\) and \(q\). Then 
\begin{enumerate}
\item both $p$ and $q$ are constants;
\item $f$ is a constant if $q \ne 0$.
\end{enumerate}
\end{lem}

\begin{proof}
Let \(f_1 = f - 1\) and \(f_2 = f + 1\). Obviously, \(\gcd( p^2 - f_1 q^2,
p^2 - f_2 q^2) = 1\) since \(\gcd (p, q) = 1\). And since \( (p^2 - f_1
q^2) (p^2 - f_2 q^2) \) is a perfect square, both \(p^2 - f_1
q^2\) and \(p^2 - f_2 q^2\) are perfect squares.

It is easy to show that if one of \(p\) and \(q\) is a constant, the other is
a constant too. On the other hand, if both are constants and $q\ne 0$,
\(f\) has to be
a constant too. Suppose that \(p, q\not\in\BC\).

Solve \(p^2 - f_i q^2 = r_i^2\) for \(i = 1,2\) and we obtain
that there exist \(q_1, q_2, q_3, q_4 \in \BC[z]\)
such that \(q = q_1 q_2 = q_3 q_4\) and
\[
p = \frac{ q_1^2 + f_1 q_2^2}{2} = \frac{q_3^2 + f_1 q_4^2}{2}.
\]
Since \(q_1 q_2 = q_3 q_4\),
there exist \( s, t, u, v\in \BC[z]\) such that \(q_1 = su\), \(q_2 = tv\),
\(q_3 = tu\) and \(q_4 = sv\). Therefore,
\[
p = \frac{s^2 u^2 + f_1 t^2 v^2}{2} = \frac{t^2 u^2 + f_2 s^2 v^2}{2}
\]
and hence
\[
\frac{u^2 - f_1 v^2}{u^2 - f_2 v^2} = \frac{s^2}{t^2}
\]
which implies that \( (u^2 - f_1 v^2) (u^2 - f_2 v^2) \) is a perfect
square.

Obviously, \(\gcd(u, v) = 1\). Since \(q\not\in \BC\), \(u\) and \(v\)
cannot both be constants. Combining with our previous argument, neither of
\(u\) and \(v\) is a constant. Therefore, \(\max (\deg u, \deg v) < \deg q
\le \max(\deg p, \deg q)\). So this procedure cannot go on forever. A
contradiction.
\end{proof}

Now let us go back to the proof of \propref{prop5} and \ref{prop5z}. 

Our proofs of both statements are based on the construction
a blowup sequence over $\overline{pq}$ and an induction on the multiplicity
of $\overline{pq}$ in $\Upsilon_0$.
It turns out that this process can be described more clearly if we study
the behaviors of $X$ and $\Upsilon$ in the analytic neighborhood of
$\overline{pq}$, or alternatively, study the formal completion of $X$ and
$\Upsilon$ along $\overline{pq}$. But we will stick to the language of
analytic geometry for it being more intuitive.

The following proposition is lengthy to state due to the fact that we
need to give a precise description of $X$ and $\Upsilon$ in the
neighborhood of $\overline{pq}$. But such description is necessary for the
purpose of induction.

\begin{prop}\label{prop5a}
Let $X$ be a flat family of analytic surfaces over disk $\Delta_t$
whose central fiber $X_0 = R_1\cup R_2$ where $R_i\isom \Delta\times\P^1$
for $i=1,2$.
Suppose that $R_1\cap R_2 = \Delta(p)\cup
\Delta(q)$ where $\Delta(p)$ and $\Delta(q)$ are disks centering at points
$p$ and $q$, respectively, and $\Delta(p)$ and $\Delta(q)$ are closed
subschemes of $R_i$ ($i=1,2$). 
Suppose that $X$ is locally given by
$\Spec \BC[[x, y, z, t]] / (xy - t^a z^b)$ and $\Spec \BC[[x, y, z,
t]]/(xy - t^a)$ at $p$ and $q$, respectively, where $a, b\in \BZ$,
$a > 0$ and $b = 0$ or $1$.

Let $z\in \Gamma(\CO_X)$ and let $z = 0$ cut out the ``banana'' curve
$\overline{pq} = \overline{pq}^{(1)}\cup \overline{pq}^{(2)}$ on $X_0$, where
$\overline{pq}^{(i)} \isom \P^1$, $\overline{pq}^{(i)}\subset R_i$
and each $\overline{pq}^{(i)}$ meets $\Delta(p)$ and $\Delta(q)$
transversely at $p$ and $q$, respectively, for $i = 1,2$.

Let $\CO_X(C)$ be a line bundle on $X$
such that the restrictions of $\CO_X(C)$ to $\overline{pq}^{(i)}$ are
$\CO_{\P^1}(1)$
and let $s_1$ and $s_2$ be two global sections of $\CO_X(C)$ which generate
$H^0(\CO_{\P^1}(1))$ when restricted to $\overline{pq}^{(i)}$ for $i=1,2$.

Let $f(s_1, s_2, z, t)\in \BC[s_1, s_2, [z, t]]$ lie in
\[
f(s_1, s_2, z, t)\in \left\{\sum_{i=0}^\infty \left(
a_i(t) s_1 z^i + b_i(t) s_2 z^i\right): a_i(t), b_i(t)\in \BC[[t]]\right\} 
\]
and $f(s_1, s_2, z, 0) \ne 0$. Suppose that $f(s_1, s_2, z, t)$ is
irreducible in $\BC[s_1, s_2, [z]]\tensor_\BC \overline{\BC[[t]]}$. 

Let $\Upsilon\subset X$ be the subscheme of $X$ cut out by $f(s_1, s_2, z,
t) = 0$. Obviously, $\Upsilon_0$ consists of a multiple of $\overline{pq}$
and two disks $\Gamma_1\subset R_1$ and $\Gamma_2\subset R_2$. Let $\mu$ be
the multiplicity of $\overline{pq}$ in $\Upsilon_0$. Suppose that $\mu > 0$.

Let $\wt{\Upsilon}$ be the nodal reduction of $\Upsilon$. Then

\begin{enumerate}
\item $\wt{\Upsilon}_0$ has at most two connected
components with $\Gamma_1$ and $\Gamma_2$ in each component, respectively;
\item $\wt{\Upsilon}_0$ is connected if $b = 0$ or $p, q\not\in\Gamma_i$
for $i=1,2$;
\item if $\wt{\Upsilon}_0$ has two connected components and $\Gamma_1\ne
\Gamma_2$, then
\begin{equation}\label{prop5a-e-1}
\delta(\Upsilon_t) = 2\mu + \Gamma_1 \cdot \Delta(p) + \Gamma_1 \cdot
\Delta(q),
\end{equation}
where the intersections are taken on $R_1$ (let $\Gamma_1 \cdot \Delta(p) =
0$ if $p\not\in \Gamma_1$ and $\Gamma_1 \cdot \Delta(q) = 0$ if
$q\not\in\Gamma_1$). 
\end{enumerate}
\end{prop}

\begin{proof}
Since $\overline{pq}^{(i)}$ (\(i = 1,2\)) meet $\Delta(p)$ transversely,
we can choose the local coordinates $(x, y, z)$ of $X$ at $p$ such that
\[
X\isom \Delta_{xyzt}^4 / \big(xy = t^a (z^b + \alpha(t))\big)
\]
at $p$ and the local function $z$ is exactly the restriction of the global
$z$ as defined in the proposition, where $\alpha(t) \in \BC[[t]]$ and
$\alpha(0) = 0$. Of course, if $b = 0$, we can make $\alpha(t)$
vanish. If $b = 1$, we can replace the global function
$z$ by $z - \alpha(t)$, which we will call a translation on $z$. So
eventually, we arrive at
\begin{equation}\label{prop5a-e0}
X\isom \Delta_{xyzt}^4 / \big(xy = t^a z^b\big)
\end{equation}
at $p$ and the local coordinates $(x, y, z)$ are chosen such that the local
function $z$ is the restriction of the global function $z$.

We may put the defining equation of $\Upsilon$ in the following form
\begin{equation}\label{prop5a-e1}
\begin{split}
f(w, z, t) &= w\prod_{i = 1}^\mu 
\bigg(z + a_i(t)\bigg)\\
&\quad + t^\beta u(z, t) \prod_{j = 1}^m
\bigg(z + b_j(t)\bigg) = 0
\end{split}
\end{equation}
where $w =(c_{11} s_1 + c_{12} s_2)/ (c_{21} s_1 + c_{22} s_2)$ for some
\[
\begin{pmatrix}  
c_{11} & c_{12}\\
c_{21} & c_{22}
\end{pmatrix}
\in\text{SL}_2(\BC[[z, t]]),
\]
$u(z, t)\in \BC[[z, t]]$, $u(0, 0) \ne 0$,
$a_i(t), b_j(t)\in \overline{\BC[[t]]}$, $a_i(0) = b_j(0) = 0$, 
$\beta > 0$ and $\mu > m$. 
And we arrange $\{a_i(t)\}$ and $\{b_j(t)\}$ in the order that
\begin{equation}\label{prop5a-e2a}
\nu(a_1(t))\le \nu(a_2(t))\le ...\le \nu(a_\mu(t))
\end{equation}
and
\begin{equation}\label{prop5a-e2b}
\nu(b_1(t))\le \nu(b_2(t))\le ...\le \nu(b_m(t)),
\end{equation}
where $\nu(c(t))$ is
the valuation of $c(t)\in\overline{\BC[[t]]}$ and we let $\nu(0) = \infty$.
If $b = 0$, we may further do a translation on $z$ and assume that
\begin{equation}\label{prop5a-e2}
\sum_{i=1}^\mu a_i(t) = 0.
\end{equation}

It is obvious that \(f(w, z, t)\) is irreducible
in $\BC[w, [z]]\times_\BC \overline{\BC[[t]]}$ since $f(s_1, s_2, z, t)$
is irreducible in $\BC[s_1, s_2, [z]]\tensor_\BC \overline{\BC[[t]]}$.

We want to blow up $X$ along a subscheme cut out by $z = t^\gamma = 0$ for
some $\gamma$. Let
\begin{equation}\label{prop5a-3}
\gamma = \min\left\{ \nu(a_1(t)), \frac{\beta}{\mu - m}, \gamma_0\right\}
\end{equation}
where
\[
\gamma_0 = \min_{1\le k\le m} 
\frac{1}{\mu-m+k}\left(\beta + \sum_{j = 1}^k \nu(b_j(t))\right).
\]
Alternatively, we may define $\gamma$ as the largest number such that
\[
t^{\mu \gamma} \big| f(w, t^\gamma z, t).
\]

After a base change, we may assume that $\gamma\in\BZ$. 

Let $\wt{X}$ be the blowup of $X$ along the subscheme cut out by $z =
t^\gamma = 0$ and let $\Upsilon'\subset \wt{X}$ be the proper transform of
$\Upsilon$. The central fiber $\wt{X}_0$ consists of four surfaces $R_1\cup
R_2 \cup (\overline{pq}^{(1)}\times \P^1)\cup (\overline{pq}^{(2)}\times
\P^1)$, where $R_i \cap (\overline{pq}^{(i)}\times \P^1) =
\overline{pq}^{(i)}$ ($i = 1,2$) and $(\overline{pq}^{(1)}\times \P^1) \cap
(\overline{pq}^{(2)}\times \P^1) = (\{p\}\times \P^1) \cup (\{q\}\times \P^1)$.

By \eqref{prop5a-e1}, $\Upsilon'$ is cut out on $\wt{X}$ by
\begin{equation}\label{prop5a-e4a}
f_1(w, z_1, t) = 0
\end{equation}
where $z_1 = z/t^\gamma$ and 
\begin{equation}\label{prop5a-e4b}
f_1(w, z_1, t) = t^{-\mu\gamma} f(w, t^\gamma z_1, t).
\end{equation}
Let
\begin{equation}\label{prop5a-e4}
f_1(w, z, 0) = \bigg(w h_1(z_1) + h_2(z_1)\bigg) 
\prod_{i=1}^l (z_1 + \alpha_i)^{\mu_i}
\end{equation}
with $z_1 = z/t^\gamma$, where $h_1(z_1), h_2(z_1)\in \BC[z_1]$, $\deg
h_1(z_1) > \deg h_2(z_1)$, $w h_1(z_1) + h_2(z_1)$ is irreducible in
$\BC[w, z_1]$,
$\alpha_1, \alpha_2, ..., \alpha_l\in \BC$ are
$l$ distinct numbers, $\mu_i\in \BZ$, $\mu_i > 0$ and
$\deg h_1(z_1) + \sum_{i=1}^l \mu_i = \mu$.

By \eqref{prop5a-e4}, we see that the central
fiber $\Upsilon_0'$ of $\Upsilon'$
consists of components $\Gamma_1, \Gamma_2, I_0, I_1, ..., I_l$ where

\begin{enumerate}
\item each $I_j$ has two components $I_j^{(1)}$ and $I_j^{(2)}$ with
$I_j^{(k)}\subset \overline{pq}^{(k)}\times \P^1$ for $j = 0, 1, ..., l$
and $k = 1, 2$;
\item $I_0^{(k)}$ is cut out by $w h_1(z_1) + h_2(z_1) = 0$ on
$\overline{pq}^{(k)}\times \P^1$; $I_j^{(k)}$ is cut out by $z_1 + \alpha_j =
0$ on $\overline{pq}^{(k)}\times \P^1$ for $j = 1,2,...,l$ and $k = 1,2$
(here we regard $(w, z)$ as the affine coordinates of
$\overline{pq}^{(k)}\times \P^1\isom\P^1 \times \P^1$);
\item $I_0^{(k)}$ projects to $\overline{pq}^{(k)}$ with a degree $\mu_0 =
\deg h_1(z_1)$ map for $k = 1,2$
(if $\mu_0 = 0$, $I_0^{(k)}$ contracts to the point 
$\Gamma_k\cap \overline{pq}^{(k)}$);
\item $I_j^{(k)}$ has multiplicity $\mu_j$ in $\Upsilon_0'$; $I_j^{(1)}$ and
$I_j^{(2)}$ meet at two points $p_j\in \{p\}\times \P^1$ and $q_j\in
\{q\}\times \P^1$ for $j = 1,2,...,l$;
\item $I_1, I_2, ..., I_l$ are disjoint from each other; $I_0^{(k)}$ meets
each $I_j^{(k)}$ at exactly one point for $j=1,2,..., l$ and $k = 1,2$;
\item $I_0^{(k)}$ meets $\Gamma_k$ at a point
$r_k\in \overline{pq}^{(k)}\times \P^1$ with coordinates $w = 1/z_1 = 0$
and $r_k\not\in I_j$ for $j=1,2,...,l$ and $k = 1,2$.
\end{enumerate}

By \lemref{lem1}, the way in which $I_0^{(k)}$ are connected to
$\Gamma_k$ ($k = 1,2$) on $\wt{\Upsilon}_0$ can be described as follows.

\begin{enumerate}
\item[$(*)$]
Either $I_0^{(k)}$ and $\Gamma_k$ are joined by curves contracting to
$r_k$ on $\wt{\Upsilon}_0$
or $I_0^{(k)}$ and $\Gamma_{3-k}$ are joined by curves contracting to
$r_k$ on $\wt{\Upsilon}_0$ for $k= 1, 2$
(the latter could happen when $r_1 = r_2\in I_0^{(1)}\cap I_0^{(2)}$).
\end{enumerate}

We will argue by induction on the pair $(\mu, b)$. We define $(\mu, b) < (\mu',
b')$ if $\mu < \mu'$ or $\mu = \mu'$ and $b < b'$.

\medskip
\noindent\underline{$l > 0$}
\medskip

Take a component $I_j$ and an analytic neighborhood
$U$ of $\wt{X}$ around $I_j$. Let $Y = U\cap \Upsilon'$. Then $Y$ and $U$
have all the properties described in the proposition. For example, $Y$ is
cut out on $U$ by $f_2(w, z_2, t) = 0$ where
\begin{equation}\label{prop5a-e5}
f_2(w, z_2, t) = f_1(w, z_2 - \alpha_j, t) = t^{-\mu\gamma} f(w, t^\gamma
(z_2 - \alpha_j), t)
\end{equation}
By \lemref{lema1}, \(f_2(w, z_2, t)\) is irreducible
in $\BC[w, [z_2]]\tensor_\BC \overline{\BC[[t]]}$ 
since \(f(w, z, t)\) is irreducible in \(\BC[w, [z]]\tensor_\BC
\overline{\BC[[t]]}\).
And in the neighborhoods of $p_j$ and $q_j$, $U$ is given by
$\Spec\BC[[x, y, z, t]]/(xy - t^{(a+\gamma)} z^c)$ and 
$\Spec\BC[[x, y, z, t]]/(xy - t^{(a+\gamma)})$, respectively, where $c = 0$
if either $b = 0$ or $\alpha_j \ne 0$ and $c = 1$ otherwise.
The central fiber $Y_0$ of $Y$ contains $I_j$ with multiplicity $\mu_j$
plus two disks $\Gamma_1'$ and $\Gamma_2'$ attached to $I_j^{(1)}$ and
$I_j^{(2)}$, respectively.
Obviously, $\Gamma_k'$ is a piece of $I_0^{(k)}$ for $k=1,2$.

To apply the induction hypothesis, we have to check that $(\mu_j, c) < (\mu,
b)$, which is an easy consequence of \eqref{prop5a-e2} and the way we choose
the number $\gamma$.

Let $\wt{Y}$ be the nodal reduction of $Y$. By the induction hypothesis,
$\wt{Y}_0$ has at most two connected components with $\Gamma_1'$ and
$\Gamma_2'$ in each component. So combining with $(*)$, we see that
$\wt{\Upsilon}_0$ has at most two connected components with $\Gamma_1$ and
$\Gamma_2$ in each component. Furthermore, as long as there exists one
$I_j$ such that the corresponding $\wt{Y}_0$ is connected,
$\wt{\Upsilon}_0$ is connected.
By the induction hypothesis, $\wt{Y}_0$ is connected if $c = 0$
or $I_0$ meets $I_j$ at points other than $q_j$ for $1\le j
\le l$. Therefore, $\wt{\Upsilon}_0$ is connected
if one of the following holds:

\begin{enumerate}
\item $b = 0$;
\item $\alpha_j \ne 0$ for some $1\le j \le l$;
\item $I_0$ meets $I_j$ at points other than $p_j$ and $q_j$ 
for some $1\le j \le l$.
\end{enumerate}

Let us deal with the remaining case that
$b = l = 1$, $\alpha_1 = 0$ and $I_0$ meets $I_1$ at $p_1$ or $q_1$. 
We necessarily have $\mu_0 > 0$ due to
the way we choose the number $\gamma$. So $I_0^{(1)}$ and $I_0^{(2)}$
should meet at (at least) one point $p_0$ on $\{p\}\times \P^1$ and one
point $q_0$ on $\{q\}\times \P^1$.

If $\{p_0, q_0\}\not\subset \{r_1, r_2, p_1, q_1\}$
for some \(p_0\in I_0^{(k)}\cap (\{p\}\times \P^1)\) and 
\(q_0\in I_0^{(k)}\cap (\{q\}\times \P^1)\),
it is not hard to see by \lemref{lem1} that $I_0^{(1)}$ and
$I_0^{(2)}$ are joined by a chain of curves contracting to $p_0$ or $q_0$
(and hence to $p$ and $q$) on $\wt{\Upsilon}_0$ and $\wt{\Upsilon}_0$
is hence connected. 

If $\{p_0, q_0\} \subset \{r_1, r_2, p_1, q_1\}$ for any \(p_0\in I_0^{(k)}\cap
(\{p\}\times \P^1)\) and \(q_0\in I_0^{(k)}\cap (\{q\}\times \P^1)\), 
then $I_0^{(1)}$ and $I_0^{(2)}$ only meet at $p_0$ and $q_0$ and

\begin{enumerate}
\item[A.] $p\in \Gamma_i$ (\(i = 1,2\)), $r_1 = r_2 = p_0$ and $q_0 = q_1$;
OR
\item[B.] $q\in \Gamma_i$ (\(i = 1,2\)), $r_1 = r_2 = q_0$ and $p_0 = p_1$.
\end{enumerate}

We need to prove \eqref{prop5a-e-1} in both cases if we further assume that
$\wt{\Upsilon}_0$ has two connected components and $\Gamma_1\ne \Gamma_2$.
Let us work on case A and case B will follow from the same argument.

By \lemref{lema2} and the induction hypothesis,
\[
\begin{split}
\delta(\Upsilon_t) &= \delta(\Upsilon_t') = \delta(\Upsilon_t', p_0) +
\delta(\Upsilon_t', I_1)\\
&= (\mu_0 + \Gamma_1\cdot \Delta(p)) + (\mu_0 + 2\mu_1)\\
&= 2(\mu_0 + \mu_1) + \Gamma_1\cdot \Delta(p) = 2\mu + \Gamma_1\cdot
\Delta(p).
\end{split}
\]

\medskip
\noindent{\underline{$l = 0$}}
\medskip

Since $\mu_0 = \mu > 0$,
$I_0^{(1)}$ and $I_0^{(2)}$ meet at (at least) one point $p_0\in
\{p\}\times\P^1$ and one point $q_0\in \{q\}\times \P^1$.

If $b = 0$ or $q_0\not\in \{r_1, r_2\}$ for some $q_0 \in I_0^{(k)}\cap
(\{q\}\times \P^1)$, then by \lemref{lem1}, $I_0^{(1)}$ and $I_0^{(2)}$ are
joined by a chain of curves contracting to either $p_0$ or $q_0$ on
$\wt{\Upsilon}_0$ and $\wt{\Upsilon}_0$ is hence connected.

If $b = 1$, $q_0 = r_1 = r_2$ and $I_0^{(k)}$ meets $\{p\}\times\P^1$ at (at
least) two different points $p_0$ and $p_0'$, then by \lemref{lem1},
$I_0^{(1)}$ and $I_0^{(2)}$ are joined by a
chain of curves contracting to either $p_0$ or $p_0'$ on $\wt{\Upsilon}_0$ and
$\wt{\Upsilon}_0$ is hence connected.

If $b = 1$,
$q_0 = r_1 = r_2$, $I_0^{(k)}$ meets $\{p\}\times\P^1$ only at one point
$p_0$ and $\mu \ge 2$, then it follows from \coref{cor3} that $I_0^{(1)}$ and
$I_0^{(2)}$ are joined by a chain of curves contracting to $p_0$ on
$\wt{\Upsilon}_0$ and $\wt{\Upsilon}_0$ is hence connected.

So the only case left is that $b=1$, $q_0 = r_1 = r_2$ and $\mu = 1$.
Obviously, it follows from $(*)$ that $\wt{\Upsilon}_0$ has at most two
connected components with $\Gamma_1$ and $\Gamma_2$ in each component,
respectively. We need to verify \eqref{prop5a-e-1} if we further assume that
$\wt{\Upsilon}_0$ has two connected components and $\Gamma_1 \ne \Gamma_2$.

By \lemref{lema2},
\[
\begin{split}
\delta(\Upsilon_t) &= \delta(\Upsilon_t, p) + \delta(\Upsilon_t, q)\\
&= 1 + (1 + \Gamma_1 \cdot \Delta(q)) = 2\mu + \Gamma_1 \cdot \Delta(q).
\end{split}
\]
\end{proof}

\propref{prop5} follows more or less directly from \propref{prop5a}.

\begin{proof}[Proof of \propref{prop5}]
Suppose that $H^0(\CO_X(C+ kF))$ is generated by $g+1$ global sections
$Y_0, Y_1, ..., Y_g$ as a free $\BC[[t]]$ module.
Let $H^0(\CO_X(F))$ be generated by two global sections $Z_0$ and $Z_1$
where $Z_1 = 0$ cuts out $\overline{pq}$ on $X_0$,
Let $W_1, W_2\in H^0(\CO_X(F))$ be two global sections of $\CO_X(C+F)$
whose restrictions to $\overline{pq}^{(i)}$ generate $H^0(\CO_{\P^1}(1))$.
All these $Y_0, Y_1, ..., Y_g, Z_0, Z_1, W_1, W_2$ exist after a base
change.

Let $s_1 = W_1/Z_0$ and $s_2 = W_2/Z_0$. When restricted to an analytic
neighborhood of $\overline{pq}$, $s_1$ and $s_2$ are holomorphic sections
of $\CO_X(C)$ and
\[
\begin{split}
& \quad\left\{\sum_{i=0}^g a_i(t) Y_i: a_i(t)\in\BC[[t]]\right\}\\
& \subset \left\{\sum_{i=0}^k (b_i(t) s_1 Z_0^i Z_1^{k-i} +
c_i(t) s_2 Z_0^i Z_1^{k-i}): b_i(t), c_i(t)\in\BC[[t]]\right\}.
\end{split}
\]

Therefore, $\Upsilon$ is locally cut out by $f(s_1, s_2, z, t) = 0$ as
described in \propref{prop5a}, where $z = Z_1/Z_0$. Since $\Upsilon_t$ is
(geometrically) irreducible, it follows from \lemref{lema1} that $f(s_1,
s_2, z, t)$ is irreducible in $\BC[s_1, s_2, [z]]\tensor_\BC
\overline{\BC[[t]]}$. Then by \propref{prop5a}, $\Gamma_1$ and $\Gamma_2$
must pass through either $p$ or $q$; otherwise, $\Gamma_1$ and $\Gamma_2$
will be joined by a chain of curves on $\wt{\Upsilon}_0$ whose images lie
in $\overline{pq}$ and we have shown in the proof of \propref{prop1b} that
$\Gamma_1$ and $\Gamma_2$ are joined by a chain of curves somewhere else on
$\wt{\Upsilon}_0$, which leads to a contradiction.

If $\Gamma_1$ and $\Gamma_2$ pass through $p$ or $q$, then by
\propref{prop5a},
\[
\begin{split}
\delta(\Upsilon_t, \overline{pq}) &= 2\mu + \intsc(\Gamma_1, E; p) +
\intsc(\Gamma_1, E; q)\\
&= \intsc(\Sigma_1, E; p) + \intsc(\Sigma_1, E; q)
\end{split}
\]
where $\mu$ is the multiplicity of $\overline{pq}$ in $\Upsilon_0$.
\end{proof}
 
\begin{prop}\label{prop5b}
Let $X\subset \Delta_{xyzt}^4$ be defined by
\begin{equation}\label{prop5b-e0}
y(y+x^2 + t^a z^b) = \lambda t^c
\end{equation}
where $a, b, c\in \BZ$, $a + b > 0$, $b = 0$ or $1$, $c > 0$,
$\lambda = \lambda(x,y, z,t)\in \BC[[x, y, z, t]]$ and
$\lambda(0, 0, 0, 0)\ne 0$.

Let $f(w, z, t)\in \BC[w, [z, t]]$ lie in
\[
f(w, z, t)\in \left\{\sum_{i=0}^\infty \left(
a_i(t) z^i + b_i(t) w z^i\right): a_i(t), b_i(t)\in \BC[[t]]\right\} 
\]
and $f(w, z, 0) = wz^\mu$ for some $\mu\in\BZ$ and $\mu\ge 0$.
Suppose that $f(w, z, t)$ is
irreducible in $\BC[w, [z]]\tensor_\BC \overline{\BC[[t]]}$.

Let $\Upsilon\subset X$ be a flat family of curves cut out by $f(w, z, t) =
0$ on $X$, where $w = w(x, y, z, t)\in \BC[[x, y, z, t]]$ satisfying
\[
w(0, 0, 0, 0) = 0 \text{ and } \frac{\partial w}{\partial x}(0, 0, 0, 0)\ne 0.
\]

Let $I_1, I_2, \Gamma_1, \Gamma_2$ be the irreducible components of
$\Upsilon_0$ where $I = I_1\cup I_2$ is cut out by $z = 0$ and $\Gamma_1
\cup \Gamma_2$ is cut out by $w = 0$ on $X_0 = R_1\cup R_2$. 
Let \(\Gamma_j\subset R_j\) and \(I_j\subset R_j\) for \(j = 1,2\). 
Suppose that $\Gamma_j$ meets \(E\) properly for \(j = 1,2\).
Let $\wt{\Upsilon}$ be the nodal reduction of $\Upsilon$. 
Then either $\Gamma_1$ and $\Gamma_2$ lie on the same
connected component of $\wt{\Upsilon}_0$ or the dual graph of
$\wt{\Upsilon}_0$ contains a circuit.
\end{prop}

\begin{proof}
If $a = 0$, the conclusion is more or less obvious. Notice that $b = 1$ if $a =
0$. Hence $\Gamma_j$ meets $E$ transversely, while $I_j$ is tangent to $E$
with multiplicity $2$. So it follows from \lemref{lem1} that $\Gamma_1$ and
$\Gamma_2$ lie on the same connected component of $\wt{\Upsilon}_0$.

It is also trivial if $\mu = 0$.

Suppose that $a > 0$ and $\mu > 0$. Again, we will
argue by induction on $(\mu, b)$.
We need to blow up $\Upsilon$ three times along subschemes on the central
fiber in order to ``bring down'' the pair $(\mu, b)$.

We may put $f(w, z, t) = 0$ in the form \eqref{prop5a-e1} and also
assume $\{a_i(t)\}$ and $\{b_j(t)\}$ to satisfy \eqref{prop5a-e2a} and
\eqref{prop5a-e2b} and to satisfy \eqref{prop5a-e2} if $b = 0$.

Since \(\Gamma_i\) meets \(E\) properly, \(w(0, 0, z, 0) \ne 0\). Let
\[
w(0, 0, z, t) = \sum_{i=1}^\eta (z + c_i(t))
\]
where \(\eta > 0\), \(c_i(t)\in \overline{\BC[[t]]}\) and \(c_i(0) =
0\). We arrange \(\{c_i(t)\}\) in the order that
\[
\nu(c_1(t)) \le \nu(c_2(t)) \le ...\le \nu(c_\eta(t))
\]

Let
\begin{equation}\label{prop5b-e2}
\gamma = \min \left\{ \frac{a}{2\eta-b},\ \frac{c}{4\eta},
\ \nu(c_1(t)),\ \nu(a_1(t)),\ \frac{\beta}{\mu- m +\eta},\ \gamma_0\right\},
\end{equation}
where
\[
\gamma_0 = \min_{1\le k\le m} 
\frac{1}{\mu-m+k+\eta}\left(\beta + \sum_{j = 1}^k \nu(b_j(t))\right).
\]
Alternatively, $\gamma$ can be defined as the largest number such that
\[
t^{(\mu+\eta)\gamma} \big | f\big(w(t^{\eta \gamma} x, t^{2\eta \gamma} y,
t^\gamma z, t),
t^\gamma z , t\big) \text{ and }
t^{4\eta\gamma} \big| g(t^{\eta \gamma} x, t^{2\eta \gamma} y, t^\gamma z, t),
\]
where we let $g(x, y, z, t) = y(y+x^2 + t^a z^b) - \lambda t^c$.

Let $\gamma\in \BZ$ after a base change. Let $\Upsilon'$ be the blowup of
$\Upsilon$ at the 0-dimensional subscheme $x = y = z = t^\gamma = 0$. Then
$\Upsilon'$ is given by
\begin{equation}\label{prop5b-e3}
y_1\left(y_1+t^\gamma x_1^2 + t^{a - (1-b)\gamma} z_1^b\right) = 
\lambda t^{c - 2\gamma}
\end{equation}
and
\begin{equation}\label{prop5b-e4}
f_1(w_1, z_1, t) = 0
\end{equation}
where $x_1 = x/t^\gamma$, $y_1 = y/t^\gamma$, $z_1 = z/t^\gamma$,
\[
w_1 = w_1(x_1, y_1, z_1, t) = t^{-\gamma} w(t^\gamma x_1, t^\gamma y_1, 
t^\gamma z_1, t)
\]
and
\[
f_1(w_1, z_1, t) = t^{-(\mu+1)\gamma} f(t^\gamma w_1, t^\gamma z_1, t).
\]

By \eqref{prop5b-e3} and \eqref{prop5b-e4},
we can describe the central fiber $\Upsilon_0'$ of $\Upsilon$ as follows.

The exceptional locus of $\Upsilon'\to \Upsilon$ is a reducible and
nonreduced curve $F$ cut out on $\P^3$ by
\begin{equation}\label{prop5b-e5}
\left\{
\begin{split}
& y_1^2 = 0\\
& f_1\left(w_1(x_1, y_1, z_1, 0), z_1, 0\right) = 0
\end{split}
\right.
\end{equation}
where \(x_1 = X_1/T_1\), \(y_1 = Y_1/T_1\) and 
\(z_1 = Z_1/T_1\) are the affine coordinates of $\P^3$ with
corresponding homogeneous coordinates \( (X_1, Y_1, Z_1, T_1) \).

It is not hard to see that $\Gamma_1\cup \Gamma_2$ meets $F$ at the point
\(p\) with coordinates
\[
\left(-\frac{\partial w}{\partial z}(0,0,0,0), 0, \frac{\partial
w}{\partial x}(0,0,0,0), 0\right)
\]
and $I_1\cup I_2$ meets $F$ at the point \(q\) with coordinates
\((1, 0, 0, 0)\). Obviously, $p\ne q$, i.e.,
the blowup $\Upsilon'\to\Upsilon$ has
separated \(\Gamma_1\cup \Gamma_2\) from \(I_1\cup I_2\).

The point \(p = \Gamma_k\cap F\) lies on a unique irreducible component
\(\Sigma\) of \(F\) and \(\Sigma\) has multiplicity 2 in \(F\).

We may continue to use \(\wt{\Upsilon}\) to denote the nodal reduction of
\(\Upsilon'\). By \lemref{lem1}, each \(\Gamma_i\) is joined to a component
dominating \(\Sigma\) on \(\wt{\Upsilon}_0\) by a chain of curves
contracting to the point \(p\) for \(i=1,2\). Therefore, in order
to show that \(\Gamma_1\) and \(\Gamma_2\) lie on the same connected
component of \(\wt{\Upsilon}_0\), it suffices to show that
\begin{enumerate}
\item[\((*)\)] all the components dominating \(\Sigma\) lie on the same
connected component of \(\wt{\Upsilon}_0\).
\end{enumerate}
This line of argument naturally leads to the second and the third blowups.
This time we need to blow up \(\Upsilon'\) along some subscheme supported along
\(F_\red\).

Let $\Upsilon''$ be the blowup of $\Upsilon'$ along the subscheme cut out
by \(x_1 = t^{(\eta -1)\gamma} = 0\). And let \(\Y\) be the blowup of
\(\Upsilon''\) along the subscheme cut out by \(y_1 = t^{(2\eta -1)\gamma}
= 0\). Notice that if \(\eta = 1\), we do not need the intermediate family
\(\Upsilon''\). Finally, we obtain \(\Y\), which is given by
\begin{equation}\label{prop5b-e6}
y_2\left(y_2+ x_2^2 + t^{a - (2\eta-b)\gamma} z_1^b\right) = 
\lambda t^{c - 4\eta\gamma}
\end{equation}
and
\begin{equation}\label{prop5b-e7}
f_2(w_2, z_1, t) = 0
\end{equation}
where \(x_2 = x_1 / t^{(\eta -1)\gamma}\), \(y_2 = y_1 / t^{(2\eta
-1)\gamma}\),
\[
w_2 = w_2(x_2, y_2, z_1, t) = \frac{1}{t^{(\eta -1)\gamma}} 
w_1(t^{(\eta -1)\gamma} x_2, t^{(2\eta -1)\gamma} y_2, z_1, t)
\]
and
\[
f_2(w_2, z_1, t) = \frac{1}{t^{(\eta -1)\gamma}} f_1(t^{(\eta -1)\gamma}
w_2, z_1, t).
\] 
Let
\begin{equation}\label{prop5b-e8}
f_2(w_2, z_1, 0) = \bigg(w_2 h_1(z_1) + h_2(z_1)\bigg)
\prod_{i=1}^l (z_1 + \alpha_i)^{\mu_i} = 0
\end{equation}
where $h_1(z_1), h_2(z_1)\in \BC[z_1]$, $\deg
h_1(z_1) > \deg h_2(z_1)$, $w_2 h_1(z_1) + h_2(z_1)$ is irreducible in
$\BC[w_2, z_1]$, $\alpha_1, \alpha_2, ..., \alpha_l\in \BC$ are
$l$ distinct numbers, $\mu_i\in \BZ$, $\mu_i > 0$ and
$\deg h_1(z_1) + \sum_{i=1}^l \mu_i = \mu$.

By \eqref{prop5b-e8}, there are curves \(\wt{\Sigma} \subset \Y_0\) and 
\(J_i \subset \Y_0\) where \(\wt{\Sigma}\) is cut out by \(w_2 h_1(z_1) +
h_2(z_1) = 0\) and \(J_i\) is cut out by \(z_1 + \alpha_i = 0\) for \(i =
1,2,..., l\). Obviously,
\(\wt{\Sigma}\) dominates \(\Sigma\) with a degree two map.

To be precise, \(\wt{\Sigma}\) and \(J_i\) are complete curves lying on the
surface \(\P\left(\CO_{\P^1}\oplus \CO_{\P^1}(1)\oplus
\CO_{\P^1}(1)\right)\). However, to argue $(*)$, we only need to study
their affine parts. So we will treat them as affine curves in 
$\A^3 = \Spec \BC[x_2, y_2, z_1]$. 
For example, when we talk about the intersections between these
curves, we are talking about their intersections in $\A^3$.

\medskip
\noindent\underline{$c = 4\eta\gamma$}
\medskip

We claim that \(\wt{\Sigma}\) is
irreducible in this case and then \((*)\) will follow immediately.

By \eqref{prop5b-e6} and \eqref{prop5b-e8}, \(\wt{\Sigma}\) is given by
\begin{equation}\label{prop5b-e9}
\left\{
\begin{split}
& y_2 (y_2 + x_2^2 + \lambda_1 z_1^b) = \lambda_0\\
& w_2(x_2, 0, z_1, 0) h_1(z_1) + h_2(z_1) = 0
\end{split}
\right.
\end{equation}
where \(\lambda_1 = 0\) or \(1\) and \(\lambda_0 = \lambda(0, 0, 0, 0)\ne
0\).

Notice that \(w_2(0, 0, z_1, 0)\) is a degree \(\eta\) polynomial in
\(z_1\) due to the way we choose \(\gamma\). So if we solve \(w_2(x_2, 0,
z_1, 0) h_1(z_1) + h_2(z_1) = 0\) for \(x_2\), we obtain \(x_2\in \BC(z_1)\)
and \(x_2\not\in \BC\).

It is not hard to see that
\(\wt{\Sigma}\) is irreducible if \( (x_2^2 + \lambda_1 z_1^b)^2 + \lambda_0
\) is a perfect square in \(\BC(z_1) \)
for some \( x_2 \in \BC(z_1) \) and \( x_2 \) nonconstant.
Such \( x_2 \) does not exist by \lemref{lema3}.

If \( c > 4 \eta \gamma \), then each of
\(\wt{\Sigma}\) and \(J_i\) have exactly two irreducible components.
Let \(\wt{\Sigma} = \wt{\Sigma}_1 \cup \wt{\Sigma}_2\) and \(J_i =
J_i^{(1)}\cup J_i^{(2)}\) for \(i = 1,2,..., l\). 

\medskip
\noindent\underline{\( c > 4 \eta \gamma \) and \( a = (2\eta -b)\gamma
\)}
\medskip

If $b = 0$, $J_i^{(1)}$ and $J_i^{(2)}$ meet at two points for all
$i$. If $b = 1$, $J_i^{(1)}$ and $J_i^{(2)}$ meet at a single point if and
only if the corresponding \(\alpha_i = 0\) and this point must be \(r =
(x_2 = y_2 = z_1 = t = 0)\). It is not hard to see that $\wt{\Sigma}_1$ and
$\wt{\Sigma}_2$ meet at (at least) one point other than $r$.

Let $G$ be the dual graph of the following components of $\wt{\Y}_0$ (let
$\wt{\Y}$ be the nodal reduction of $\Y$)

\begin{enumerate}
\item $\wt{\Sigma}_1$ and $\wt{\Sigma}_2$;
\item the components dominating $J_i^{(1)}$ or $J_i^{(2)}$
for $J_i^{(1)}$ and $J_i^{(2)}$ that meet at two points;
\item the contractible components which contract to a point in
\( (\wt{\Sigma}_1\cap \wt{\Sigma}_2)\cup (J_1^{(1)}\cap J_1^{(2)}) 
\cup (J_2^{(1)}\cap J_2^{(2)}) \cup
...\cup (J_l^{(1)}\cap J_l^{(2)}) \backslash \{r\} \).
\end{enumerate}

Then \lemref{lem1} tells us that $\deg([\wt{\Sigma}_k]) \ge 1$ for $k=1,2$
and all the other vertices of $G$ has degree at least two. Therefore,
either $[\wt{\Sigma}_1]$ and $[\wt{\Sigma}_2]$ lie on the same component of
$G$ or $G$ contains a circuit.

Let \( \mu_0 = \deg h_1(z_1) \).

\medskip
\noindent\underline{\( c > 4 \eta \gamma \), \( a > (2\eta -b)\gamma \) and 
either \(\mu_0 > 0\) or \(l > 1\)}
\medskip

Let $p_1 \in \wt{\Sigma}_1\cap \wt{\Sigma}_2$. If no $J_i$ passes through
$p_1$, then $(*)$ follows directly from \lemref{lem1}. Otherwise suppose
that $p_1\in J_1$. By \eqref{prop5b-e6} and \eqref{prop5b-e7}, 
$\Y$ is locally defined by  
\begin{equation}\label{prop5b-e10}
y \left(y + x^2 + t^{a - (2\eta-b)\gamma} z^{b'}\right) = 
\lambda t^{c - 4\eta\gamma}
\end{equation}
and
\begin{equation}\label{prop5b-e11}
f_3(w_2, z, t) = 0
\end{equation}
at $p_1$, where \( f_3(w_2, z, t) = f_2(w_2, z - \alpha_1, t) \), $b' =
0$ if $b = 0$ or $\alpha_1\ne 0$ and $b' = 1$ otherwise. Notice that
$f_3(w_2, z, t)$ is irreducible in $\BC[w_2, [z]]\tensor_\BC
\overline{\BC[[t]]}$ by \lemref{lema1}
and $(\mu_1, b') < (\mu, b)$. Apply the induction
hypothesis and we are done.

\medskip
\noindent\underline{\( c > 4 \eta \gamma \), \( a > (2\eta -b)\gamma \), 
\(\mu_0 = 0\), \(l = 1\) and \(\alpha_1 \ne 0\)}
\medskip

Notice that this case happens only if $b = 1$ due to \eqref{prop5a-e2}.
Following the argument of the previous case, we observe that $b' = 0$ in
\eqref{prop5b-e10}. So $(\mu, b') < (\mu, b)$ and the induction hypothesis
still applies.

\medskip 
\noindent\underline{\( c > 4 \eta \gamma \), \( a > (2\eta -b)\gamma \), 
\(\mu_0 = 0\), \(l = 1\) and \(\alpha_1 = 0\)}
\medskip

We necessarily have $\gamma = \nu(c_1(t))$ in this case. So there exists
$p_1\in \wt{\Sigma}_1\cap \wt{\Sigma}_2$ such that $p_1\not\in J_1$. Then
$(*)$ follows immediately from \lemref{lem1}.
\end{proof}

\begin{prop}\label{prop5c}
Let $X$ be a flat family of analytic surfaces over disk $\Delta_t$
whose central fiber $X_0 = R_1\cup R_2$ where $R_i\isom \Delta\times\P^1$
for $i=1,2$. Suppose that $R_1\cap R_2 = \Delta(r)$ where
$\Delta(r)$ is a disk centering at point $r$ and $\Delta(r)$ is a closed
subscheme of $R_i$ ($i=1,2$). Suppose that $X$ is locally defined by
\eqref{prop5b-e0} at $r$ when embedded to $\Delta_{xyzt}^4$.

Let $z\in \CO_X$ and let $z = 0$ cut out the curve
$\overline{r} = \overline{r}^{(1)}\cup \overline{r}^{(2)}$ on $X_0$, where
$\overline{r}^{(i)} \isom \P^1$, $\overline{r}^{(i)}\subset R_i$
and each $\overline{r}^{(i)}$ meets $\Delta(r)$
at $r$ with multiplicity 2, for $i = 1,2$.

Let $\CO_X(C)$ be a line bundle on $X$
such that the restrictions of $\CO_X(C)$ to $\overline{r}^{(i)}$ are
$\CO_{\P^1}(1)$
and let $s_1$ and $s_2$ be two global sections of $\CO_X(C)$ which generate
$H^0(\CO_{\P^1}(1))$ when restricted to $\overline{r}^{(i)}$ for $i=1,2$.

Let $f(s_1, s_2, z, t)\in \BC[s_1, s_2, [z, t]]$ lie in
\[
f(s_1, s_2, z, t)\in \left\{\sum_{i=0}^\infty \left(
a_i(t) s_1 z^i + b_i(t) s_2 z^i\right): a_i(t), b_i(t)\in \BC[[t]]\right\} 
\]
and $f(s_1, s_2, z, 0) \ne 0$. Suppose that $f(s_1, s_2, z, t)$ is
irreducible in $\BC[s_1, s_2, [z]]\tensor_\BC \overline{\BC[[t]]}$. 

Let $\Upsilon\subset X$ be the subscheme of $X$ cut out by $f(s_1, s_2, z,
t) = 0$. Obviously, $\Upsilon_0$ consists of a multiple of $\overline{r}$
and two disks $\Gamma_1$ and $\Gamma_2$. Let $\mu$ be the multiplicity of
$\overline{r}$ in $\Upsilon_0$. Suppose that $\mu > 0$ and $\Gamma_i$ meets
$\Delta(r)$ properly if $\Gamma_i$ passes through $r$ for $i=1,2$.

Let $\wt{\Upsilon}$ be the nodal reduction of $\Upsilon$. Then either
$\Gamma_1$ and $\Gamma_2$ lie on the same connected component of
$\wt{\Upsilon}_0$ or the dual graph of $\wt{\Upsilon}_0$ contains a circuit.
\end{prop}

\begin{proof}
Our argument proceeds almost identically to that for \propref{prop5a}.

Just as in the proof of \propref{prop5a}, we can choose local coordinates
$(x, y, z, t)$ of $X$ at $r$ such that $X$ is defined by \eqref{prop5b-e0}
at $r$ and the local function $z$ is the restriction of the global function
$z$.

We may put the defining equation in the form \eqref{prop5a-e1} and also
assume $\{a_i(t)\}$ and $\{b_j(t)\}$ to satisfy \eqref{prop5a-e2a} and
\eqref{prop5a-e2b} and to satisfy \eqref{prop5a-e2} if $b = 0$.

Let $\gamma$ be the number defined by \eqref{prop5a-3} and let $\gamma\in
\BZ$ after a base change.

The case that $\Gamma_i$ passes through $r$ has been covered by
\propref{prop5b}. Suppose that $\Gamma_i$ does not pass through $r$ 
for $i=1,2$. Again, we will argue by induction on $(\mu, b)$).

Let $\wt{X}$ be the blowup of $X$ along the subscheme cut out by $z =
t^\gamma = 0$ and let $\Upsilon'\subset \wt{X}$ be the proper transform of
$\Upsilon$. The central fiber $\wt{X}_0$ consists of four surfaces $R_1\cup
R_2 \cup (\overline{r}^{(1)}\times \P^1)\cup (\overline{r}^{(2)}\times
\P^1)$, where $R_i \cap (\overline{r}^{(i)}\times \P^1) =
\overline{r}^{(i)}$ ($i = 1,2$) and $(\overline{r}^{(1)}\times \P^1) \cap
(\overline{r}^{(2)}\times \P^1) = \{r\}\times \P^1$.

We have the same defining equations \eqref{prop5a-e4a} for $\Upsilon'$ and
\eqref{prop5a-e4} for $\Upsilon_0'$. Let $I_0 = I_0^{(1)}\cup I_0^{(2)}, 
I_1= I_1^{(1)}\cup I_1^{(2)}, ..., I_l = I_l^{(1)}\cup I_l^{(2)}\subset
\Upsilon_0'$ be the components of $\Upsilon_0'$ defined in the same way as
in the proof of \propref{prop5a}. Let $p_j = I_j^{(1)}\cap I_j^{(2)}$.

Let $r_k$ be the intersection between $\Gamma_k$ and $I_0^{(k)}$ with
coordinates $w = 1/z_1 = 0$ on $\overline{r}^{(k)}\times\P^1\isom
\P^1\times \P^1$. By \lemref{lem1},

\begin{enumerate}
\item[$(*)$]
$I_0^{(k)}$ and $\Gamma_k$ are joined by curves contracting to
$r_k$ on $\wt{\Upsilon}_0$ for $k=1,2$.
\end{enumerate}

\medskip
\noindent{\underline{$l > 0$}}
\medskip

Take a component $I_j$ and an analytic neighborhood
$U$ of $\wt{X}$ around $I_j$. Let $Y = U\cap \Upsilon'$. Then $Y$ and $U$
have all the properties described in the proposition. 
And in the neighborhood of $p_j$, $U$ is given by
$y(y + x^2 + t^{a+\gamma} z^{b'}) = t^c$. 
where $b' = 0$
if either $b = 0$ or $\alpha_j \ne 0$ and $b' = 1$ otherwise.
The central fiber $Y_0$ of $Y$ contains $I_j$ with multiplicity $\mu_j$
plus two disks $\Gamma_1'$ and $\Gamma_2'$ attached to $I_j^{(1)}$ and
$I_j^{(2)}$, respectively.
Obviously, $\Gamma_k'$ is a piece of $I_0^{(k)}$ for $k=1,2$.

To apply the induction hypothesis, we have to check that $(\mu_j, b') < (\mu,
b)$, which is an easy consequence of \eqref{prop5a-e2} and the way we choose
the number $\gamma$. Also we need to check that $\Gamma_k'$ meets
$\{r\}\times\P^1$ properly if $\Gamma_k'$ passes through $p_j$. This is
trivially true because $\{r\}\times\P^1 \not\subset I_0$. So the induction
hypothesis applies. Combining with $(*)$, we obtain the statement of the
proposition.

\medskip
\noindent{\underline{$l = 0$}}
\medskip

Obviously, $I_0^{(1)}$ and $I_0^{(2)}$ meet at some point $p\in
\{r\}\times\P^1$. Since $\Gamma_k$ does not pass through $r$, $p\ne r_k$
for $k = 1, 2$. Therefore, $I_0^{(1)}$
and $I_0^{(2)}$ are joined by curves contracting to $p$ on
$\wt{\Upsilon}_0$. Combining with
$(*)$, we conclude that $\Gamma_1$ and $\Gamma_2$ lie on the same connected
component of $\wt{\Upsilon}_0$.
\end{proof}

\propref{prop5z} follows directly from \propref{prop5c}. We will leave the
details to the readers.

\subsection{Proof of \thmref{thm3}}

Before we proceed, we would like to raise a simple question.

Let $\pi: X\to Y$ be a proper and dominant map between two irreducible
varieties $X$ and $Y$. Assume that $Y$ is normal. Let $q\in Y$ and $p\in
\pi^{-1}(q)$ be a point on the fiber over $q$. Is it true that for any
analytic neighborhood $U$ of $p$, $\pi(U)$ contains a neighborhood $V$ of
$q$?

In general, this is not true. For example, we may take $X$ to be the blowup
of $Y = \P^2$ at a point $p$ and $q$ to be a point on the exceptional
divisor. The statement is true for $\pi$ finite. More generally, by using
Stein factorization, we have the following.

\begin{prop}[Open Mapping Principle]\label{prop6}
Let $\pi: X\to Y$ be a proper and dominant map between two irreducible
varieties $X$ and $Y$. Assume that $Y$ is normal. Let $q\in Y$ and let
$W$ be a connected component of $\pi^{-1}(q)$. Then for any analytic
neighborhood $U$ of $W$, $\pi(U)$ contains a neighborhood $V$ of $q$.
\end{prop}

Suppose that for a general point $[R^s]\in \wt{\CR_g}$, there is a point
\[
\left( [R^s], [\Upsilon_0], w_1, w_2, ..., w_\beta, \delta_1, \delta_2,
..., \delta_\beta, \nu_1, \nu_2, ..., \nu_\beta\right) \in \wt{\Z_g}
\]
lying on the fiber of $\wt{\Z_g}\to \wt{\H_g}$ over $[R^s]$ such that
$\Upsilon_0$ contains a Type III chain $\overline{r_{ij}} \cup
\overline{r_{ij} q_1}\cup \overline{q_1 q_2} \cup ...\cup \overline{q_{m-1}
q_m}$. Let $W_0$ be the index set such that
\[
\{w_1, w_2, ..., w_\beta\}\cap (\overline{r_{ij}} \cup
\overline{r_{ij} q_1}\cup \overline{q_1 q_2} \cup ...\cup \overline{q_{m-1}
q_m}) = \{w_i: i\in W_0\}.
\]
\thmref{thm3} is equivalent to saying that $\sum_{i\in W_0} \delta_i \ge
\sum_{l = 0}^m \intsc(\Sigma_1, E; q_l)$ (let $q_0 = r_{ij}$) and if the
equality holds, $\nu_i = 1$ for $i\in W_0$.

As indicated in \eqref{thm3-e3}, we will move $R^s$ in four steps.

\medskip
\noindent \underline{Step 1. $[R^s]\in \wt{\CR_g}\imply [R^s]\in \wt{\CR_g^0}$}
\medskip

When a $R^s$ in $\wt{\CR_g}$ degenerates to a general $R^s$ in $\wt{\CR_g^0}$, 
the points $r_{ij} = q_0, q_1, q_2, ..., q_m$ on the Type III chain 
``collapse'' to the point $r_{1j} = r_{2j}$. At the same time,
\begin{equation}\label{thm3-e4a}
\sum_{i=0}^m \intsc(\Sigma_1, E; q_i) = \intsc(\Sigma_1, E; r_{1j})
\end{equation}
where the $\Sigma_1$ on the LHS refers to the component of $\Upsilon_0$
lying on the old $R^s$ and the $\Sigma_1$ on the RHS lies on the new $R^s$.

Meanwhile, the limiting singularities of the old
$R^s$ which lie on the Type III chain will degenerate to the points on
$\overline{r_{1j}}\cup \overline{r_{2j}}$. 
So for $[R^s]\in \wt{\CR_g^0}$, there is a point
\[
\left( [R^s], [\Upsilon_0], w_1, w_2, ..., w_\beta, ... \right) 
\in \wt{\Z_g}
\]
lying on the fiber of $\wt{\Z_g}\to \wt{\H_g}$ over $[R^s]$ such that
$\Upsilon_0$ contains $\overline{r_{1j}}\cup \overline{r_{2j}}$ with
multiplicity $\mu$ and 
\[
\{w_1, w_2, ..., w_\beta\}\cap (\overline{r_{1j}}\cup \overline{r_{2j}})
= \{w_i: i\in W_0\}.
\]
In the proof of \propref{prop1b}, we have shown that $\Gamma_1\cup
\Gamma_2$ does not pass through $r_{1j} = r_{2j}$ on the new $R^s$ in
$\wt{\CR_g^0}$. So by \eqref{thm3-e4a},
\begin{equation}\label{thm3-e4}
2\mu = \sum_{l = 0}^m \intsc(\Sigma_1, E; q_l).
\end{equation}

\medskip
\noindent \underline{Step 2. $[R^s]\in \wt{\CR_g^0}\imply [S]\in \Y_g$}
\medskip

Let us consider the fiber $(\wt{\Z_g})_{[R^s]}$ of
$\wt{\Z_g}\to\wt{\H_g}$ over a point $[R^s]\in \wt{\CR_g^0}$. Let the
point 
\[
\left( [R^s], [\Upsilon_0], w_1, w_2, ..., w_\beta, \delta_1,
\delta_2, ..., \delta_\beta, \nu_1, \nu_2, ..., \nu_\beta\right)
\]
with the property described above 
lie on a connected component $Z$ of $(\wt{\Z_g})_{[R^s]}$. Then
for any point
\[
\left( [R^s], [\Upsilon_0'], w_1', w_2', ..., w_\beta', \delta_1',
\delta_2', ..., \delta_\beta', \nu_1', \nu_2', ..., \nu_\beta'\right)
\]
lying on the same component $Z$, we necessarily have that

\begin{enumerate}
\item $\Upsilon_0'$ contains $\overline{r_{1j}}\cup \overline{r_{2j}}$ with
multiplicity $\mu$;
\item it follows from \propref{prop3b} that
\[
\{w_1', w_2', ..., w_\beta'\}\cap (\overline{r_{1j}}\cup \overline{r_{2j}})
= \{w_i' : i\in W_0\}.
\]
\end{enumerate}

If we take an analytic neighborhood $U$ of $Z$ and project it to
$\wt{\H_g}$, by \propref{prop6}, the image $\pi(U)$ will contain a
neighborhood $V$ of $[R^s]$, where $\pi$ is the map $\wt{\Z_g}\to
\wt{\H_g}$. By \propref{prop3c}, $V$ contains general points of $\Y_g$.

For a general point $[S]\in \Y_g\cap V$,
the fiber of $U\to \pi(U)$ over $[S]$ consists of points
\[
([S], [D], w_1, w_2, ..., w_\beta, ...)
\]
where

\begin{enumerate}
\item by \propref{prop4} and \ref{prop5},
$D$ contains two connected components $D_1, D_2\in |\CO_S(F)|$ with
multiplicities $\mu_1$ and $\mu_2$, respectively, where $\mu_1 + \mu_2 =
\mu$;
\item
\[
\{w_1, w_2, ..., w_\beta\}\cap (D_1\cup D_2) = \{w_i: i\in W_0\}.
\]
\end{enumerate}

Basically, if we have a family of surfaces $S_t$ approach $R^s$, by
\propref{prop5z}, the corresponding family of curves $D_t\subset S_t$ will have
components in $|\CO_S(F)|$ with total multiplicities $\mu$ degenerating to
$\overline{r_{1j}}\cup \overline{r_{2j}}$. And by \propref{prop4}, there
are exactly two rational curves in $D_1, D_2\in |\CO_S(F)|$ in the
neighborhood of $\overline{r_{1j}}\cup \overline{r_{2j}}$. So $D$ will
contain $D_1$ and $D_2$ with a total multiplicity $\mu$.
Let $G$ be the irreducible component of $D$ in $|\CO_S(C+k'F)|$ for some
$k'$.

Let $W_1$ and $W_2$ be the index sets such that
\[
\{w_1, w_2, ..., w_\beta\}\cap D_j = \{ w_i: i\in W_j \}\text{ for } j = 1,2.
\]
Of course, $W_0 = W_1\cup W_2$ and $W_1\cap W_2 = \emptyset$.

\medskip
\noindent \underline{Step 3. $[S]\in \Y_g\imply
[R^s] \in \wt{\CR_g^0}$}
\medskip

If we degenerate an elliptic K3 surface $S$ in $\Y_g$ to a general $R^s$ in
$\wt{\CR_g^0}$ with
the corresponding $D\subset S$, described as above, degenerating to
$\Upsilon_0\subset R$, by \propref{prop3} and \ref{prop4}, we can make $D_1$
degenerate to $\overline{p_1 q_1} = \overline{p_1 q_1}^{(1)}\cup
\overline{p_1 q_1}^{(2)}$ and make $D_2$ degenerate to
$\overline{p_2 q_2} = \overline{p_2 q_2}^{(1)}\cup\overline{p_2 q_2}^{(2)}$.
So $\Upsilon_0$ will contain $\overline{p_j q_j}$ with multiplicity at
least $\mu_j$ for $j = 1,2$. Let $\Upsilon_0$ contain $\overline{p_j q_j}$ 
with multiplicity $\mu_j + \mu_j'$. Namely, $G_t\subset S_t$ will
degenerate to a curve on $R$ containing $\overline{p_j q_j}$ with multiplicity
$\mu_j'$.

Let $W_1'$ and $W_2'$ be the index sets such that
\[
\{w_1, w_2, ..., w_\beta\}\cap \overline{p_j q_j} = \{
w_i: i\in W_j\} \cup \{
w_i: i\in W_j'\}\text{ for } j = 1,2
\]
where $([R^s], [\Upsilon_0], w_1, w_2, ..., w_\beta, ...)\in \wt{\Z_g}$ is
the corresponding limit when $S$ degenerates to $R^s$.

More intuitively, those $w_i$'s for $i\in W_1'\cup W_2'$ are the extra
limiting singularities we get due to the extra multiplicities $\mu_1'$ and
$\mu_2'$. We have to control the extra $\delta$-invariant
\[
\sum_{i\in W_1'\cup W_2'} \delta_i.
\]

By \propref{prop5} and \lemref{lema2},
\begin{equation}\label{thm3-e5}
\begin{split}
\sum_{i\in W_1'\cup W_2'} \delta_i &= \delta(G_t, \overline{p_1q_1}) +
\delta(G_t, \overline{p_2q_2})\\
&= 2\mu_1' + \intsc(\Gamma_1, E; p_1) + \intsc(\Gamma_1, E; q_1)\\
&\quad + 2\mu_2' + \intsc(\Gamma_1, E; p_2) + \intsc(\Gamma_1, E; q_2).
\end{split}
\end{equation}

\medskip
\noindent \underline{Step 4. $\wt{\CR_g^0}\imply [R^s]\in [R^s]\in \wt{\CR_g}$}
\medskip

By using the similar argument to that used in Step 2, we can show that for a
$[(R^s)']\in \wt{\CR_g}$ lying in a neighborhood of $[R^s]\in
\wt{\CR_g^0}$, there is a point
\[
\big([(R^s)'], [\Upsilon_0'], w_1', w_2', ..., w_\beta', ...\big)\in
\wt{\Z_g}
\]
where $\Upsilon_0'$ contains two Type II chains $C_1$ and $C_2$ such that

\begin{enumerate}
\item $C_1$ and $C_2$ degenerate to $\overline{p_1q_1}$ and
$\overline{p_2q_2}$, respectively, when $(R^s)'$ degenerates to $R^s$;
\item
\begin{equation}\label{thm3-e6}
\begin{split}
\sum_{q\in C_j\cap E} \intsc(\Sigma_1', E; q) &= \intsc(\Sigma_1, E; p_j) +
\intsc(\Sigma_1, E; q_j)\\
&= 2(\mu_j + \mu_j') + \intsc(\Gamma_1, E; p_j) +
\intsc(\Gamma_1, E; q_j)
\end{split}
\end{equation}
for $j = 1, 2$, where $\Sigma_1'$ is the component of
$\Upsilon_0'$ defined as usual;
\item
\[
\{w_1', w_2', ..., w_\beta'\}\cap C_j = \{ w_i: i\in W_j \cup W_j'\}
\]
for $j = 1,2$.
\end{enumerate}

So applying \thmref{thm2} to $\Upsilon_0'\subset (R^s)'$, we have
\begin{equation}\label{thm3-e7}
\sum_{i\in W_j \cup W_j'} \delta_i \ge \sum_{q\in C_j\cap E} 
\intsc(\Sigma_1', E; q),
\end{equation}
for $j = 1, 2$.
Combining with \eqref{thm3-e5} and \eqref{thm3-e6}, we have
\begin{equation}\label{thm3-e8}
\sum_{i\in W_1} \delta_i + \sum_{i\in W_2} \delta_i 
\ge 2(\mu_1 + \mu_2) = 2\mu.
\end{equation}
If the equality in \eqref{thm3-e8} holds, the equalities in \eqref{thm3-e7}
have to hold. So by \thmref{thm2}, $\nu_i = 1$ for $i\in W_1\cup W_2\cup
W_1' \cup W_2'$. Combining with \eqref{thm3-e4}, we have proved
\thmref{thm3}.

\subsection{Completion of the Proof of \thmref{T:NOD}}

With \thmref{thm1}, \ref{thm2} and \ref{thm3} in place,
\thmref{T:NOD} is more or less obvious.

Let $C_0, C_1, ..., C_l\subset \Upsilon_0$ be all the $F$-chains in
$\Upsilon_0$ where $C_0$ is the Type I chain. It follows from
\thmref{thm1}, \ref{thm2} and \ref{thm3} that
\begin{align*}
\delta(\Upsilon_t) &= \sum_{i=0}^l \delta(\Upsilon_t, C_i)\\
&\ge \sum_{i=0}^l \sum_{q\in C_i\cap E} \intsc(\Sigma_1, E; q) - 1\\
&= \Sigma_1\cdot E - 1 = g.
\end{align*}

Obviously, $\delta(\Upsilon_t) = g$. Hence we must have
\[
\delta(\Upsilon_t, C_i) = \sum_{q\in C_i\cap E} \intsc(\Sigma_1, E; q)
\]
for $i = 1, 2, ..., l$. Consequently, 
by \thmref{thm2} and \ref{thm3}, $\Upsilon_t$ only
has nodes as singularities in the neighborhood of $C_i$ for
$i=1,2,...,l$. And by \thmref{thm1}, $\Upsilon_t$ only has nodes as
singularities in the neighborhood of $C_0$. Therefore, all the singularities
of $\Upsilon_t$ are nodes.

\end{document}